\documentclass[3p,times]{elsarticle}

\usepackage{mathtools}
\usepackage{amssymb}
\usepackage{amsfonts}

\usepackage{dsfont}
\usepackage{amsthm}
\usepackage{cases}
\usepackage{bm}
\usepackage{geometry}
\usepackage{indentfirst}
\usepackage{tabularx}
\usepackage{booktabs}
\usepackage{graphicx}
\usepackage{graphicx}
\usepackage{ulem}
\usepackage{xspace}
\usepackage{xcolor}
\usepackage{aliascnt}
\usepackage{multirow}
\usepackage{xspace}

\newcommand{\x}{\mathbf{x}} 
\newcommand{\Q}{\mathbf{Q}} 
\newcommand{\V}{\mathbf{V}}
\newcommand{\F}{\mathbf{F}}
\newcommand{\G}{\mathbf{G}}

\newcommand{\T}{\mathbf{T}}
\newcommand{\A}{\mathbf{A}}
\newcommand{\B}{\mathbf{B}}
\renewcommand{\H}{\mathbf{H}}
\newcommand{\jbeta}{j,\,\beta\,}
\newcommand\smallO{
  \mathchoice
    {{\scriptstyle\mathcal{O}}}% \displaystyle
    {{\scriptstyle\mathcal{O}}}% \textstyle
    {{\scriptscriptstyle\mathcal{O}}}% \scriptstyle
    {\scalebox{.7}{$\scriptscriptstyle\mathcal{O}$}}%\scriptscriptstyle
  }

\newcommand{\n}{\mathbf{n}}
\renewcommand{\v}{\mathbf{v}}
\newcommand{\vv}{\tilde{\mathbf{v}}}
\newcommand{\vel}{\pmb{\mathsf{v}} }
\renewcommand{\u}{\mathit{u}}
\newcommand{\q}{\mathbf{q}}
\newcommand{\dt}{\Delta t}

\newcommand{\Path}{\boldsymbol{\Psi}}

\usepackage[tracking=true]{microtype}

\newcommand{\Puno}{$P_1$\xspace}
\newcommand{\Pdue}{$P_2$\xspace}
\newcommand{\Ptre}{$P_3$\xspace}

\newcommand{\be}{\begin{equation} \begin{aligned} }
\newcommand{\ee}{\end{aligned} \end{equation}}

\usepackage{cancel}
\usepackage{soul}
\setstcolor{red}

\newtheorem{theorem}{Theorem}

\newtheorem{proposition}[theorem]{Proposition}

\newcommand{\wb}{\textcolor{black}}

\journal{Computer Methods in Applied Mechanics and Engineering}
\begin{document}
	
	%!=========================================================================
	%!
	%!      F R O N T    M A T T E R
	%!
	\begin{frontmatter}
		%-------------------------------------------------------
		% TITLE
		\title{Discontinuous Galerkin schemes for hyperbolic systems in non-conservative variables: quasi-conservative formulation with subcell finite volume corrections}
		%-------------------------------------------------------
		%-------------------------------------------------------
		% AUTHORS
		\author[univr]{Elena Gaburro$^*$}
		\ead{elena.gaburro@univr.it}
		\cortext[cor1]{Corresponding author}
		
		\author[usb]{Walter Boscheri}
		\ead{walter.boscheri@univ-smb.it}
		
		\author[koeln]{Simone Chiocchetti}
		\ead{simone.chiocchetti@uni-koeln.de}
				
		\author[inria]{Mario Ricchiuto}
		\ead{mario.ricchiuto@inria.fr}

		%-------------------------------------------------------
		% INSTITUTIONS		
		\address[univr]{Department of Computer Science, University of Verona, Strada le Grazie 15, Verona, 37134, Italy}		
		\address[usb]{Laboratoire de Mathématiques UMR 5127 CNRS, Universit{\'e} Savoie Mont Blanc, 73376 Le Bourget du Lac, France}
		\address[koeln]{Division of Mathematics, University of Cologne, Weyertal 86-90, Cologne, 50931, Germany}
		\address[inria]{Inria, Univ. Bordeaux, CNRS, Bordeaux INP, IMB, UMR 5251, 200 Avenue de la Vieille Tour, 33405 Talence cedex, France}
		%-------------------------------------------------------

		%-------------------------------------------------------
		% ABSTRACT 
		\begin{abstract}  % N.B. US English!
			We present a novel quasi-conservative arbitrary high order accurate ADER discontinuous Galerkin (DG) method 
			allowing to efficiently use a non-conservative form of the considered partial differential system, 
			so that the governing equations can be solved directly in the most physically relevant set of variables. 
			This is {particularly interesting} for multi-material flows with moving interfaces and steep, large magnitude contact discontinuities, 
			as well as in presence of highly non-linear thermodynamics.
			However, the non-conservative formulation of course introduces a conservation error which would normally lead to a wrong approximation of shock waves.			
			Hence, from the theoretical point of view, we give a formal definition of the conservation defect of non-conservative schemes and 
			we analyze this defect providing a local quasi-conservation condition, which allows us to prove a \textit{modified Lax-Wendroff theorem}. 
			Within this formalism, we also reformulate classical results concerning smooth solutions, contact discontinuities and moving interfaces.				
			Then, to deal with shock waves in practice, 
			we exploit the framework of the so-called \textit{a posteriori} subcell finite volume~(FV) limiter,
			so that, in troubled cells appropriately detected, we can incorporate a local conservation correction. 
			Our corrected FV update entirely removes the local conservation defect, allowing, at least formally, 
			to fit in the hypotheses of the proposed modified Lax-Wendroff theorem. 					
			Here, the shock-triggered troubled cells are detected by combining 
			physical admissibility criteria, a discrete maximum principle and a shock sensor inspired by Lagrangian hydrodynamics.			
			  
			To prove the capabilities of our novel approach, first, we show that we are able to 
			recover the same results given by conservative schemes on classical benchmarks for the single-fluid Euler equations. 
			We then conclude the presentation by showing the improved reliability of our 
			scheme on the \textit{multi-fluid Euler system} on examples like the interaction of a shock with a helium bubble for which we are 
			able to avoid the development of any spurious oscillations.
		\end{abstract}
		%-------------------------------------------------------
		
		%-------------------------------------------------------
		% KEY WORDS
		%-------------------------------------------------------
	\end{frontmatter}
	%!=========================================================================
	
%--------- SECTION --------------------------------------------------------
\section{Introduction} \label{sec_intro}

The benefits, as well as the dangers, concerning the use of non-conservative formulations
of hyperbolic conservation and balance laws are well understood since the nineties.
Indeed, {regrettably} the use of non-conservative formulations may lead to 
inadmissible results in presence of genuinely non-linear discontinuities 
(shocks)~\cite{karni1992viscous,lefloch94}, 
but conversely, these formulations are also known to be better suited to treat interfaces and in
 general contact discontinuities.
This issue had already been remarked in~\cite{zbMATH04083250} and it is especially relevant for multi-material/phase flows, 
but can also affect simulations of single phase flows with strong contact discontinuities.  
In fact, non-conservative formulations provide the advantage of evolving directly physical variables, 
and allow us to avoid evaluating complex thermodynamics across interfaces, 
which is a well known source of spurious oscillations, which may have catastrophic
consequences on the simulations~\cite{abgrall1996prevent,abgrall2003efficient,meng2024high}.
Reducing the number of evaluations of the thermodynamics also allows considerable reductions in 
computational costs for very complex systems, as discussed e.g in~\cite{dumbser2014high,ZanottiDumbser2016}.

The seminal work of Karni in~\cite{karni1992viscous} has shown that using the non-conservative formulation 
is indeed possible provided that the numerical viscosity of the numerical scheme 
is appropriately designed to match the one of a conservative discretization. 
The notion of well controlled dissipation was later generalized 
by Le Floch, Mishra and collaborators (see e.g.~\cite{beljadid_etal2017} and references therein). 
More recently, a general correction strategy of non-conservative schemes 
allowing instead to enforce given jump conditions has been proposed e.g. in~\cite{ABGRALL201810} and subsequent works.
Another general constructive methodology is provided by the path conservative approach~\cite{Castro2006,CASTRO2017131,RENAC20191}. 
Here, the obtained results naturally depend on the choice of the path which could allow, 
in principle, to also include consistent viscous profiles. 
But, for simplicity when dealing with complex systems, especially those not admitting a conservative form, 
the most of the time, the segment path is used: 
this choice is in general not consistent with the correct jump conditions, when available, 
and thus may provide wrong evolution of solutions containing shocks, 
unless some additional correction is included~\cite{ABGRALL20102759,CHALONS2017592,chalons2020path,pimentel2022cell}.
Nevertheless, the simplicity of the path-conservative formalism with the segment path allows it
to be used in the most disparate variety of cases, 
even for novel systems of partial differential equations for which such jump conditions are rather complex~\cite{carlino2023well,dumbser2024well,gaburro2021well,cao2023flux}
especially if the size of the intercell jumps is small.

Adaptive or hybrid approaches have also been proposed to ensure the local consistency with a given set of
the Rankine-Hugoniot relations.  
Well known works are for example those due to Karni~\cite{karni96}, and Abgrall and Karni~\cite{ABGRALL2001594}. 
The first proposes the idea of a local blending two updates of the pressure, 
a fully non-conservative one at interfaces, and the one based on the thermodynamics relations applied to the updated conserved variables.
The second paper presents a well known ghost-fluid approach in which
the single phase thermodynamics are independently applied on the two sides of a material interface, 
avoiding spurious pressure oscillations at the cost of a small mismatch in energy fluxes.
Then, in~\cite{abgrall_kumar_2014} a blending between Roe's method with Glimm's scheme
has been shown to be able to reproduce, even for systems with non-conservative products, 
consistent jump relations when available. 
Other approaches can be found in literature. 
Related to the above methods are the so called double flux, or multiple flux methods~\cite{BILLET20031473,LV2014105}  
which combine updates using different fluxes to satisfy both conservation in shocks and enhance the approximation of contacts.
In a similar spirit, the authors of~\cite{adaptive-FV} propose an effective local blending 
between a non-conservative method fully based on physical variables, 
with a conservative Roe-type scheme only used for shock waves.

%\begin{itemize}
%\item It is well known that using non-conservative  equations/methodes can lead to bad results  ~\cite{karni1992viscous} - ~\cite{lefloch94} 
%\item it also well known that non conservative formulations are better  suited to treat interfaces, and in general contact discontinuities.
%This  issue is especially relevant for multi-material/phase but also relevant in single phase with strong contact 
% discontinuities ~\cite{abgrall1996prevent} -~\cite{abgrall2003efficient}
%\item known features of discontinuous solutions: contacts vs shocks,  pressure/velocities vs density/entropy~\cite{karni1992viscous} 
%\item use of physical variables instead of conservative also often allows to avoid cumbersome inversion of
 % non-linear thermodynamic relastions~\cite{ZanottiDumbser2016}  
%\item  using non-conservative form based on appropriate design of viscosity~\cite{karni1992viscous} - ~\cite{beljadid_etal2017}
%\item  using non-conservative form based on hybrid/adaptive schemes or of corrections allowing the correct discrete jump relations
%~\cite{karni96}- ~\cite{ABGRALL201810}
%  -~\cite{Tang-Zhou99} %-~\cite{adaptive-FV} %  -~\cite{abgrall_kumar_2014}	
%\end{itemize}	

\medskip

The work proposed here can be considered a hybrid approach, 
and falls into the framework of high order discontinuous Galerkin (DG) methods. 
The philosophy behind our strategy exploits the fact {that conservation 
is only expressed} in terms of the evolution of cell averages. 
This suggests that a non-conservative formulation can be safely used everywhere for the higher order modes, 
as long as a consistent conservative update is guaranteed for the averages. 
Next, we exploit the fact that exact conservation is strictly needed only in some areas (shocks), 
and a fully non-conservative update can be used everywhere else. 
We thus propose a numerical scheme based on a fully non-conservative DG method 
which we use to evolve the high degree polynomials discretizing 
the selected non-conservative variables (chosen according to physical and/or computational convenience).
Then, the first idea we exploit consists in making the method locally conservative
by means of a correction of the constant mode of the primitive variables, 
consistent with a finite volume (FV) update of the cell averages of the conserved variables.
The second original idea is to insert this approach within an \textit{a posteriori} limiter,
applied to the cells in which a shock is found, which we flag by  
using a divergence-based shock sensor, similar to those proposed and studied e.g. in 
\cite{DUCROS1999517,FUJIMOTO2019264}. 
This allows to keep a fully non-conservative update away from shocks, 
while guaranteeing {some form} of local conservation where necessary.

These two steps fit beautifully in the framework of an ADER-DG time stepping method with local space-time predictors,
and finite volume subcell limiters~\cite{Dumbser2008,DGLimiter2,Gaburro2021PNPMLimiter}.
The locality of the predictors naturally lends itself to exploit a non-conservative formulation of the PDEs, 
as already suggested in~\cite{ZanottiDumbser2016}.
In addition, we extend the previous work by 
i) computing also the corrector step on the non-conservative variables, 
so to retain the overall consistency of the discretization 
and to fully exploit the advantages of the non-conservative form at the interfaces, 
{and ii) with a hybrid formulation in which} the consistency is enforced thanks to a conservative update 
applied on shock-flagged elements within the \textit{a posteriori} subcell finite volume limiter.  
The resulting discretization, while evolving polynomials discretizing the most convenient set of physical variables, 
embeds a fully conservative evolution of subcell averages of the conserved variables in correspondence of shocks.

To give further justification to our approach, and in particular to the  
\textit{a posteriori} modification of the cell-average update, 
we propose a notion of quasi-conservation based on a vanishing defect, 
which leads to a modified formulation of the Lax-Wendroff theorem. 
This characterization is described specifically for our ADER-DG setting, 
to which some known results~\cite{abgrall1996prevent,abgrall2003efficient} are adapted. 
 
The obtained scheme is tested thoroughly on benchmarks involving smooth solutions 
as well as strong shocks, Riemann problems, and shock-contact interactions.
All the numerical results confirm the high accuracy of the proposed approach for smooth solutions, 
a fully conservative and non-oscillatory approximation of very non-linear discontinuities, 
and the capability of computing strong contact discontinuities free of any spurious numerical artifacts.

%  
%\begin{itemize}
%\item Philosophy: conservation is expressed by the evolution of cell averages, while a lot of freedom exists for the other modes
%\item  Underlying method: evolution of  high degree polynomials  constructed from values of   non-conservative variables selected on physical principles
%\item  First original idea: introduce an a-posteriori  correction  of the constant mode to ensure conservation
%\item In practice: correction using conservative finite volume based on space time primitive variable expansion anchored to previous time step conservative average
%\item Second original idea: subcell formulation of the  a-posteriori correction allows combination with a-posteriori limiter leading to a 
% local correction only in troubled cells
%\item In practice: shock filter based on divergence sensor (similar to the Ducros et al sensor) allows a local correction only in 
% shocks~\cite{DUCROS1999517} -~\cite{FUJIMOTO2019264} 
%\item Resulting scheme is conservative across shocks, while allowing use of physically selected non-conservative form everywhere else
%\item 
%\end{itemize}	

\medskip

The rest of the paper is organized as follows. 
We first introduce the structure of the hyperbolic partial differential equations~(PDE) considered in this work 
and we present both their conservative and non-conservative formulation, see Section~\ref{sec_pde}.  
As application, we consider in this section 
a two-species formulation of the Euler equations for a perfect gas.
Next, in Section~\ref{sec_numscheme}, 
we describe the non-conservative ADER-DG framework with \textit{a posteriori} subcell limiter. 
Then, Section~\ref{sec_conservation} is entirely devoted to our novel approach. 
We start by proposing a discussion on the issue of conservation errors, 
and we introduce a Lax-Wendroff theorem with vanishing conservation defect. 
Afterwards, {we recast in our setting some known examples allowing a non-conservative treatment,
namely} smooth flows, and contact discontinuities for single and multi-species flow. 
And furthermore, we present a strategy to embed a local conservation correction in the \textit{a posteriori} limiter, 
in such a way to remove the conservation defect in appropriate cells.  
In Section~\ref{sec_validation}, we show a large set of numerical results 
covering benchmarks coming from the Euler equations of gasdynamics 
and challenging simulations modeled by the multi-material Euler equations. 
We demonstrate that our scheme achieves up to fourth order of accuracy in space and time, 
correctly captures shock waves and avoid spurious oscillations at the (smeared) interface between different gases. 
The numerical results are commented and compared with available reference solutions wherever possible.
Finally, we close the paper with some remarks and an outlook to future works in Section~\ref{sec_concl}.
	
%--------- END OF SECTION -------------------------------------------------	

%--------- SECTION --------------------------------------------------------
\section{Governing equations} \label{sec_pde}	

In this paper we consider the numerical approximate solutions of hyperbolic systems of conservation laws
\begin{equation}
	\label{eqn.PDE}
	\partial_t \Q  + \nabla \cdot \mathbb{F}(\Q) = \mathbf{0}, \qquad \x \in \Omega \subset \mathds{R}^d, 
	\quad t \in \mathds{R}_0^+, \quad \Q \in \Omega_{\Q} \subset \mathds{R}^\nu,     
\end{equation}
where $d=2$ is the number of space dimensions, $\x=(x,y)$ is the vector of 
spatial coordinates in the physical domain $\Omega$ and $t$ denotes the time. 
The vector of conserved variables $\Q$ is defined in the space of the admissible 
states $\Omega_{\Q}\subset \mathds{R}^\nu$, 
while $\mathbb{F}(\Q)=\left( \F(\Q),\G(\Q) \right)$ represents the non-linear flux tensor.
The governing PDE~\eqref{eqn.PDE} can also be formulated in terms of a set of  physical (or primitive) 
variables $\V$, which belong to the admissible states $\Omega_{\V} \subset \mathds{R}^\nu$. 
The Jacobian  matrix of the transformation $\V(\Q)$ is denoted by   $\T$, 
while with $\A$ and $\B$ we denote   the Jacobian matrices of the conservative fluxes
\begin{equation}
	\label{eqn.Jac}
	\T=\frac{\partial \V}{\partial \Q}, \qquad \A=\frac{\partial \F}{\partial \Q}, \qquad \B=\frac{\partial \G}{\partial \Q}.
\end{equation} 
Multipling  from the left  the conservative form~\eqref{eqn.PDE} by   $\T$ and introducing the Jacobian matrices defined above, we obtain
\begin{eqnarray}
	\label{eq.passagescons2prim}
	&&\T \, \frac{\partial \Q}{\partial t} + \T \, \left( \A \frac{\partial \Q}{\partial x} + \B \frac{\partial \Q}{\partial y} \right) = \mathbf{0} \nonumber \\
	&&\frac{\partial \V}{\partial t} + \T \, \A \, \T^{-1} \, \frac{\partial \V}{\partial x} + \T \, \B \T^{-1} \, \frac{\partial \V}{\partial y} = \mathbf{0}.
\end{eqnarray}
A non-conservative form of   system~\eqref{eqn.PDE} using primitive variables therefore reads  
\begin{equation}
	\label{eqn.PDE_prim}
\frac{\partial \V}{\partial t} + \A_{\V} \, \frac{\partial \V}{\partial x} + \B_{\V} \, \frac{\partial \V}{\partial y} = \mathbf{0},
\end{equation}
with the definitions
\begin{equation}
	\label{eqn.Jac_prim}
	\A_{\V}:=\T \, \A \, \T^{-1}, \qquad \B_{\V}:=\T \, \B \, \T^{-1}.
\end{equation}
Let us remark that the  matrices $(\A_{\V},\B_{\V})$ involved in the non-conservative form  are similar 
by construction to the Jacobian matrices of the conservative system $(\A,\B)$, hence they share the 
same eigenvalues. 
Also, we assume that the mapping $\Q(\V)$  is always invertible, and  verifies the relation
\begin{equation}\label{eqn.constoprim-diff}
\Q(\V_2)-\Q(\V_1)=\T^{-1}(\V_1,\V_2)(\V_2 - \V_1),
\end{equation}
for some appropriately defined (Roe-like) average of $\T^{-1}$\wb{, i.e. $\T^{-1}(\V_1,\V_2)$}. Explicit examples will be provided later.
Also, we assume that   all the above matrices are continuous \wb{with respect to} their arguments, namely that
for any relevant choice of the  norms involved
\begin{equation}
	\label{eqn.Jac_prim_cont}
\| \T(\Q) - \T(\tilde \Q)\|  \le  K_{\T} \|\Q- \tilde \Q\|\;,\;\;	
\| \T^{-1}(\V) - \T^{-1}(\tilde \V)\|  \le  K_{\T^{-1}} \|\V- \tilde \V\|\;,\;\;	
%\|\A(\Q)-\A(\tilde \Q)\|  \le  K_{\A} \|\Q- \tilde \Q\|\;,\;\;
\|\A_{\V}(\V)-\A_{\V}(\tilde \V)\|  \le  K_{\A_{\V}} \|\V- \tilde \V\|\;,
\end{equation}
with bounded constants $K_{\T}$, $K_{\A}$, $K_{\A_{\V}}$.

\smallskip

As discussed in the introduction, due to the non-conservative nature of~\eqref{eqn.PDE_prim}  particular 
care must be devoted to the construction of shock capturing numerical schemes for its  discretization.

\subsection{Multi-material compressible Euler equations}\label{sec_pde_Euler}

As motivational example for the scheme proposed in this paper, 
{as well as for all the presented numerical applications}, 
we consider a multi-material formulation of the compressible Euler equations, 
modeling the flow of two immiscible perfect gases,
following e.g.~\cite{quirk1996dynamics,karni96,abgrall1996prevent,ABGRALL2001594,abgrall2003efficient,gouasmi2020formulation,li2022high}. 
The system can be written in the conservative form~\eqref{eqn.PDE} with 
\begin{equation}
	\label{eq.multimaterialCons}
	\Q = \left( \begin{array}{c} \rho   \\ \rho u  \\ \rho v  \\ \rho E \\ \rho \chi  \end{array} \right), \quad
	\mathbb{F} = \left( \begin{array}{ccc}  
		\rho u       & \rho v        \\ 
		\rho u^2 + p & \rho u v          \\
		\rho u v     & \rho v^2 + p      \\ 
		u(\rho E + p) & v(\rho E + p)   \\
		u \rho \chi & v\rho\chi
%		\tilde{\gamma}      & v \rho \tilde{\gamma}
	\end{array} \right),
\end{equation}
where the conserved variables $\Q$ involve the fluid density $\rho$, 
the momentum density vector $\rho \vel=(\rho u, \rho v)$ and the total energy density $\rho E$. 
We have also introduced a marker function $\chi(x,y,t)$ which is advected with the material interface.
One can then evaluate the equation of state based on the knowledge of $\chi$. 
This variable is assumed to take discrete values corresponding to the different materials 
present in the {flow}. In this paper we will only consider the case of two materials, 
but everything can be generalized to more components. 
So for this paper we take $\chi\in \{\chi_1,\chi_2\}$ and its numerical
approximation will lay in the continuous interval $[\chi_1,\chi_2]$. 
{The variable $\tilde{\gamma}$ is the (material-dependent) adiabatic index $\tilde{\gamma}= \gamma\chi$}, 
which can vary in space if there are interfaces between different gases, 
with $\gamma=1.4$ the value taken for air in standard atmosphere.  
More precisely, the fluid pressure $p$ is related to the conserved quantities via the ideal gas equation of state (EOS)  
\begin{equation}
	\label{eqn.eos} 
	 p = (\tilde{\gamma}-1)\left(\rho E -  \rho \kappa \right),
\end{equation}
where $\kappa=u^2/2 + v^2/2$ is the kinetic energy.
The speed of sound associated to the system is $c=\sqrt{\tilde\gamma p/\rho}$.
{Let 
\begin{equation}\label{eqn.lambda}
\tilde \lambda =  \dfrac{1 }{\tilde{\gamma}-1},
\end{equation}
a classical choice, related to the results of~\cite{abgrall1996prevent,abgrall2003efficient}, is to set  for example  
\begin{equation}
	\label{eqn.eos1a} 
	\tilde \lambda(\chi)  =  \dfrac{ \chi  }{\gamma -1} 
\end{equation}
in which case the values $\chi_1$ and $\chi_2$ represent the ratios $\chi_j=(\gamma-1)/(\gamma_j-1)$.
In practice, one obtains excellent results also using the simpler  definition 
\begin{equation}
	\label{eqn.eos1} 
		\tilde \lambda(\chi) =  \dfrac{1}{\tilde\gamma(\chi)-1}\;,\quad
 \tilde{\gamma}(\chi) =\chi \gamma 
\end{equation}
in which case $\chi_j=  \gamma_j/\gamma$. 
The  Euler equations for a single perfect gas are recovered for constant $\chi$, and/or by neglecting the last equation.
}

Moreover, the system \eqref{eq.multimaterialCons} can be re-written in the non-conservative 
form~\eqref{eqn.PDE_prim} by considering for example the vector of physical variables
\begin{equation}
		\label{eq.multimaterialPrim0}
		\V = \left( \begin{array}{l} \rho  \;\; u  \;\; v  \;\; p \;\; \chi  \end{array} \right),
\end{equation}			
for which the matrices in~\eqref{eqn.PDE_prim} become
\begin{equation}
		\label{eq.multimaterialPrim}
		\A_{\V} = \left( \begin{array}{ccccc}  	
			u & \rho & 0 & 0 & 0 \\
			0 & u & 0 & 1/\rho & 0 \\
			0 & 0 & u & 0  & 0 \\
			0 & \tilde{\gamma} p & 0 & u & 0 	\\
			0 & 0 & 0 & 0 & u  \\
		\end{array} \right), \quad 
		\B_{\V} =  \left( \begin{array}{ccccc} 
			v & 0 & \rho & 0 & 0\\
			0 & v & 0 & 0 & 0 \\
			0 & 0 & v & 1/\rho & 0  \\
			0 & 0 & \tilde{\gamma} p & v  & 0 \\
			0 & 0 & 0 & 0 & v  \\
		\end{array} \right).
\end{equation}		
%Finally, we remark that for all the test cases presented in this Section, for what concerns the shock indicator 
%of Equation~\eqref{eqn.flattener}, we have set $m_1=15$. 	
The transformation from conservative to primitive variables $\T$ depends on the computation of the equation of state. 
We can write the transformation in a general form as 
\begin{equation}
		\label{eq.constoprim}
		\T^{-1}(\V) = \left( \begin{array}{ccccc}  	
			1 & 0 & 0 & 0 & 0 \\
			u & \rho & 0 & 0 & 0 \\
			v & 0 & \rho & 0  & 0 \\
			\kappa & \rho u & \rho v & \tilde \lambda(\chi) & \beta(p,\chi,\gamma)  	\\
			\chi & 0 & 0 & 0 & \rho  \\
		\end{array} \right),
\end{equation}	
having denoted by $\kappa=u^2/2 + v^2/2$ the kinetic energy.
When the formula \eqref{eqn.eos1a} is used, {one can easily check that
\begin{equation}
		\label{eq.constoprima}
		 \beta(p,\chi,\gamma)  = \dfrac{p}{\gamma-1}\,,
\end{equation}	
while in the case of   \eqref{eqn.eos1} we have a slightly more complex formula 
\begin{equation}
		\label{eq.constoprimb}
		 \beta(p,\chi,\gamma) = -\gamma p \tilde\lambda(\chi)^2 = - \gamma\dfrac{p}{(\tilde\gamma(\chi)-1)^2}  \,.
\end{equation}	
}
%\begin{equation}
%		\label{eq.constoprim}
%		\T^{-1}(\V) = \left( \begin{array}{ccccc}  	
%			1 & 0 & 0 & 0 & 0 \\
%			u & \rho & 0 & 0 & 0 \\
%			v & 0 & \rho & 0  & 0 \\
%			\kappa & \rho u & \rho v & \dfrac{1}{\tilde\gamma -1} &  -\gamma\dfrac{p}{(\tilde\gamma-1)^2} 	\\
%			\tilde\gamma & 0 & 0 & 0 & \rho  \\
%		\end{array} \right),
%\end{equation}	
For this system it is easy to verify that \eqref{eqn.constoprim-diff} holds with 
	\begin{equation}
		\label{eqn.constoprim-diff2}
	\widehat \T^{-1}(\V_1,\V_2)= \left( \begin{array}{ccccc}  	
			1 & 0 & 0 & 0 & 0 \\
			\overline u_{12} &\overline \rho_{12}& 0 & 0 & 0 \\
			\overline v_{12} & 0 &\overline \rho_{12} & 0  & 0 \\
			\overline \kappa_{12}&\overline \rho_{12}\,\overline u_{12}  & \overline \rho_{12}\, 
			\overline v_{12}&   \dfrac{\overline{ \tilde\gamma}_{12} -1}{(\tilde\gamma_1 -1)(\tilde\gamma_2 -1)} &    \bar \beta_{12} 	\\
			\overline{ \tilde\gamma}_{12} & 0 & 0 & 0 &\overline  \rho_{12}  \\
		\end{array} \right),
\end{equation}	
where 
$$
\overline{(\cdot)}_{12}=\dfrac{(\cdot)_1+(\cdot)_2}{2}
$$
denotes the classical arithmetic average. The averaged thermodynamic   coefficient $\beta$
is 
$$
 \bar \beta_{12} = \dfrac{\overline p_{12}}{\gamma-1 }
$$
for  \eqref{eqn.eos1a}, and 
$$
 \bar \beta_{12} = -\dfrac{\gamma\overline p_{12}}{(\tilde\gamma_1 -1)(\tilde\gamma_2 -1)}
$$
when using \eqref{eqn.eos1}.

%--------- END OF SECTION -------------------------------------------------	

%--------- SECTION --------------------------------------------------------
\section{ADER-DG framework for hyperbolic systems in non-conservative form} \label{sec_numscheme}	

To discretize the general hyperbolic system written in the non-conservative form~\eqref{eqn.PDE_prim}  
we consider a framework combining a high order discontinuous Galerkin approximation,
with the fully discrete ADER (Arbitrary high order DErivative Riemann problem) predictor-corrector
approach, originally due to~\cite{toro1,toro3,toro4} and proposed in its modern form in~\cite{Dumbser2008}. 
Moreover, to integrate the non-conservative form we follow the so called \textit{path-conservative} formulation 
{to define fluctuations at faces}~\cite{Pares2006,dumbser2011simple}.

\subsection{Geometry and notation} 

Let us consider a discretization of the \textit{spatial} computational domain $\Omega$ by means of 
a tessellation composed of non-overlapping control volumes of arbitrary polygonal shape.
We denote by $\omega_i$ the generic mesh element, and by $\partial \omega_i$ its boundary.
The perimeter and surface of an element are instead denoted by $|\partial \omega_i|$ and
$|\omega_i|$ respectively~\cite{Gaburro2020Arepo,boscheriAFE2022}. 
Given two adjacent cells $\omega_i$ and $\omega_j$,
each edge $\partial \omega_{ij}$ shared between the cells  is oriented by $\n_{ij}$, 
the  unit normal vector  pointing outwards with respect to  $\omega_i$. 
In each cell we define the   barycentric coordinates 
\begin{equation}
\label{eq.barycenters}
\x_{{b}_i} = \frac{1}{N_{k,i}} \, \sum \limits_{k=1}^{N_{k,i}} \x_{k,i}, \qquad \forall \omega_i \in \Omega,
\end{equation}
with $N_{k,i}$  the number of vertices of  $\omega_i$, of coordinates $\x_{k,i}$. The barycenter
is used to define  the characteristic length $h_i$ of each cell, obtained as the   radius of gyration
\begin{equation}
	h_i %= \sqrt{\frac{J_{xx} + J_{yy}}{|\omega_i|}} 
	= \frac{1}{\sqrt{|\omega_i|}}\,{\left(\int_{\omega_i} \left((x - x_{{b}_i})^2 + (y - y_{{b}_i})^2 \right)\,dx\,dy\right)}^{1/2}.
\end{equation}
Note that all  the integrals over the cells are evaluated using Gaussian quadrature formulae of suitable accuracy on 
triangles~\cite{stroud}. This means that the polygonal cell $\omega_i$ is split into a total number 
of $N_{k,i}$ sub-triangles that are obtained by connecting the barycenter with each vertex. In particular,
the surface of a cell $|\omega_i|$ is trivially obtained as the sum of triangles' surfaces.  
Integrals along the segments forming cell  boundaries $\partial \omega_i$ are computed using 
one-dimensional Gauss quadrature formulae of appropriate degree.

The \textit{time} domain $[0, t_f]$ is subdivided in {slabs} $[t^n,\,t^{n+1}]$, with % I like slabs, it is very correct and precise
$t^n$ the level at which the solution is known and $t^{n+1}$ the temporal level of the unknown solution.
The time step is denoted by $\dt=t^{n+1}-t^n$, with no super/sub script referring to the specific level to simplify the notation.
The computation of $\dt$ is done according to the usual CFL-type stability condition for explicit DG schemes 
\begin{equation}
\label{eq.CFL}
\Delta t < \textnormal{CFL} \left( 
\frac{ |\omega_i| }{ (2N+1) \, d \, |\lambda_{\max,i}| \, |\partial\omega_i| } 
\right), 
\qquad \textnormal{CFL}\leq 1, 
\qquad \forall \omega_i \in \Omega,  
\end{equation}
where $N$ is the polynomial degree used for the approximation of the numerical solution 
and $|\lambda_{\max,i}|$ is the spectral radius of the Jacobian matrices of the system~\eqref{eqn.Jac_prim},
evaluated {using the average solution in cell $i$ and the point values of the polynomial at cell corners}.
\smallskip 

Within each $\omega_i$ we consider $\v_h^n$, the \textit{numerical solution} at the current time level $t^n$, 
represented by piecewise polynomials of degree $ N\geq 0$ by a modal expansion directly defined in physical space
\begin{equation}
\mathbf{v}_h^n:=\mathbf{v}_h(\x,t^n) = \sum \limits_{\ell=0}^{\mathcal{N}-1} \phi_\ell(\x) \,  \hat{\mathbf{v}}^{n}_{\ell}, \qquad \x \in \omega_i,
\label{eqn.uh}
\end{equation}
where $\mathcal{N}=\textnormal{dof}\left(N,d\right)$  is the total number of expansion coefficients $\hat{\mathbf{v}}^{n}_{\ell}$, 
usually referred to  as degrees of freedom. This number  is determined according to the formula
\begin{equation}
	\textnormal{dof}\left(N,d\right) = \frac{1}{d!} \prod \limits_{m=1}^{d} \left(N+m\right).
\end{equation}
The $\mathcal{N}$ basis functions $\phi_\ell(\x)$ in $\omega_i$ are explicitly given by
\begin{equation}
\phi_\ell(\x) = \left( \frac{x-x_{{b}_i}}{h_i} \right )^{p_{\ell}} \, \left( \frac{y-y_{{b}_i}}{h_i}\right)^{q_\ell}, 
\qquad \ell=0, \ldots, \mathcal{N}-1, \qquad 0 \leq p_{\ell}+q_{\ell} \leq N,
\label{eq.space_basis}
\end{equation}
with $\ell=(p_\ell,q_\ell)$ a multi-index notation for each function. 
%%$h_i$ the characteristic length of each $\omega_i$ simply defined as
%its radius of gyration
%\begin{equation}
%	h_i = \sqrt{\frac{J_{xx} + J_{yy}}{|\omega_i|}} = \frac{1}{\sqrt{|\omega_i|}}\,{\left(\int_{\omega_i} \left((x - x_{{b}_i})^2 + 
% (y - y_{{b}_i})^2 \right)\,dx\,dy\right)}^{1/2}.
%\end{equation}
%The total number 
%

Finally, we introduce the   predictor solution   $\vv_h$  consisting of  a local  space-time polynomial 
approximation  of the solution within the {slab} $\tilde{\omega_i}= \omega_i \times [t^n;t^{n+1}]$:
%
%Similarly, let us also introduce the discretization of  within each space-time cell:
\begin{equation}
\vv_h:=\vv_h(\x,t) = \sum \limits_{\ell=0}^{\mathcal{L}-1} \theta_\ell(\x,t) \,  \hat{\vv}_{\ell}, \qquad (\x,t) \in \tilde{\omega_i},
\label{eqn.qh}
\end{equation}
with a total number $\mathcal{L}=\textnormal{dof}\left(N+1,d\right)$ of degrees of 
freedom $\hat{\vv}^{n}_{\ell}$. The $\mathcal{L}$ basis functions $\theta_\ell(\x,t)$ indexed by the multi-index $\ell =(p_\ell,q_\ell,r_\ell)$ 
are given again by a rescaled Taylor expansion in space-time around the cell barycenter $\tilde{\x}_{{b}_i}=({\x}_{{b}_i},t^n)$, that is
\begin{equation}
\theta_\ell(\x) = \left( \frac{x-x_{{b}_i}}{h_i} \right )^{p_\ell} \, 
\left(\frac{y-y_{{b}_i}}{h_i}\right)^{q_\ell} \, 
\left(\frac{t-t^n}{h_i}\right)^{r_\ell}, 
\qquad \ell=0,\ldots,\mathcal{L}-1, \qquad 0 \leq p_{\ell}+q_{\ell}+r_{\ell} \leq N.
\label{eq.spacetime_basis}
\end{equation}

\paragraph{Remark} The degrees of freedom and the basis functions are element-dependent, thus for cell $\omega_i$ one should 
write $\phi_{\ell,i}$ and $\hat{\mathbf{v}}^{n}_{\ell,i}$, and similarly for   the space-time 
bases which refer to the space-time cell $\tilde{\omega_i}$ .
% and the associated expansion 
%coefficients we have $\theta_{\ell,i}$ and $\hat{\mathbf{q}}^{n}_{\ell,i}$. 
To simplify the notation the element subscript has been  omitted. 
%To lighten the notation we drop the subscript referring  to the element, bearing in mind that the cell dependency is implied. 

\subsection{Local space-time ADER predictor in physical variables}
\label{ssec.predictor}
In order to obtain a space-time representation of the predictor solution of the form~\eqref{eqn.qh}, 
a weak formulation of the PDE~\eqref{eqn.PDE_prim} is derived upon multiplication with a space-time 
test function $\theta_k(\x,t)$, which is of the same form of the basis 
functions $\theta_{\ell}$~\eqref{eq.spacetime_basis}, and subsequent integration over the space-time cell $\tilde{\omega}_i$:
\begin{equation}
\int \limits_{\tilde{\omega}_i} \theta_k \, \frac{\partial \V}{\partial t} \, d\x \, dt \ + \int 
\limits_{t^n}^{t^{n+1}} \int \limits_{\tilde{\omega}_i} \theta_k \, \left( \A_{\V} \, \frac{\partial 
\V}{\partial x} + \B_{\V} \, \frac{\partial \V}{\partial y} \right) \, d\x \, dt \ =\  \mathbf{0}.
\label{eqn.PDEweak}
\end{equation}
Let us now introduce the abbreviation
\begin{equation}
\H(\x,t) = \A_{\V} \, \frac{\partial \V}{\partial x} + \B_{\V} \, \frac{\partial \V}{\partial y}
\label{eqn.Hnoncons}
\end{equation}
and for its  space-time  approximation
\begin{equation}
\tilde{\H}_h:= \tilde{\H}_h(\x,t) = \sum \limits_{\ell=0}^{\mathcal{L}-1} \theta_\ell(\x,t) \,  
\hat{\tilde{\mathbf{H}}}_{\ell}, \qquad (\x,t) \in \tilde{\omega_i}.
\label{eqn.Hh}
\end{equation}
By inserting the ansatz~\eqref{eqn.qh} for the numerical solution $\V$ and~\eqref{eqn.Hh} for the 
non-conservative terms $\H(\x,t)$ and integrating by parts the first term on the left hand side 
of~\eqref{eqn.PDEweak}, we obtain the following element-local non-linear algebraic equation system 
for the  array $\vv$ containing
the coefficients $\hat{\vv}_{\ell}$ of the unknown space-time expansion
\begin{equation}
\mathbf{K}_1 \, \hat{\vv}  = \mathbf{F}_0 \, \hat{\mathbf{v}}^{n}   - \mathbf{M} \, \hat{\tilde{\mathbf{H}}}.
\label{eqn.ADERpred}
\end{equation}
Here above we have employed the following definitions
\begin{eqnarray*}
%\mathbf{K}_1 &=& \int \limits_{{\omega}_i} \theta_k(\x,t^{n+1}) \, \theta_{\ell}(\x,t^{n+1}) \, d\x \ -
 % \int \limits_{\tilde{\omega}_i}  \frac{\partial \theta_{k}(x,t)}{\partial t} \theta_{\ell}(\x,t) \, d\x \, dt, \nonumber \\
%\mathbf{F}_0 &=& \int \limits_{{\omega}_i} \theta_k(\x,t^{n}) \, \theta_{\ell}(\x,t^{n}) \, d\x, \\
%\mathbf{M} &=& \int \limits_{\tilde{\omega}_i} \theta_k(\x,t) \, \theta_{\ell}(\x,t) \, d\x \, dt \nonumber.
\mathbf{K}_1 = \int \limits_{{\omega}_i} \theta_k(\x,t^{n+1}) \, \theta_{\ell}(\x,t^{n+1}) \, d\x\! &-&\!\! 
\int \limits_{\tilde{\omega}_i}  \frac{\partial \theta_{k}(x,t)}{\partial t} \theta_{\ell}(\x,t) \, d\x \, dt, \nonumber \\
\mathbf{F}_0 = \int \limits_{{\omega}_i} \theta_k(\x,t^{n}) \, \theta_{\ell}(\x,t^{n}) \, d\x,\,\;\;&\,&\;\;
\mathbf{M} = \int \limits_{\tilde{\omega}_i} \theta_k(\x,t) \, \theta_{\ell}(\x,t) \, d\x \, dt, 
\end{eqnarray*}
with $\hat{\tilde{\mathbf{H}}}$ the array of the non-conservative product coefficients $\hat{\tilde{\mathbf{H}}}_{\ell}$
and $\hat{\mathbf{v}}^{n}$ the arrays of the degrees of freedom at time $t^n$ $\hat{\mathbf{v}}^{n}_{\ell}$, 
which allow to account for the initial condition of the local Cauchy problem~\eqref{eqn.PDEweak}.
% is taken into account by the contribution of matrix $\mathbf{F}_0$ with the known degrees of freedom at
 % time $t^n$, namely $\hat{\mathbf{v}}^{n}_{\ell}$. 
The non-linear system~\eqref{eqn.ADERpred} is then solved by means of a fixed point iterative method
%, namely
%\begin{equation}
%\hat{\mathbf{q}}^{n,s+1}_{\ell} = \mathbf{K}_1^{-1} \left( \mathbf{F}_0 \, \hat{\mathbf{v}}^{n}_{\ell} -
 % \mathbf{M} \, \hat{\mathbf{H}}^{n,s}_{\ell} \right),
%\label{eqn.nlsys_ader}
%\end{equation}
%where $s$ denotes the iteration number. The iterative scheme stops when the residual of the numerical solution is less than a 
with the  prescribed tolerance between two iterates $\vv_h^{n,s+1}$ and $\vv_h^{n,s}$
%which in practice is implemented as %$|\hat{\mathbf{q}}^{n,s+1}_{\ell}-\hat{\mathbf{q}}^{n,s}_{\ell}| < \delta$, with $\delta=10^{-12}$.
\be 
 {\left|\left| \vv_h^{n,s+1} - \vv_h^{n,s} \right|\right|_{L_2}}  < \delta\, {\left|\left| 
 \vv_h^{n,s} \right|\right|_{L_2}}, \quad \text{with} \quad \delta=10^{-7}.
\ee

\subsection{{Path-integrated} discontinuous Galerkin corrector} 
\label{ssec.corrector}
The space-time predictor $\vv_h$ is used to devise a fully discrete one-step DG scheme. 
The governing equations~\eqref{eqn.PDE_prim} are multiplied by a spatial test function $\phi_{k}(\x)$, which
is of the same form of the basis functions $\phi_{\ell}$~\eqref{eq.space_basis}, and then integrated in space and time. %:
%\begin{equation}
%\int \limits_{t^n}^{t^{n+1}} \int \limits_{\omega_i} \phi_{k} \frac{\partial \V}{\partial t} \, d\x \, dt \ +
 % \int \limits_{t^n}^{t^{n+1}} \int \limits_{\omega_i} \phi_{k} \, \left(\A_{\V} \, \frac{\partial \V}{\partial x} + \B_{\V} \, 
 % \frac{\partial \V}{\partial y} \right) \, d\x \, dt \ = \mathbf{0}.
%\end{equation}
The formulation in primitive variables  only involves non-conservative products 
as shown in Section~\ref{sec_pde}. These terms are easily  integrated
on the cell volumes, while a {path-integral} formulation is used to obtain face fluctuations, as
in~\cite{Toumi1992,Pares2006,Castro2006,Castro2008}. We thus obtain 
\begin{equation}
\left( \int \limits_{\omega_i} \phi_{k} \phi_{\ell} \, d\x \right) \, \left( \hat{\mathbf{v}}_{\ell}^{n+1} - 
\hat{\mathbf{v}}_{\ell}^{n} \right) \ + \int \limits_{t^n}^{t^{n+1}} \int \limits_{\partial \omega_i} 
\phi_{k} \mathbf{D} \cdot \n \, dS \, dt + \int \limits_{t^n}^{t^{n+1}} \int \limits_{ \omega_i 
\backslash \partial \omega_i}\phi_{k}  \H(\vv_h) \, d\x \, dt = \mathbf{0},
\label{eqn.DGscheme}
\end{equation}
where we use the approximation of the numerical solution~\eqref{eqn.uh} and the abbreviation for the 
non-conservative terms~\eqref{eqn.Hnoncons}. 
The  term $\mathbf{D}$ takes into account {solution jumps}   
%jumps of the solution %$\V$ 
on the element boundaries $\partial \omega_i$ according to the path integral
\begin{eqnarray}
%\mathbf{D} \cdot  \mathbf{\tilde n}   &=& \half \, \int \limits_0^1 
%\left[ \A_{\V}\left(\Path(\q_h^-,\q_h^+,s)\right), \ \B_{\V}\left(\Path(\q_h^-,\q_h^+,s)\right) \right] 
% \cdot \n \, \frac{\partial \Path}{\partial s} \, ds \nonumber \\
%&=& \half \ \left( \int \limits_0^1 
%\left[ \A_{\V}\left(\Path(\q_h^-,\q_h^+,s)\right), \ \B_{\V}\left(\Path(\q_h^-,\q_h^+,s)\right) \right] 
% \cdot \n \, ds \right) \, \left( \q_h^+ - \q_h^- \right),
\mathbf{D} \cdot  \mathbf{ n}  &=& % \half \,
 \int \limits_0^1 
\left\{\left[ \A_{\V}\left(\Path(\vv_h^-,\vv_h^+,s)\right), \ \B_{\V}\left(\Path(\vv_h^-,\vv_h^+,s)\right) 
\right] \cdot \n\right\}^- \, \frac{\partial \Path}{\partial s} \, ds 
= %\half  
\mathsf{D}_{ \mathbf{ n}}^- 
 \left( \vv_h^+ - \vv_h^- \right),
\label{eqn.pathint1} 
\end{eqnarray}
having set
\begin{eqnarray}
 \mathsf{D}_{ \mathbf{ n}}^- := \int \limits_0^1\left(
\left[ \A_{\V}\left(\Path(\vv_h^-,\vv_h^+,s)\right), \ \B_{\V}\left(\Path(\vv_h^-,\vv_h^+,s)\right) \right] \cdot \n \right)^-\!\!ds ,
\label{eqn.pathint2} 
\end{eqnarray}
where the matrix sign is computed using standard eigenvalue decomposition 
and having used the usual segment path approximation
%where a simple straight line segment has been chosen as integration path, namely
%\begin{equation}
$\Path = \Path(\vv_h^-,\vv_h^+,s) = \vv_h^- + s (\vv_h^+ - \vv_h^-)$. 
%\label{eqn.segpath} 
%\end{equation} 
The values $(\vv_h^-,\vv_h^+)$ are the traces of  the predictor solutions $\vv_h$ of neighboring 
elements $\omega_i$ and  $\omega_j$ on their shared   boundary.

\subsection{\textit{A posteriori} subcell finite volume limiter}
\label{ssec.limiter}

To guarantee the non-oscillatory  character of the solution, we use an \textit{a posteriori} limiting strategy
relying on a subcell finite volume update.  
This approach involves four main steps: 
i) the construction of a conformal sub-triangulation of each cell,
ii) the definition of a troubled cell indicator,
iii) a non-linear monotonicity preserving subcell finite volume update and 
iv) a reconstruction operator to recover a DG polynomial from the subcell averages.
We briefly recall these steps hereafter, and refer the interested reader 
to~\cite{DGLimiter3,ALEDG,Gaburro2021PNPMLimiter,kuzmin2022limiter} and references therein for more details.  
%\cite{SonntagDG,DGLimiter1,DGLimiter3,vilar2019posteriori,Gaburro2021PNPMLimiter,kuzmin2020subcell,rueda2021subcell,ramirezgassner-limiter2022}. 

\paragraph{Subcell definition and averaging} We introduce a sub-triangulation of each $\omega_i$ 
composed by a set of $N_{\tau, i}$ non-overlapping   sub-triangles $\tau_{i, \alpha}$ with $\alpha \in [1, N_{\tau, i}]$.
To construct them, we start by subdividing each polygonal cell in $N_{k,i}$ sub-triangles (the same 
used for integration purposes) by connecting the polygon's vertices $\x_{k,i}$ with the cell 
barycenter~$\x_{{b}_i}$. Then, we further subdivide these triangles in $N^2$ 
sub-triangles  obtaining a total number of $N_{\tau, i}=N_{k,i} \cdot N^2$ subcells. 
Note that the number of subcells
grows with the order of accuracy $N$ of the method in order to always ensure a proper subcell resolution
and allow the local reconstruction of a continuous 
polynomial starting from cell averages (cf. later). In particular, in  each 
sub-triangle $\tau_{i,\alpha}$ we define  the subcell average of the  solution at time $t^n$  
\begin{equation}
	\mathbf{v}_{i,\alpha}^n  = \frac{1}{|\tau_{i,\alpha}|} \int_{\tau_{i,\alpha}} \mathbf{v}_{h}^n(\x,t^n) \, d\x = \frac{1}{|\tau_{i,\alpha}|}
	\sum \limits_{\ell=0}^{\mathcal{N}-1}\Big[\int_{\tau_{i,\alpha}} \phi_\ell(\x) \, d\x \Big]
	 \, \hat{\mathbf{v}}^{n}_{\ell}:=\mathcal{P}(\mathbf{v}_h^n) \qquad \forall \alpha \in [1,N_{\tau, i}],
	\label{eqn.vh}
\end{equation}
where $|\tau_{i,\alpha}|$ denotes the volume of subcell $\tau_{i,\alpha}$ of element $\omega_i$ and 
the definition $\mathcal{P}(\mathbf{v}_h)$ is a $L_2$ projection operator
onto the space of degree zero polynomials. 
Similarly, starting from the solution  $\mathbf{v}_h^{n+1}$ obtained in the corrector update of 
Section~\ref{ssec.corrector}, 
we obtain the \textit{candidate} subcell averages of the numerical solution at time $t^{n+1}$ as 
$\mathbf{v}_{i,\alpha}^{n+1,*} = \mathcal{P}(\mathbf{v}_h^{n+1})$. 
%Note that we have denoted the numerical solution computed so far by $\mathbf{v}^{n+1,*}_{h}(\x,t^{n+1})$ 
%as it should be considered a \textit{candidate} solution. 
%Indeed, we now need to check it against a set of detection criteria before validating it. \\

\paragraph{Troubled cell detection} 
We use a mixed criterion combining thermodynamic admissibility with a local maximum principle.
In practice, we require that 
% the requirement that the computed solution is physically acceptable, i.e. belongs to the phase space of the conservation 
% law being solved. For instance, if the compressible Euler equations for gas dynamics are considered, 
density and pressure should be greater than a prescribed tolerance $\epsilon=10^{-12}$ and that the following
relaxed discrete maximum principle (DMP) should be verified
\begin{equation}
	\min_{j \in \mathcal{V}(\omega_i)} \left( \, \min_{\beta\in[1,{N}_{\tau, j}]} (\mathbf{v}_{\jbeta}^n) \right) 
	- \delta \leq \ \mathbf{v}^{n+1,*}_{i,\alpha} \ \leq 
	\max_{j \in \mathcal{V}(\omega_i)} \left( \, \max_{\beta\in[1,{N}_{\tau, j}]} (\mathbf{v}_{\jbeta}^n) \right) 
	+ \delta \qquad \forall \alpha \in [1,{N}_{\tau, i}],
	\label{eqn.RDMP}
\end{equation}
where $\mathcal{V}(\omega_i)$ is the set of all the neighbors of $\omega_i$ and $\delta$ is a parameter defined as
\begin{equation}
	\delta = \max \, \Biggl( \delta_0\, , \ \epsilon \cdot 
	\biggl[
	\max_{j \in \mathcal{V}(\omega_i)} \biggl(\, \max_{\beta\in[1,N_{\tau, j}]} (\mathbf{v}_{\jbeta}^n) \biggr) -
	\min_{j \in \mathcal{V}(\omega_i)} \biggl(\, \min_{\beta\in[1,N_{\tau, j}]} (\mathbf{v}_{\jbeta}^n) \biggr) 
	\biggr]
	\Biggr),
	\label{eqn.deltaRDMP}
\end{equation}
with $\delta_0=10^{-4}$ and $\epsilon=10^{-3}$.  

A cell is marked as troubled if any of the above conditions is violated.

%\subsection{\textit{A posteriori} conservative finite volume limiter}

\paragraph{Subcell finite volume update}   In troubled subcells  the candidate  averaged 
values $\mathbf{v}^{n+1,*}_{i,\alpha}$
are replaced by  $\mathbf{v}^{n+1}_{i,\alpha}$ obtained using  a path-conservative 
MUSCL-Hancock~\cite{vanleer1974, vanleer1979, torobook} finite volume method, 
equivalent to a  second order ADER-FV scheme.   The method consists of two steps.
%
%scheme that evolve the cell averages $\mathbf{v}_{i,\alpha}^n$ associated to each troubled polygons is organized as follows:
%in both cases, the scheme can be seen to evolve the primitive subcell data $\mathbf{v}_{i,\alpha}^n$
%to the next time level generating $\mathbf{v}^{n+1}_{i,\alpha}$ for each subcell. 
%
 First we evaluate in each subcell  $\tau_{i,\alpha}$ a local predictor  
\begin{equation}\label{eqn.local.predictor}
\vv_{i,\alpha}(x, y, t) = \mathbf{v}_{i,\alpha}^n - \left(\mathbf{A}_\mathbf{V}\left(\mathbf{v}_{i,\alpha}^n\right)\,\mathbf{S}_x + 
\mathbf{B}_\mathbf{V}\left(\mathbf{v}_{i,\alpha}^n\right)\,\mathbf{S}_y\right)\,(t - t^n) +\mathbf{S}_x\,(x - x_{i,\alpha}) + \mathbf{S}_y\,(y - y_{i,\alpha}).
\end{equation}
where $\mathbf{S}_x$ and $\mathbf{S}_y$ are slopes computed using the primitive subcell data $\mathbf{v}_{i,\alpha}^n$,  
along with  the  values  $\mathbf{v}_{\jbeta}^n$ relative to  neighboring subcells $\tau_{\jbeta}$
sharing at least a vertex with the  $\tau_{i,\alpha}$. These slopes are obtained from a standard ENO reconstruction operator, 
selecting among a set of candidate reconstruction polynomials the one such that $\mathbf{S}_x^2 + \mathbf{S}_y^2$ is minimized.
This stencil selection is carried out independently variable by variable.
  
%Then a predictor solution is evaluated at each of the three space-time midpoints of the space-time faces
%of the element $\tau_{i,\alpha}$. The spatial location of these quadrature nodes are denoted with the coordinates $x_{i,\alpha,k}$ and $y_{i,\alpha,k}$,
%and their time coordinate is simply $t^{n+1/2} = (t^n + t^{n+1})/2$.
%Additionally, we denote $\tau_{\jbeta}$ the sub-triangle sharing the face midpoint $x_{i,\alpha,k}$ with the sub-triangle $\tau_{i,\alpha}$, 
%which has area $|\tau_{i,\alpha}|$ and for edge lengths denoted by $|\tau_{i,\alpha,k}|$.
%Formally the predictor solution for each space-time coordinate within the element $\tau_{i,\alpha}$ is computed as 
%\begin{equation}
%\mathbf{q}_{i,\alpha}(x, y, t) = \mathbf{v}_{i,\alpha}^n - \left(\mathbf{A}_\mathbf{V}\left(\mathbf{v}_{i,\alpha}^n\right)\,\mathbf{S}_x + 
%\mathbf{B}_\mathbf{V}\left(\mathbf{v}_{i,\alpha}^n\right)\,\mathbf{S}_y\right)\,(t - t^n) + \mathbf{w}_{i,\alpha}(x,\ y).
%\end{equation}

Once all the subcell predictors  have been computed,  the primitive cell average values are obtained as  
\begin{equation}\label{eqn.subcell-corrector}
\begin{aligned}
	\mathbf{v}_{i,\alpha}^{n+1} = \mathbf{v}_{i,\alpha}^{n} &- \Delta t\,\left(
	\mathbf{A}_\mathbf{V}\left(\vv_{i,\alpha}(x_{i,\alpha},\, y_{i,\alpha},\, t^{n+1/2}\right)\,\mathbf{S}_x + 
	\mathbf{B}_\mathbf{V}\left(\vv_{i,\alpha}(x_{i,\alpha},\, y_{i,\alpha},\, t^{n+1/2}\right)\,\mathbf{S}_y\right) + \\
	&- \frac{\Delta t}{|\tau_{i,\alpha}|} \sum_{k=1}^3  
	|\partial\tau_{i,\alpha,k}|\, \mathbf{D}\left(\vv_{i,\alpha}(x_{i,\alpha,k},\, y_{i,\alpha,k},\, t^{n+1/2}),\, 
	\vv_{j,\beta}(x_{i,\alpha,k},\, y_{i,\alpha,k},\, t^{n+1/2})\right)\cdot\mathbf{n},
\end{aligned}
\end{equation}
where we denote with $\x_{i,\alpha,k}$ the mid-point of the k-th face of 
subcell $\tau_{i,\alpha}$, with  $\partial\tau_{i,\alpha,k}$ the length of the corresponding edge 
and $\mathbf{D}\left(\vv^{-}, \vv^{+}\right)\cdot \mathbf{n}$ is computed according to \eqref{eqn.pathint1}.

\paragraph{DG polynomial reconstruction} Once valid cell averages $\mathbf{v}^{n+1}_{i,\alpha}$ have been computed with the above procedure,  
the DG polynomial $\mathbf{v}_{h}(\x,t^{n+1})$ in the polygonal cell $\omega_i$ is obtained 
via a reconstruction operator $\mathcal{R}(\mathbf{v}_{i,\alpha}^{n+1})(\x)$ such that 
\begin{equation}
	\int_{\tau_{i,\alpha}} \mathbf{v}_{h}(\x,t^{n+1}) \, d\x := \int_{\tau_{i,\alpha}} \mathcal{R}(\mathbf{v}_{i,\alpha}^{n+1})(\x) \, d\x=
	|\tau_{i,\alpha}| \mathbf{v}_{i,\alpha}^{n+1}
%	\int_{\tau_{i,\alpha}} \mathbf{v}_{i,\alpha}^{n+1}(\x,t^n) \, d\x :=\mathcal{R}(\mathbf{v}_{i,\alpha}^{n+1}(\x,t^n)) 
	 \qquad \forall \alpha \in [1,N_{\tau, i}],
	\label{eqn.intRec}
\end{equation}
where by construction we also have 
\begin{equation}
	\int_{\omega_i} \mathbf{v}_{h}(\x,t^{n+1}) \, d\x = \sum_{\alpha\in[1,N_{\tau, i}]}  |\tau_{i,\alpha}| 
	\int_{\tau_{i,\alpha}} \mathbf{v}_{i,\alpha} \, d\x.
	\label{eqn.LSQ}
\end{equation}

Although the projection and reconstruction operators \eqref{eqn.vh} and~\eqref{eqn.intRec}-\eqref{eqn.LSQ}
 formally  satisfy  the property $\mathcal{P} \cdot \mathcal{R}=\mathcal{I}$, with $\mathcal{I}$ being the identity 
operator, their linearity does not fully guarantee the absence of oscillations in the reconstructed solution. 
As   a safety mesure,  the  subcell averages $\mathbf{v}_{i,\alpha}^{n+1}$ are stored  as long as 
the cell $\omega_i$ is not marked as valid again. If the cell is detected to be troubled 
for a second time step in a row, then  the stored subcell averaged are used instead of those 
obtained via the projection.

Note also that, similarly to what is done in~\cite{Gaburro2021PNPMLimiter}, no corrections are 
introduced to make sure that boundary integrals involving limited and non-limited cells
match. To account for this, and for robustness reasons,   when a shock lies near the boundary of a limited cell
the troubled cell flag is extended also to all the direct neighbours of the elements flagged as troubled by the admissibility 
criteria.  The issue of non-matching boundary integrals in the subcell a posteriori limiting is discussed later in this paper.

\section{Quasi-conservative formulation via \textit{a posteriori} local conservative subcell corrections}\label{sec_conservation}

% Mario
%The method described in the previous section has many attractive features: 
%it is arbitrary high order accurate for smooth solutions;
%it embeds a very high resolution non-oscillatory approximation of discontinuities thanks to the subcell limiter;
%by allowing to directly evolving physical variables, 
%it is  free of all the spurious numerical effects recalled  in the
%introduction  when approximating interfaces (e.g. in multi-material flows) and contact discontinuities.

% Elena
The method described in the previous section has many attractive features: 
indeed i) it is arbitrary \textit{high order accurate} for smooth solutions, 
ii) at the same time, thanks to the subcell limiter, 
it also provides very high resolution, \textit{free from oscillatory behaviors}, on discontinuities and, 
iii) by allowing to evolving directly the physical variables, 
it is not affected by all the spurious numerical artifacts, 
recalled in the introduction, 
when approximating \textit{interfaces} (e.g. in multi-material flows) and \textit{contact discontinuities}.

In this section, we propose local corrections allowing to retain also a correct approximation of \textit{shocks}. 
Our hybrid approach has {commonalities} with other works mentioned in the introduction, 
which also enforce conservation only locally at shocks. 
The originality of our method is that it is naturally formulated
within the very high order ADER predictor-corrector approach with \textit{a posteriori} limiting, 
in which we now include {modifications} that restore the conservation.
To give meaning to our methodology, 
we start by proposing a modified {view of the conditions to approximate} weak solutions, 
and in particular on the local conservation conditions 
which serve as a base for the usual proofs of the Lax-Wendroff theorem.

\subsection{Local conservation, Lax-Wendroff theorem, and  vanishing conservation defect}

The classical characterization {of \textit{conservation}, for the purposes of discretizing hyperbolic systems} %non mi piace come suona questa fr
% ase ... per me mancano gli articoli o l'ordine delle parole stona ...
is associated to a local discrete balance which is usually written as
\begin{equation}
	|T_i|\dfrac{\bar{\mathbf{q}}_i^{n+1} - \bar{\mathbf{q}}_i^n}{\Delta t}  +  \sum\limits_{{\ell}} \mathbf{F}_{\n_{i\ell}} =0,
	\label{eqn.localConservation0}
\end{equation}
where   $\bar{\mathbf{q}}_i$ are  local  averages of the conservative variables $\Q$ in a generic cell $T_i$,  
$\mathbf{F}_{\n_{i\ell}}$ are face averaged numerical fluxes
and $\n_{i\ell}$ face normal vectors scaled by the length of the face between cell $i$ and its $\ell^\text{th}$ neighbor.

For the usual ADER-DG methods based on the conservative form  of the equations  in the corrector 
step (see e.g.~\cite{Dumbser2008,ZanottiDumbser2016,Gaburro2021PNPMLimiter}), the above statement can be written
for all  the cell averages
% purtroppo in tutti que paper qh e' il predictor, non il polinomio a t^n
\begin{equation}
\bar{\mathbf{q}}_{i} := \dfrac{1}{|\omega_i|}\int\limits_{\omega_{i}}\mathbf{q}_hd\mathbf{x},
	\label{eqn.localConservation1}
\end{equation}
with $\mathbf{q}_h$ the spatial polynomial approximation~\eqref{eqn.uh}
of the conservative variables written with the basis functions~\eqref{eq.space_basis}. 
%The above statement is naturally associated to the first mode in \eqref{eq.space_basis}, which represents the evolution of the cell 
% average. \wb{NO: the basis \eqref{eq.space_basis} is not conservative in the sense that $\int_{\omega_i} \phi \, d\mathbf{x} \neq 1$ 
% for $N>1$, so I would remove this last sentence.}
The fluxes appearing in~\eqref{eqn.localConservation0} are, in this case,
given by the integral over a time interval with
the intermediate values in time being obtained from the space-time predictor.

In the context of schemes \wb{using} subcell finite volume limiters, 
one can also consider a local conservation written at the level of the subcells
\begin{equation}
|\tau_{i,\alpha}|\dfrac{\bar{\mathbf{q}}_{i,\alpha}^{n+1} -\bar{\mathbf{q}}_{i,\alpha}^{n} }{\Delta t} 
 +  \sum\limits_{{\ell_{i,\alpha}}} \mathbf{F}_{\n_{{\ell_{i,\alpha}}}} =0, 
	\label{eqn.localConservation2}
\end{equation}
where $\bar{\mathbf{q}}_{i,\alpha}$ is the sub-element average
\begin{equation}
\bar{\mathbf{q}}_{i,\alpha} := \dfrac{1}{|\tau_{i,\alpha}|}\int\limits_{\tau_{i,\alpha}}\mathbf{q}_hd\mathbf{x}.
	\label{eqn.localConservation3}
\end{equation}
Clearly   if \eqref{eqn.localConservation2} is verified  $\forall\alpha$, so is \eqref{eqn.localConservation0} due to the identity  
\begin{equation}\label{eqn.subtocellaverage}
|\omega_i| \bar{\mathbf{q}}_{i}  = \sum_{\alpha\in[1,N_{\tau, i}]}  |\tau_{i,\alpha}| \bar{\mathbf{q}}_{i,\alpha}
\end{equation}
and the fact that numerical fluxes on internal subcell faces cancel out when summing on $\alpha$.

The local conservation condition \eqref{eqn.localConservation0} can be used to show 
a \textit{Lax-Wendroff theorem} guaranteeing that, if convergent, 
the discrete solution converges to a weak solution of the conservation law. 
For details on this result on several types of discretizations we refer e.g. to~\cite{zbMATH02026806,SHI20183}. 

{The main idea} is to prove that,
under some relatively classical approximation results, 
plus a bounded variation assumption (which will be recalled shortly),
the satisfaction of a slightly generalized version of~\eqref{eqn.localConservation0}, accounting also for a conservation defect, 
also implies that
\begin{equation}\label{eqn.LW0}
\int_{\Omega}\mathbf{q}_h^0\phi + \int_{\mathbb{R}^+} \int_{\Omega}\mathbf{q}_h \partial_t\phi 
d\mathbf{x}\,dt+  \int_{\mathbb{R}^+} \int_{\Omega}\mathbb{F}_h(\mathbf{q}_h )\cdot \nabla\phi d\mathbf{x}\,dt
+ \smallO(1) =0,
\end{equation}
for any function $\phi\in C^{\infty}_0(\Omega)$ {with $\mathbf{q}_h^0$ being the discrete initial condition} and $\smallO(1)$ a quantity vanishing when $h \rightarrow 0$.
This is important because~\eqref{eqn.LW0} implies, in particular, that if $\mathbf{q}_h \rightarrow \mathbf{q}$ 
then the latter is by definition a weak solution of the PDE system~\eqref{eqn.PDE}. 

For the sake of completeness, we would like to recall that classically~\eqref{eqn.LW0} can be obtained from~\eqref{eqn.localConservation0}
by studying the quantity
\begin{equation} 
\sum\limits_n\sum\limits_{i} |T_i|(\bar{\mathbf{q}}_i^{n+1} - \bar{\mathbf{q}}_i^n)\phi_i^n + 
\sum\limits_n\sum\limits_{i} \Delta t\sum\limits_{{\ell}} \phi_i^n\mathbf{F}_{\n_{i\ell}},
\end{equation} 
with $\phi_i^n$ the values of $\phi$ at the cell centers at time $t^n$. 
A similar procedure will be also used in what follows.

\smallskip

Indeed, we would like to suggest a simple modification of the local conservation condition~\eqref{eqn.localConservation0},
{which we may think of writing as} 
\begin{equation}
	|\omega_i|\dfrac{\bar{\mathbf{q}}_i^{n+1} - \bar{\mathbf{q}}_i^n}{\Delta t}  +  
	\sum\limits_{{\ell}} \mathbf{F}_{\n_{i\ell}} = 	|\omega_i|\Delta_i^n
	\label{eqn.localQuasiConservation0}
\end{equation}
where $\Delta_i^n$ is a \textit{local conservation defect}. 
It is now easy to show the following result. 

\begin{proposition}[Local quasi-conservation condition and  Lax-Wedroff theorem]\label{prop:quasi-cons}
Consider a scheme  whose discrete solution $\mathbf{q}_h$ verifies   the modified local 
conservation statement \eqref{eqn.localQuasiConservation0}.
Assume that the discrete  initial condition $\mathbf{q}_h^0$ verifies the weak convergence estimate
$$
\int\limits_{\Omega}(\mathbf{q}_h^0 - \mathbf{q}_0)\phi \rightarrow 0\;\text{ as }\; h\rightarrow  
0\;\;\forall \phi\in C^{\infty}_0(\Omega)  \ \text{ with } \ \wb{\mathbf{q}_0=\mathbf{q}(\mathbf{x},t=0)}
$$
and that the discrete solution also verifies, in a ball $B_{i} $ centered around the center of the cell $\omega_i$,
the total-variation-{boundedness} condition in the appropriate norm
\begin{equation}\label{eqn.TVB}
h^2\sum_i\sum_{\alpha}\max_{\mathbf{x} \in B_{i,\alpha} }\|\mathbf{q}_h^n(\mathbf{x}) - 
\mathbf{q}_h^n(\mathbf{x}_{i,\alpha})  \| \rightarrow 0\;\text{ as }\; h\rightarrow  0\;\;\forall n\ge 0.
\end{equation}
Then, on {shape-regular meshes}, provided that
\begin{equation}\label{eqn.quasi-cons-cond}
\sup_{i,n} \|\Delta_{i}^n \| = \smallO(1),
\end{equation}
  if   $\mathbf{q}_h^n$ converges to $\mathbf{q}$,  then $\mathbf{q}$ is a weak solution of the 
  system of conservation laws in the sense of \eqref{eqn.LW0}. 
 \begin{proof}
The proof is obtained  following  a classical strategy.
First, we multiply~\eqref{eqn.localQuasiConservation0} by the values
$\phi_{i}^n$ of \wb{an arbitrary $\phi\in C^{\infty}_0(\Omega)$ taken at the cell centers} and time $t^n$ 
and we take a sum over $n$ and $i$. 
Then, following step by step the analysis of~\cite{SHI20183} Section 3.3 (or Appendix~A of~\cite{zbMATH02026806}, page 31), 
we can show that \eqref{eqn.localQuasiConservation0} 
implies  the satisfaction of
$$
\int_{\Omega}\mathbf{q}_h^0\phi + \int_{\mathbb{R}^+} \int_{\Omega}\mathbf{q}_h \partial_t\phi 
d\mathbf{x}\,dt+  \int_{\mathbb{R}^+} \int_{\Omega}\mathbb{F}_h(\mathbf{q}_h )\cdot \nabla\phi d\mathbf{x}\,dt +
\sum_{i,n}\Delta t \, |\omega_i| \, \phi_{i}^n \, \Delta_{i}^n
+ \smallO(1)  =0.
$$
On {shape-regular meshes}, the additional term w.r.t.~\eqref{eqn.LW0} can be easily bounded as
$$
\left \| \, \sum_{i,n}\Delta t \, |\omega_i| \, \phi_{i}^n \, \Delta_{i}^n \, \right  \| 
\ \le \
| \Omega| \, t_f  \left \| \phi \right \|_{L^{\infty}(\Omega\times[0,t_f])} \, \sup_{i,n} \left \|\Delta_{i}^n \right \|
$$
and thus, {by hypothesis \eqref{eqn.quasi-cons-cond}}, can be included in the $\smallO(1)$ remainder 
of the Lax-Wendroff theorem, 
allowing to recover~\eqref{eqn.LW0}. 
 \end{proof}
\end{proposition} 
This proposition justifies the construction of ({possibly} non-linear) discretizations  
which allow conservation defects throughout the computational domain 
whenever the above conditions on the smallness of this defect apply.
Note that the non-linearity of the discretizations  
is an important {element guaranteeing that~\eqref{eqn.TVB} 
is verified through some form of oscillation control, as e.g.} 
the \textit{a posteriori} subcell limiter discussed in Section~\ref{ssec.predictor}. % la frase è interessante, ma non sta in piedi.
The discretizations obtained in such a way and verifying Proposition~1, 
will be defined as \textit{locally quasi-conservative}.

\smallskip

Now, {our objective of the next sections is 
to first} recast non-conservative discretizations in the setting of the above proposition. 
{Then, we aim at proposing} some heuristics allowing to define a scheme for 
which~\eqref{eqn.localQuasiConservation0} holds true with a non-zero conservation defect only
in cases in which~\eqref{eqn.quasi-cons-cond} is expected to be true.

\subsection{Conservation defect for non-conservative schemes} 

Let us now consider a scheme {which provides discrete evolution equations starting from 
some non-conservative form.}
We will essentially focus on the temporal update defined by the {path-integral-based} 
ADER-DG scheme~\eqref{eqn.DGscheme}, but the discussion can be easily adapted to other approaches. 

We begin by considering {the evolution} of the cell average for {every} cell $\omega_i$, 
which can be computed explicitly as 
\begin{equation}
|\omega_i|\dfrac{\bar{\mathbf{q}}_{i}^{n+1} -\bar{\mathbf{q}}_{i}^{n} }{\Delta t}  =  
\dfrac{1}{\Delta t}\left(\,\int\limits_{\omega_i} \mathbf{Q}(\mathbf{v}_h^{n+1}) d\mathbf{x} -
\int\limits_{\omega_i} \mathbf{Q}(\mathbf{v}_h^{n}) d\mathbf{x}
\right)
=   \dfrac{1}{\Delta t}\left(\,\int\limits_{\omega_i} \mathbf{Q}(\mathbf{v}_h^{n} +   
\dfrac{\Delta t}{|\omega_i|}\delta\mathbf{v}_h   ) d\mathbf{x} -
\int\limits_{\omega_i} \mathbf{Q}(\mathbf{v}_h^{n}) d\mathbf{x}
\right),
	\label{eqn.localConservation3a}
\end{equation}
where $\mathbf{Q}(\V)$ is the mapping between physical and conservative variables 
discussed in Section~\ref{sec_pde} 
and the scaled evolution operator  $\delta\mathbf{v}_h $   can be  deduced 
from \eqref{eqn.DGscheme}. % (see Appendix \ref{app:ADER-DG_loc1} for the detailed expression).
Let us now  consider for each element the time averaged flux balance 
\begin{equation}\label{eqn.flux-balance}
\Phi_{i} :=   \sum\limits_{\ell_{i}} \mathbf{F}_{\n_{{\ell_{i}}}} = \dfrac{1}{\Delta t} 
\int\limits_{t^n}^{t^{n+1}} \int\limits_{\partial\omega_{i}} \mathbf{F}_{\n}(\tilde{\mathbf{v}}_h^-,\, \tilde{\mathbf{v}}^+_h)dS, 
\end{equation}
where the space-time polynomials of the physical variables $\tilde{\mathbf{v}}_h$ 
can be obtained \wb{by means of} some explicit predictor strategy,
as e.g the one discussed in Section~\ref{ssec.predictor}.  
The $\mathbf{F}_{\n_{{\ell_{i}}}} $ are now temporal and edge averaged consistent numerical fluxes,
to be defined  using  values of the local predictors  within $\omega_i$ and its neighbors. 
To cast the schemes in the setting of Proposition \ref{prop:quasi-cons},  we can now define 
the \textit{cell conservation defect} as follows
 \begin{equation}\label{eqn.defect}
\Delta_{i}^n:=  \dfrac{1}{|\omega_{i}|} \dfrac{1}{\Delta t}\left(\,\int\limits_{\omega_{i}} 
\mathbf{Q}(\mathbf{v}_h^{n} +   \dfrac{\Delta t}{|\omega_i|}\delta\mathbf{v}_h   ) d\mathbf{x} -
\int\limits_{\omega_{i}} \mathbf{Q}(\mathbf{v}_h^{n}) d\mathbf{x}
\right)
+ \dfrac{\Phi_{i}}{|\omega_{i}|}.
\end{equation}
The non-conservative scheme \eqref{eqn.DGscheme} can be easily shown to  
verify  \eqref{eqn.localQuasiConservation0} with the above definition
of $\Delta_i^n$.

Note that, to account for the effects of subcell limiting, we can perform a very similar 
construction for the subcell triangulation, and introduce (with similar notations)
subcell conservation defects
 \begin{equation}\label{eqn.defectSub}
\Delta_{i,\alpha}^n:=  \dfrac{1}{|\tau_{i,\alpha}|} \dfrac{1}{\Delta t}\left(\,\int\limits_{\tau_{i,\alpha}} 
\mathbf{Q}(\mathbf{v}_h^{n} +   \dfrac{\Delta t}{|\tau_{i,\alpha}|}\delta\mathbf{v}_h   ) d\mathbf{x} -
\int\limits_{\tau_{i,\alpha}} \mathbf{Q}(\mathbf{v}_h^{n}) d\mathbf{x}
\right)
+ \dfrac{\Phi_{i,\alpha}}{|\tau_{i,\alpha}|},
\end{equation}
where 
\begin{equation}\label{eqn.flux-balanceSub}
\Phi_{i,\alpha} :=  \sum\limits_{{\ell_{i,\alpha}}} \mathbf{F}_{\n_{{\ell_{i,\alpha}}}} = 
\dfrac{1}{\Delta t} \int\limits_{t^n}^{t^{n+1}} \int\limits_{\partial \tau_{i,\alpha}} \mathbf{F}_{\n}(\tilde{\mathbf{v}}_h^-,\, \tilde{\mathbf{v}}^+_h)dS. 
\end{equation}
We can now readily write the local quasi-conservation statement
\begin{equation}
|\omega_i|\dfrac{\bar{\mathbf{q}}_{i,\alpha}^{n+1} -\bar{\mathbf{q}}_{i,\alpha}^{n} }{\Delta t}   + 
  \sum\limits_{{\ell_{i,\alpha}}} \mathbf{F}_{\n_{\ell_{i,\alpha}}}  = |\tau_{i,\alpha}|\Delta_{i,\alpha}^n,
	\label{eqn.localConservation4a0}
\end{equation} 
where now the numerical fluxes $\mathbf{F}_{\n_{\ell_{i,\alpha}}}$ are time averaged values 
arising  naturally in the ADER framework, and 
implicitly defined in \eqref{eqn.flux-balanceSub}. 
Using \eqref{eqn.subtocellaverage}, the above relation immediately translates in a cell local (quasi-)conservation statement 
\begin{equation}
|\omega_i|\dfrac{\bar{\mathbf{q}}_{i}^{n+1} -\bar{\mathbf{q}}_{i}^{n} }{\Delta t}   +   
 \sum\limits_{{\ell}} \mathbf{F}_{\n_{i\ell}}  = \sum_{\alpha}|\tau_{i,\alpha}|\Delta_{i,\alpha}^n .
	\label{eqn.localConservation4a}
\end{equation}  

The main idea of our construction is to modify our scheme described in Section~\ref{sec_numscheme} in such a way that, 
in appropriately flagged elements, 
essentially those containing shocks, the conservation defect is removed. 
To this end we will first briefly  study the behavior of $\Delta_i^n$ and $\Delta_{i,\alpha}^n$ for some special and well known cases:
regular solutions and contact discontinuities.  
Indeed, under appropriate hypotheses, we can show that in these cases the defect vanishes or is
 small enough for Proposition \ref{prop:quasi-cons} to hold.
Then, we will discuss a heuristic construction to remove the error defects {at} shocks.

\subsection{Vanishing conservation defect: examples and heuristics}

\subsubsection{Regular  solutions.}
When the solution is smooth, the classical strong form of the PDE is just as {relevant} as the weak form. 
{In some sense, in this case a Lax-Wendroff theorem is not really necessary}. 
Nevertheless, it is interesting to see {how} this case falls naturally in the formalism of Proposition \ref{prop:quasi-cons}.
In particular, the continuity properties of the transformations between different sets of variables, 
explicitly recalled in section \ref{sec_pde},  
allow to derive an explicit estimate on the conservation defect which is summarized in the property below. 
Note that, everywhere it is assumed that the integrals
arising from the {variational} formulations are evaluated with quadrature formulas 
exact w.r.t. the degree of the employed spatial polynomials.

\begin{proposition}[Conservation defect:  locally regular solutions\label{prop.smooth}]
Let $B_{i} $ be  a ball centered around the barycenter $\x_i$ of the cell $\omega_i$, 
and such that each $\x$ in $\omega_i$ or $\partial\omega_i$ is contained in $B_i$. 
Assume that 
for some bounded $C > 0$, $C_{\partial}> 0 $
the discrete  solution verifies the smoothness bounds
\begin{equation}\label{eqn.smooth-HP}
\max\limits_{\x\in B_{i} } \max(\|\v_h(\x)-\v_h(\x_{i})\|,\, h\|\nabla \v_h\| ) \le C\, h 
\end{equation}
as well as 
\begin{equation}\label{eqn.smooth-HP1}
\|\v_h(\x^+)-\v_h(\x^-) \| \le C_{\partial}\, h%^{1+\beta}
\;,\quad \forall \x \in\partial \omega_i\,.
\end{equation}
Then, under the hypotheses of Proposition \ref{prop:quasi-cons}  there exist    numerical fluxes  verifying the classical
consistency and continuity assumptions such that 
 the non-conservative ADER-DG scheme \eqref{eqn.DGscheme}    is
% and  the conservative ADER-DG corrector \eqref{eqn.DGscheme-cons} with locally mismatching fluxes are
   quasi-conservative in the sense of Proposition \ref{prop:quasi-cons}. More precisely the associated  subcell conservation defect \eqref{eqn.defect} 
   %and \eqref{eqn.defect-cons}    verify 
   verifies  $\Delta^n_{i} =o(1)$. 
 \begin{proof}
 See \ref{app:prop-smooth-proof}. 
 \end{proof}
\end{proposition}

This result puts the treatment of smooth areas with non-conservative terms in the modified setting of the Lax-Wendroff theorem
associated to Proposition \ref{prop:quasi-cons}.  The hypotheses  \eqref{eqn.smooth-HP}  and  \eqref{eqn.smooth-HP1}  on the discrete solution
correspond to bounded first  derivatives. 
In practice, {at elements interfaces}, for smooth data, one could expect small jumps of order 
$\|\v_h(\x^+)-\v_h(\x^-) \| =\mathcal{O}(h^{p+1}) $  for an approximation of degree $p$.  
These hypotheses are in any case stronger conditions than the TVB requirement \eqref{eqn.TVB}  used in the usual  Lax-Wendroff theorem proof by~\cite{SHI20183}.  \\
 
 The above result can be extended in two ways. The first relates to the case in which the subcell limiter is activated 
so the non-conservative  update is not given  by  \eqref{eqn.DGscheme}, but by the subcell path conservative
finite volume corrector \eqref{eqn.subcell-corrector}. In this case, we can use the definitions  \eqref{eqn.defectSub} of the conservation defects 
in each subcell, for which we can prove the following property. 

\begin{proposition}[Conservation  defect:  regular solutions and subcell non-conservative update\label{prop.smooth.sub}] 
Let $B_{i,\alpha} $ be  a ball    centered around the center of the subcell $\tau_{i,\alpha}$ and containing at least its boundary $\partial\tau_{i,\alpha}$.
 Assume that  for some bounded $C > 0$ %, $C_{\partial}> 0$ and $\beta>0$ 
the discrete  solution verifies the smoothness bounds  
\begin{equation}\label{eqn.smooth-HP-sub}
\max\limits_{\x\in B_{i,\alpha} } \max(\|\v_h(\x)-\v_h(\x_{i})\|,\, h\|\nabla \v_h\| ) \le C\, h
\end{equation}
as well as 
\begin{equation}\label{eqn.smooth-HP1-sub}
\|\v_h(\x^+)-\v_h(\x^-) \| \le C_{\partial}\, h\;,\quad \forall \x \in\partial\tau_{i,\alpha} \,.
\end{equation}
Then, under the hypotheses of Proposition \ref{prop:quasi-cons}  there exist    numerical fluxes  verifying the classical
consistency and continuity assumptions such that   the non-conservative subcell update  \eqref{eqn.subcell-corrector}   
can be characterized by the local \mbox{(quasi-)}conservation property \eqref{eqn.localConservation4a0}, where
the conservation defect is defined by \eqref{eqn.defectSub}, and verifies   $\Delta^n_{i,\alpha} =o(1)$. 
In particular, if the above is true $\forall\alpha$, then on shape regular meshes  \eqref{eqn.localQuasiConservation0} is true with 
$$
\Delta_i^n := \sum_{\alpha}\dfrac{|\tau_{i,\alpha}|}{|\omega_i|}\Delta^n_{i,\alpha} =o(1)
$$
and the scheme is  quasi-conservative in the sense of Proposition \ref{prop:quasi-cons}. 
 \begin{proof}
We can show that   $\Delta^n_{i,\alpha} =\mathcal{O}(h)$ following the exact same steps  used for the evolution 
of the averages in the proof of   Proposition  \ref{prop:quasi-cons},
and  it results in  the same definition of consistent numerical flux allowing to satisfy the property.  
We refer also to \ref{app:prop-smooth-proof} for further details. 
The remainder of the proof, is a consequence of the boundedness of the ratios  $ |\tau_{i,\alpha}|/|\omega_i| $ on any shape regular mesh.
 \end{proof}
\end{proposition}

A second interesting remark  is related to the definition of the flux balance terms \eqref{eqn.flux-balance} and \eqref{eqn.flux-balanceSub},
and to an extent already accounted for within the proof reported in \ref{app:prop-smooth-proof}. More particularly,
the Lax-Wendroff proof  assumes the usual continuity of  numerical fluxes $\mathbf{F}_{\mathbf{n}_{i\ell}}$ which we can write as
$$
\mathbf{F}_{\mathbf{n}_{i\ell}}  + \mathbf{F}_{\mathbf{n}_{\ell i}} =0.
$$
However,  the above condition may be broken in some cases, even when discretizing the conservative form of the equations. 
%as long as \eqref{eqn.smooth-HP1} (or equivalently \eqref{eqn.smooth-HP1-sub}) hold one can allow 
%a mismatch in  numerical fluxes. 
For example, this can occur   if numerical fluxes in neighboring cells are not evaluated simultaneously
and with the same left and right polynomial approximations. This is always the case in \textit{a posteriori} limiting
unless some \textit{ad hoc} corrections are performed on non-troubled cells, {as proposed in~\cite{ALEDG}}. 
In this case one may consider a local quasi-conservation statement of the form
$$
\dfrac{\q_i^{n+1}-\q_i^n}{\Delta t} + \sum_{\ell}  \bar{\mathbf{F}}_{\mathbf{n}_{i\ell}} = |\omega_i|\Delta_i^n
$$
with $\bar{\mathbf{F}}_{\mathbf{n}_{i\ell}}$   the average of the (possibly mismatching) numerical fluxes
evaluated in neighboring cells sharing the face $i\ell$, and where now
$$
\Delta_i^n:= \dfrac{1}{|\omega_i|}\dfrac{1}{2}\int_{t^n}^{t^{n+1}}\int_{\partial \omega_i}
\bigg(
\mathbf{F}_{\mathbf{n}} (\tilde \v^+, \tilde \v^-) -
\mathbf{F}_{\mathbf{n}}(\tilde \v^+, \hat \v^-)
\bigg)\,dS\,dt .
$$
Now if we use the classical Lipschitz continuity property of any numerical flux function  we have
$$
\|\mathbf{F}_{\mathbf{n}_{i\ell}}(\tilde \v^+, \tilde \v^-)-\mathbf{F}_{\mathbf{n}_{i\ell}}(\tilde \v^+, \hat \v)\|
\le   K_{\mathbf{F} } \|\tilde \v^--\hat \v^- \|.
$$
In this case we can prove easily that if $ \|\tilde \v^--\hat \v^- \|=\mathcal{O}(h^{1+\beta})$ for some $\beta >0$, then
$\Delta_i^n=o(1)$, and the scheme remains quasi-conservative, in the sense of Proposition \ref{prop:quasi-cons}. 
This confirms the observation that extending sufficiently the troubled cells indicator, to a layer of regular elements around a detected problematic cell, 
it is enough to retain consistency and convergence to the correct solutions, 
thus eliminating the need of introducing corrections in the update of {non-flagged} cells~\cite{DGLimiter3}.

This can be summarized as follows.

\begin{proposition}[Conservation  defect:   data mismatch on cell boundaries\label{prop.smooth.mismatch}]
Consider a discretization for which in each cell  $\omega_i$ the local quasi-conservation balance \eqref{eqn.localQuasiConservation0} holds
for some definition of the numerical flux in \eqref{eqn.flux-balance} based on data $\tilde \v_h^-$, within $\omega_i$,  and $\tilde \v_h^+$ in its neighbors. 
Assume that for some  cell  $\omega_{\ell}$
neighboring of $\omega_i$ the same \wb{is true}, with however a numerical flux using  a discrete approximation  $\hat\v_h$ for the data in $\omega_i$ different than
$\tilde \v_h^-$. Provided that the discrete solution verifies the regularity estimate
 $$
   \tilde \v^-_h(\x)-\hat \v_h(\x)=\mathcal{O}(h^{1+\beta}) \;,\quad\forall\x\in \partial \omega_i, \qquad \wb{\beta>0},
$$
then the scheme is quasi-conservative in the sense  of Proposition \ref{prop:quasi-cons}.
\end{proposition}

%
%
% typical example could be  the case in which different polynomial approximations
%are used to evaluate the numerical flux on two sides of the same face, e.g.
%$$
%\mathbf{F}_{\mathbf{n}_{i\ell}}= \mathbf{F}_{\mathbf{n}_{i\ell}}(\tilde \v_h , \tilde \v_h^-)\;,\quad
%\mathbf{F}_{\mathbf{n}_{\ell i}} =\mathbf{F}_{\mathbf{n}_{\ell i}}( \hat \v_h^+, \tilde \v_h^-)
%$$
%where  $\hat \v_h^+$ is some local po 

 \subsubsection{Contact discontinuities for the multi-material Euler equations}% MARIUZ-contact

Contact discontinuities are another interesting, and well known, case to look at. We focus here on the particular system
used in the numerical tests, namely the multi-material Euler equations  introduced in section \S2.1. 
We consider a situation in which the flow can be  locally considered  one dimensional, and in particular we consider
contact  problems for which for a given cell $i$ there exist a ball $B_i$ centered in $\x_{b_i}$, and containing both $\omega_i$ and 
$\partial\omega_i$ such that $\forall\x\in B_i$ we have  for the flow velocity and for the pressure:
\begin{equation}\label{eqn.contact.q1d}
\pmb{\mathsf{v}} = \pmb{\mathsf{v}}_0 + c_v(\x)h^{1+\epsilon}\;,\;\;
p=p_0 +  c_p(\x)h^{1+\epsilon}\;,\quad \wb{\epsilon >0,}
\end{equation}
with $\pmb{\mathsf{v}}_0$ and $p_0$ constant.
Density, temperature, and the marker function $\chi$ can vary, and in particular \wb{they} can even be discontinuous, but bounded.  
We prove the following result, which is a generalization of the classical result by~\cite{abgrall1996prevent} (see also e.g.~\cite{abgrall2003efficient}).

\begin{proposition}[Contact discontinuities and mass/momentum/energy conservation  defect\label{prop.contact}]
Consider a quasi one dimensional flow  verifying \eqref{eqn.contact.q1d} in a ball  $B_{i} $      centered around the center of the cell 
$\omega_i$ and such that $\x$ in $\omega_i$ or $\partial\omega_i$ is contained in $B_i$. 
If $N$ denotes the degree of the polynomial approximation used for  $\tilde \v_h$, assume that the numerical quadrature is exact for polynomials
of degree $N-1$ within each cell, and exact for polynomials of the degree $N$ on cell boundaries.
Then,  provided that the update of the non-conservative scheme \eqref{eqn.DGscheme} 
preserves the scaling \eqref{eqn.contact.q1d} and  under the hypotheses of Proposition \ref{prop:quasi-cons},  
  there exist    numerical fluxes  verifying the classical
consistency and continuity assumptions such that the non-conservative scheme \eqref{eqn.DGscheme}  is
\begin{enumerate}
\item quasi-conservative for mass/momentum/energy in the sense of Proposition \ref{prop:quasi-cons},  for single material flows ($\chi=\text{const}$);
% and  the conservative ADER-DG corrector \eqref{eqn.DGscheme-cons} with locally mismatching fluxes are
\item quasi-conservative  for mass/momentum/energy in the sense of Proposition \ref{prop:quasi-cons},  for multi-material flows when
using approximation \eqref{eqn.eos1a} 
\item   quasi-conservative for mass/momentum/energy in the sense of Proposition \ref{prop:quasi-cons},  for multi-material flows when
using approximation \eqref{eqn.eos1} provided the following additional discrete equivalence is true $\forall\x \in B_i$
\begin{equation}\label{eqn.contact.gammaHP}
 \dfrac{1}{\gamma(\chi_h)-1} -    \dfrac{\rho(\chi_h)}{\gamma - 1} = \mathcal{O}(h^{1+\delta})\;,\quad\delta >0
\end{equation}
with  $\rho(\chi_h)$ the linear map allowing to pass from \eqref{eqn.eos1a}  to \eqref{eqn.eos1}: 
$$
\rho = \dfrac{\gamma-1}{\gamma_1-1} + (  \gamma-1)  \dfrac{\chi_h-\chi_1}{\chi_2-\chi_1}\bigg( \dfrac{1}{\gamma_2-1} -\dfrac{1}{\gamma_1-1} \bigg).
$$
\end{enumerate}
\begin{proof}
The proof is obtained by  evaluation of the conservation defect \eqref{eqn.defect} for the mass/momentum/energy equation under the hypotheses made. 
Details are reported in \ref{app:prop-contact-proof}. 
\end{proof}
\end{proposition}

The above result implies that the weak form of the conservation equation is verified on fine enough meshes, and in particular, 
if mesh convergence is achieved, then the resulting solution verifies the appropriate jump conditions for mass, momentum, and  energy.
The  first and second points are well known. The second is essentially a reformulation in our setting of the analysis 
made in~\cite{abgrall1996prevent,abgrall2003efficient}
showing that if one solves an advection equation for $1/(\tilde \gamma-1)$ then mass, momentum and energy jump conditions 
are met across material interfaces.
The last result  requires postulating that the two different approximations of $\tilde\lambda$ obtained with \eqref{eqn.eos1a}  and \eqref{eqn.eos1}
given the same discrete approximation of $\chi$ converge to the same  thing as $h\rightarrow 0$ at a rate sufficient for the 
conditions of  Proposition \ref{prop:quasi-cons}  to be true.
In practice, the numerical results obtained using  \eqref{eqn.eos1a} will show a correct resolution for contacts as well as  
 shock-contact interactions confirming the above result.

\subsection{A hybrid scheme based on \textit{a posteriori} subcell conservative corrections}

Our purpose is now to exploit the \textit{a posteriori} limiting framework to remove completely the conservation defect $\Delta_i^n$ in cells
which may fall outside the heuristics of the previous sections. This essentially means removing  the conservation defect in cell potentially
containing genuinely non-linear discontinuities, namely shock waves.   We also need to   extend this correction
to a region sufficiently large such that the mismatching treatment at cell boundaries may fall in the hypotheses of  Property \ref{prop.smooth.mismatch}. 

To summarize our approach, we will introduce first some heuristics to appropriately flag cells potentially corresponding to a shock wave,
to be combined with those already presented in Section~\ref{ssec.limiter}. 
Next, we want to introduce a correction method which completely removes the conservation defect in these cells 
so that wherever instead the non-conservative update is kept
the results of the previous sections, and in particular Proposition~\ref{prop:quasi-cons}, apply. 
In this respect, the \textit{a posteriori} limiting used in this work provides a {suitable} context
since it relies on a new modified update for the averages in troubled cells, which is precisely the type of quantity we want to control.
How this is achieved is described hereafter.

\paragraph{Modified troubled cell detection} Here, we proceed initially exactly as described in Section~\ref{ssec.limiter}. 
First, we need to introduce the sub-triangulation of each $\omega_i$ composed by the sub-triangles  $\{\tau_{i, \alpha}\}_{\alpha \in [1, N_{\tau, i}]}$.
On each subcell we compute the averages of the non-conservative polynomials as in \eqref{eqn.vh}, for both the solution at time $t^n$,
and for the non-conservative corrector obtained by~\eqref{eqn.DGscheme}. 

To construct them, we start by subdividing each polygonal cell in $N_{\tau, i}=N_{k,i}$ sub-triangles (indeed the same used for 
integration purposes) by connecting the polygon's vertices $\x_{k,i}$ with the cell barycenter~$\x_{{b}_i}$. Then, we further 
subdivide these triangles in $N^2$ 
{small sub-triangles} reaching a total number of $N_{\tau, i}=N_{k,i} \cdot N^2$ subcells which thus grows with the order of 
accuracy $N$ of the method in order to always ensure a proper subcell resolution. As before, the latter averages are 
{denoted with} $\mathbf{v}_{i,\alpha}^{n+1,*}(\x,t^{n+1})$,
which are to be considered as  \textit{candidate} values. The {troubled cell detection} now proceeds as follows:
\begin{enumerate}
\item Evaluate the physical admissibility of the solution, 
namely that density and pressure are strictly positive (we verify in particular that $p,\,\rho \ge \epsilon=10^{-12}$). 
If this criterion is violated the cell is marked as troubled.
\item Evaluate the validity of the relaxed discrete maximum principle \eqref{eqn.RDMP}-\eqref{eqn.deltaRDMP}. If this 
criterion is violated the  cell is also marked as troubled.
\item  Evaluate a \textit{shock detector} based on the approximated undivided divergence 
\begin{equation}
	\widetilde{\nabla \cdot \mathbf{v}}|_{i,\alpha} = \frac{h_{i,\alpha}}{|\tau_{i,\alpha}|}\sum \limits_{f \in \partial \tau_{i,\alpha} 
	}{ |\partial \tau_{i,\alpha,f}|\left(\mathbf{v}^+ + \mathbf{v}^- \right) \cdot \n_{f} }, \qquad c_{s,\min} = \min \limits_{f\in \partial 
	\tau_{i,\alpha}}{\left(c_{s}^-,c_{s}^+\right)},
	\label{eqn.divV}
\end{equation}
with $h_{i,\alpha}=\sqrt{|\tau_{i,\alpha}|}$  the characteristic size of the sub-triangle. The shock indicator $\sigma$ is defined as   
\begin{equation}
	\sigma_{i,\alpha} = \frac{\widetilde{\nabla \cdot \v}|_{i,\alpha} + m_1 h_{i,\alpha} c_{s,\min}}{m_1 h_{i,\alpha} c_{s,\min}},
	\label{eqn.flattener}
\end{equation}
where we employ  $m_1 \in [1;15]$. If $\sigma_{i,\alpha}<0$ then the subcell is marked as a \textit{shock-triggered troubled cell}.
\end{enumerate}
The last criterion  was originally proposed to design positivity preserving schemes for hydrodynamics and magneto-hydrodynamics~\cite{BalsaraFlattener}
and has similarities with other {shock detection criteria} proposed e.g.~\cite{DUCROS1999517,YU18,FUJIMOTO2019264} and references therein.

Whenever a cell or one of its subcells are marked as troubled or shock-triggered troubled, 
the candidate solution $\mathbf{v}^{n+1,*}_{h}$  is rejected, and a 
new solution $\mathbf{v}^{n+1}_{h}$ is recomputed using a limited second order ADER-FV scheme 
based on a predictor-corrector procedure, see~\eqref{eqn.local.predictor} and~\eqref{eqn.subcell-corrector}.

For what concerns the predictor step~\ref{eqn.local.predictor}, it is computed in the same way on any kind of troubled cells,
so we refer to Section~\ref{ssec.limiter} for further details.

\paragraph{Locally conservative finite volume correction} 
The limited ADER-FV corrector step depends instead on the type of troubled cell.
For the cells that do not satisfy {one or both of the first two criteria} described above, 
we apply the corrector formula reported in~\eqref{eqn.subcell-corrector}.
 
On shock-triggered troubled cell, we first set $\q_{i,\alpha}=\Q(\mathbf{v}_{i,\alpha})$, and then 
we apply a conservative update as 
\begin{equation}
\label{eqn.conslimier}
	\mathbf{q}_{i,\alpha}^{n+1} = \mathbf{q}_{i,\alpha}^{n} - \frac{\Delta t}{|\tau_{i,\alpha}|} \sum_{k=1}^3  
	|\tau_{i,\alpha,k}|\, \mathbf{F}\left(\mathbf{v}_{i,\alpha}(x_{i,\alpha,k},\, y_{i,\alpha,k},\, t^{n+1/2}),\, 
	\mathbf{v}_{j,\beta}(x_{i,\alpha,k},\, y_{i,\alpha,k},\, t^{n+1/2})\right)\cdot\mathbf{n}
\end{equation}
with $\mathbf{F}\left(\mathbf{q}^{-}, \mathbf{q}^{+}\right)\cdot\mathbf{n}$ the numerical flux obtained from a
standard two-state Riemann solver. 
% applied to the inner and outer states $\mathbf{q}^-$ and $\mathbf{q}^+$ 
%in the outward normal direction $\mathbf{n}$ (see the book of Toro~\cite{torobook}).
In particular, throughout this work we employed the simple Rusanov~\cite{Rusanov:1961a} approximate flux.

Once the conservative update is performed,  we obtain the primitive subcell average values  
as $\mathbf{v}_{i,\alpha}^{n+1}=\V(\mathbf{q}_{i,\alpha}^{n+1})$.

At this point we proceed as in Section~\ref{ssec.limiter} to reconstruct a primitive variable 
cell polynomial using the reconstruction operator \eqref{eqn.intRec}.
Note that by construction, the subcell averages of the obtained polynomial coincide with the one 
of $\mathbf{v}_{i,\alpha}^{n+1}$. 
As such the conservative average $\Q( \mathbf{v}_{i,\alpha}^{n+1})$ is given by $\mathbf{q}_{i,\alpha}^{n+1} $ 
obtained via the conservative update~\eqref{eqn.conslimier}, which is by construction locally conservative. 

The extension of the flagging in space and time follows the  
indications provided at the end of Section~\ref{ssec.limiter}. 
In particular, in the interest of safeguarding the overall robustness and making sure the unflagging is sufficiently far from shocks, 
the troubled cell flag is extended also to all the direct neighbors of the elements flagged as troubled by the admissibility 
criteria. In particular neighbors of cells where shocks have been detected, share the same 
treatment with regards to the type of limiter scheme being applied to them.

%--------- END OF SECTION -------------------------------------------------	

%--------- SECTION --------------------------------------------------------
\section{Numerical results} \label{sec_validation}	

In this section, we present some numerical results in order both 
to validate our quasi-conservative approach on well-established benchmarks in gasdynamics, 
and illustrate its capabilities and potential applicability in the field of multi-material modeling. 
In particular, we first verify the order of convergence and the shock capturing properties of our strategy
by comparing our numerical results with exact or reference solutions; 
then, we show that our novel approach completely avoids spurious oscillation 
and easy the computation of primary thermodinamical variables, like pressure, 
at the (smeared) interface between different materials.
All our numerical tests are carried out on unstructured two-dimensional 
polygonal tessellations and the CFL stability value in~\eqref{eq.CFL} is taken to be $0.1$ for the 
first $10$ time steps of each simulations (as suggested in~\cite{torobook}) and then is always chosen CFL~$ = 0.5$. 

Finally, we remark that for all the test cases presented in this Section for what concerns the shock indicator 
of~\eqref{eqn.flattener} we have set $m_1=1$ for the classical single phase perfect gas Euler equations, 
while we have $m_1=15$ when using the multi-material EOS with spatially variable parameters.

\subsection{Euler equations with constant EOS parameter}

\subsubsection{Moving isentropic Shu-type vortex}
\begin{table*}[!t]
	\caption{Isentropic Shu-type vortex. Numerical convergence results for our ADER predictor-corrector 
	scheme applied on the non-conservative formulation of the Euler equations from the second order \Puno to the fourth order \Ptre scheme. 
	The error norms refer to the primitive variables $\rho$ and $p$ at time $t=1.0$ in the $L_2$ norm.} 
	\label{tab.orderOfconvergenceDG_shu}
	\centering
	\begin{tabular}{ccccc|ccccc|ccccc} 
		\hline 
		\multicolumn{5}{c}{\Puno$\rightarrow \mathcal{O}2$} & \multicolumn{5}{c}{\Pdue$\rightarrow \mathcal{O}3$} & \multicolumn{5}{c}{\Ptre$\rightarrow \mathcal{O}4$}    \\ 
		\hline
		\!\!$h$ \!\!\!&\!\!\! $\epsilon_{L_2}(\rho)$ \!\!\!&\!\!\!\! $\!\!\mathcal{O}(\rho)\!\!\!$ \!\!\!&\!\!\! $\epsilon_{L_2}(p)$ \!\!\!&\!\!\!\! $\!\!\mathcal{O}(p)\!\!\!$ \!\!\!&  $h$ \!\!\!&\!\!\! $\epsilon_{L_2}(\rho)$ \!\!\!&\!\!\! $\!\!\mathcal{O}(\rho)\!\!\!$ \!\!\!&\!\!\! $\epsilon_{L_2}(p)$ \!\!\!&\!\!\!\!\! $\!\!\!\mathcal{O}(p)\!\!\!$ \!\!\!& $h$ \!\!\!&\!\!\! $\epsilon_{L_2}(\rho)$ \!\!\!&\!\!\!\! $\!\!\mathcal{O}(\rho)\!\!\!$ \!\!\!&\!\!\! $\epsilon_{L_2}(p)$ \!\!\!&\!\!\!\! $\!\!\mathcal{O}(p)\!\!\!\!$  \\
		\hline
		\!\!\!1.15    \!\!\!&\!\!\! 1.28E-1  \!\!\!&\!\! -   \!\!\!&\!\!\! 1.59E-1  \!\!\!&\!\! -   \!\!&\!\!1.44    \!\!\!&\!\!\! 4.43E-2 \!\!\!&\!\! -   \!\!\!&\!\!\! 5.70E-2 \!\!\!&\!\! -    \!\!&\!\!1.93 \!\!\!&\!\!\!   3.55E-2 \!\!\!&\!\!  -  \!\!\!&\!\!\!   4.10E-2 \!\!\!&\!\!   - \!\!\!\\ 
		\!\!\!7.38E-1 \!\!\!&\!\!\! 4.62E-2  \!\!\!&\!\! 2.3 \!\!\!&\!\!\! 5.33E-2  \!\!\!&\!\! 2.5 \!\!&\!\!9.30E-1 \!\!\!&\!\!\! 1.41E-2 \!\!\!&\!\! 2.6 \!\!\!&\!\!\! 1.71E-2 \!\!\!&\!\! 2.8  \!\!&\!\!1.44 \!\!\!&\!\!\!   1.34E-2 \!\!\!&\!\! 3.3 \!\!\!&\!\!\!   1.58E-2 \!\!\!&\!\!  3.3\!\!\!\\
		\!\!\!4.73E-1 \!\!\!&\!\!\! 1.56E-2  \!\!\!&\!\! 2.4 \!\!\!&\!\!\! 1.83E-2  \!\!\!&\!\! 2.4 \!\!&\!\!5.76E-1 \!\!\!&\!\!\! 3.51E-3 \!\!\!&\!\! 2.9 \!\!\!&\!\!\! 4.38E-3 \!\!\!&\!\! 2.9  \!\!&\!\!1.15 \!\!\!&\!\!\!   6.00E-3 \!\!\!&\!\! 3.6 \!\!\!&\!\!\!   6.75E-3 \!\!\!&\!\!  3.8\!\!\!\\
		\!\!\!3.85E-1 \!\!\!&\!\!\! 9.86E-3  \!\!\!&\!\! 2.2 \!\!\!&\!\!\! 1.17E-2  \!\!\!&\!\! 2.2 \!\!&\!\!4.73E-1 \!\!\!&\!\!\! 2.06E-3 \!\!\!&\!\! 2.7 \!\!\!&\!\!\! 2.58E-3 \!\!\!&\!\! 2.7  \!\!&\!\!0.92 \!\!\!&\!\!\!   2.73E-3 \!\!\!&\!\! 3.8 \!\!\!&\!\!\!   3.16E-3 \!\!\!&\!\!  3.6\!\!\!\\		
		\hline 
	\end{tabular}       
\end{table*}

We start our set of benchmarks by considering a smooth isentropic vortex flow inspired to the one 
proposed in~\cite{HuShuVortex1999}, 
which represents a stationary equilibrium of the Euler system and allows to easily compute the order 
of convergence of our scheme. 
The initial computational domain is the square $\Omega=[-20;30]\times[-20;30]$ with wall boundary conditions. 
The initial condition is given by some perturbations $\delta$ that are superimposed onto a homogeneous background field $\V_0=(\rho,u,v,p)=(1,0,0,1)$. 
The perturbations for density and pressure are
\begin{equation}
	\label{rhopressDelta}
	\delta \rho = (1+\delta T)^{\frac{1}{\gamma-1}}-1, \quad \delta p = (1+\delta T)^{\frac{\gamma}{\gamma-1}}-1, 
\end{equation}
with the temperature fluctuation $\delta T = -\frac{(\gamma-1)\epsilon^2}{8\gamma\pi^2}e^{1-r^2}$ and the 
vortex strength is $\epsilon=5$.
The velocity field is affected by the following perturbations
\begin{equation}
	\label{ShuVortDelta}
	\left(\begin{array}{c} \delta u \\ \delta v  \end{array}\right) = \frac{\epsilon}{2\pi}e^{\frac{1-r^2}{2}} \left(\begin{array}{c} 
	-(y-5) \\ \phantom{-}(x-5)  \end{array}\right).
\end{equation}
Due to the smooth character of this test case the limiter is never needed, 
thus the equations are solved everywhere by employing our ADER predictor-corrector~\eqref{eqn.ADERpred}-\eqref{eqn.DGscheme} 
scheme applied on the non-conservative formulation. The obtained order of convergence is 
reported in Table~\ref{tab.orderOfconvergenceDG_shu} and it is achieved as expected.

\subsubsection{Planar shock}

In this subsection and the next three, we run some essentially 1D Riemann problems 
to assess the ability of our method to capture one-dimensional simple waves.

We start with some planar right moving shock waves respectively traveling at {Mach} $M=5, 10$ and $20$ 
on the domain $\Omega = [0,1.2]\times[0,0.1]$ with initial conditions such that
\begin{equation}
(\rho, u, v, p)(\x) = 
\begin{cases} 
	(\gamma,   0, 0, 1       ),  \ \text{ if } \  x > 0.2, \\
	(\gamma \delta \rho,   \delta u , 0, \delta p    ),  \ \text{ if } \  x \le 0.2, \\	
\end{cases}
\end{equation}
with
$\delta \rho = \frac{M^2(\gamma + 1 )}{ 2 + M^2( \gamma - 1 ) }$,
$\delta p = \frac{ 2 \gamma M^2  - \gamma + 1 }{ \gamma + 1 } $
and 
$\delta u =  
\frac{2( M - \frac{1}{M} )}{ \gamma + 1}\sqrt{ \frac{\gamma\,p}{\rho} } $.

\begin{figure}[!b]
	\centering
	\includegraphics[width=0.5\linewidth,trim=95 5 5 5,clip]{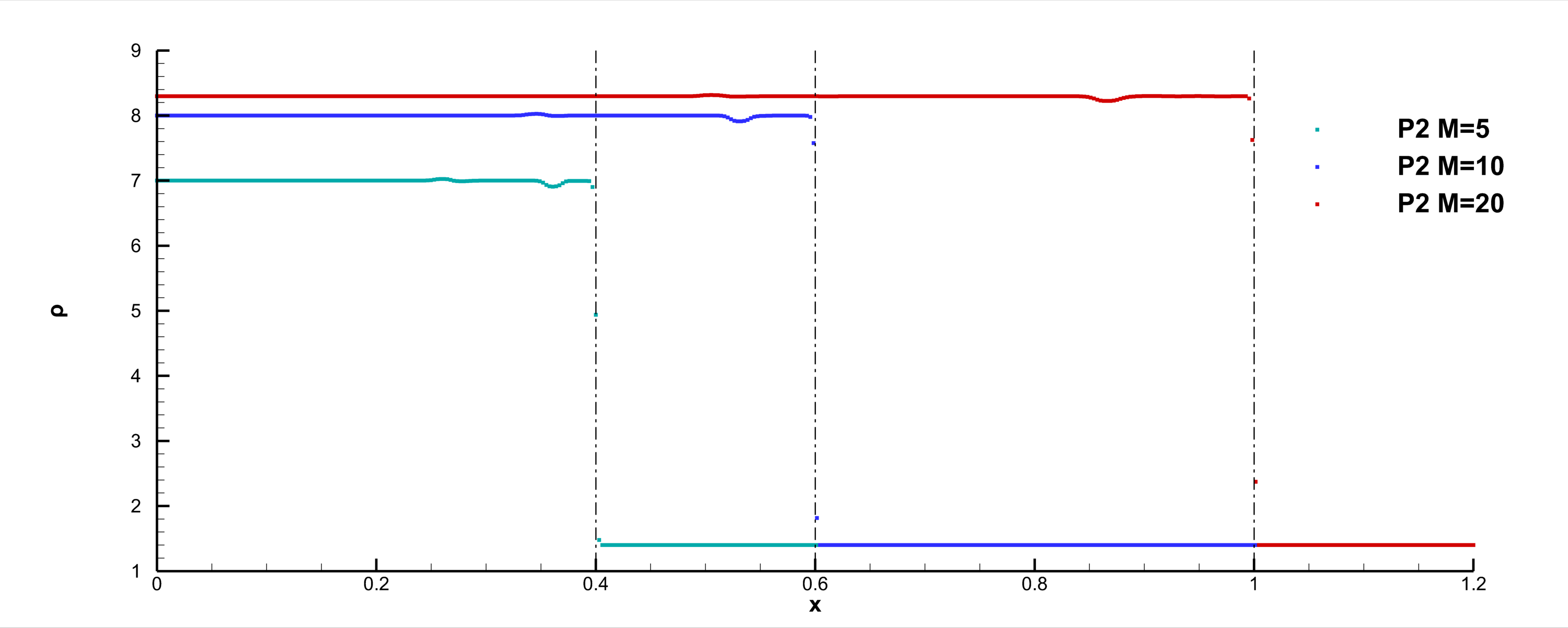}%
	\includegraphics[width=0.5\linewidth,trim=95 5 5 5,clip]{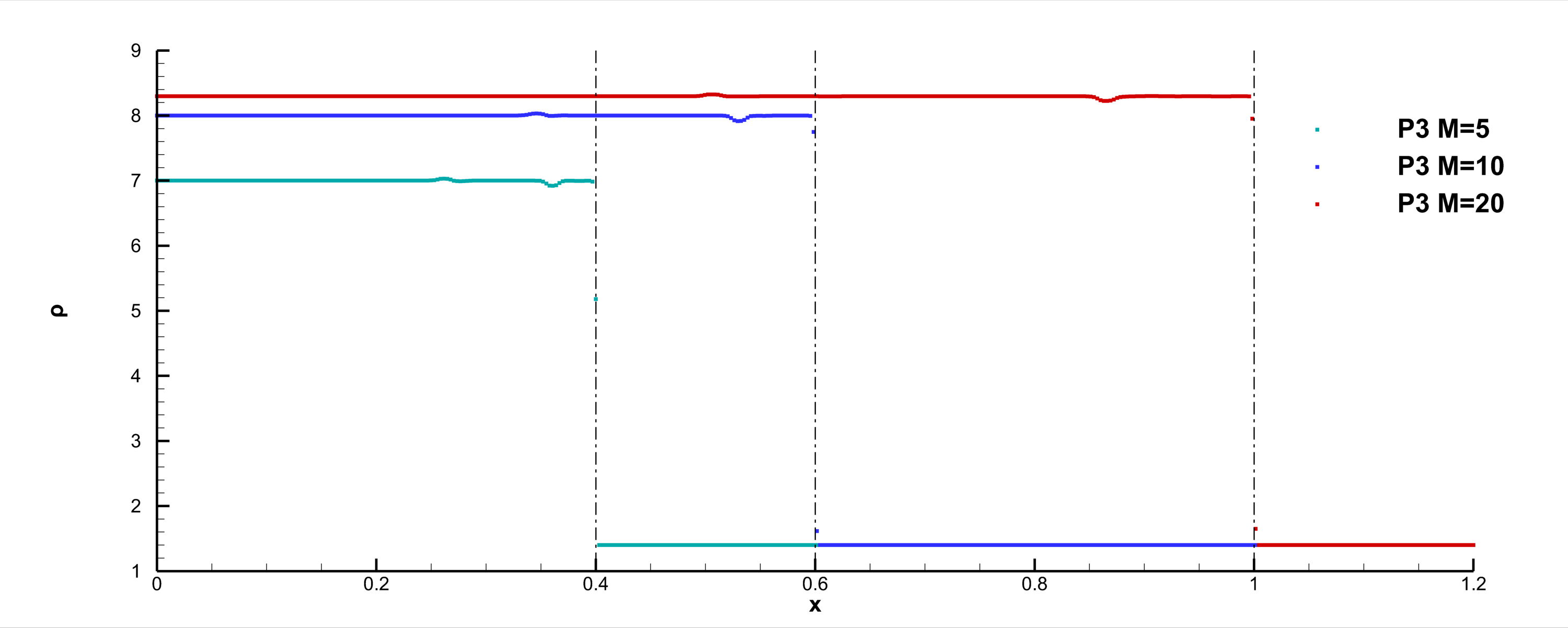}%
	\caption{Numerical results obtained with our third \Pdue (left) and fourth \Ptre (right) order 
	scheme for the planar shock wave propagation with Mach $M=5$~{(dark-green)}, 
	$M=10$~(blue) and $M=20$~(red). We also report with dashed black lines the expected position of 
	each shock wave at the final time $t_f=0.04$ to show that our scheme perfectly deals with shock waves modeling.}
	\label{fig.planarshock}
\end{figure}

The results obtained with our quasi-conservative approach of order three and four at $t_f=0.04$ 
on a polygonal tessellation with characteristic mesh size $h=1/180$ are reported in 
Figure~\ref{fig.planarshock} together with the expected shock position for each propagation speed. % h = 5.46E-3
We note that both  shock position and magnitude are   always perfectly captured by our scheme.
We can also see the appearance of the start-up  perturbations  downstream of the shock which 
are related to the small mismatch between the exact one dimensional
Rankine-Hugoniot conditions, and the  discrete jump conditions verified by the solution on the 
two-dimensional grid. Similar features have been observed and analysed
with fully conservative schemes as well e.g. in~\cite{RICCHIUTO2005249}.

\subsubsection{Lax shock tube}
\label{ssec.lax}

\begin{figure}[!b]
	\centering
	\includegraphics[width=0.63\linewidth]{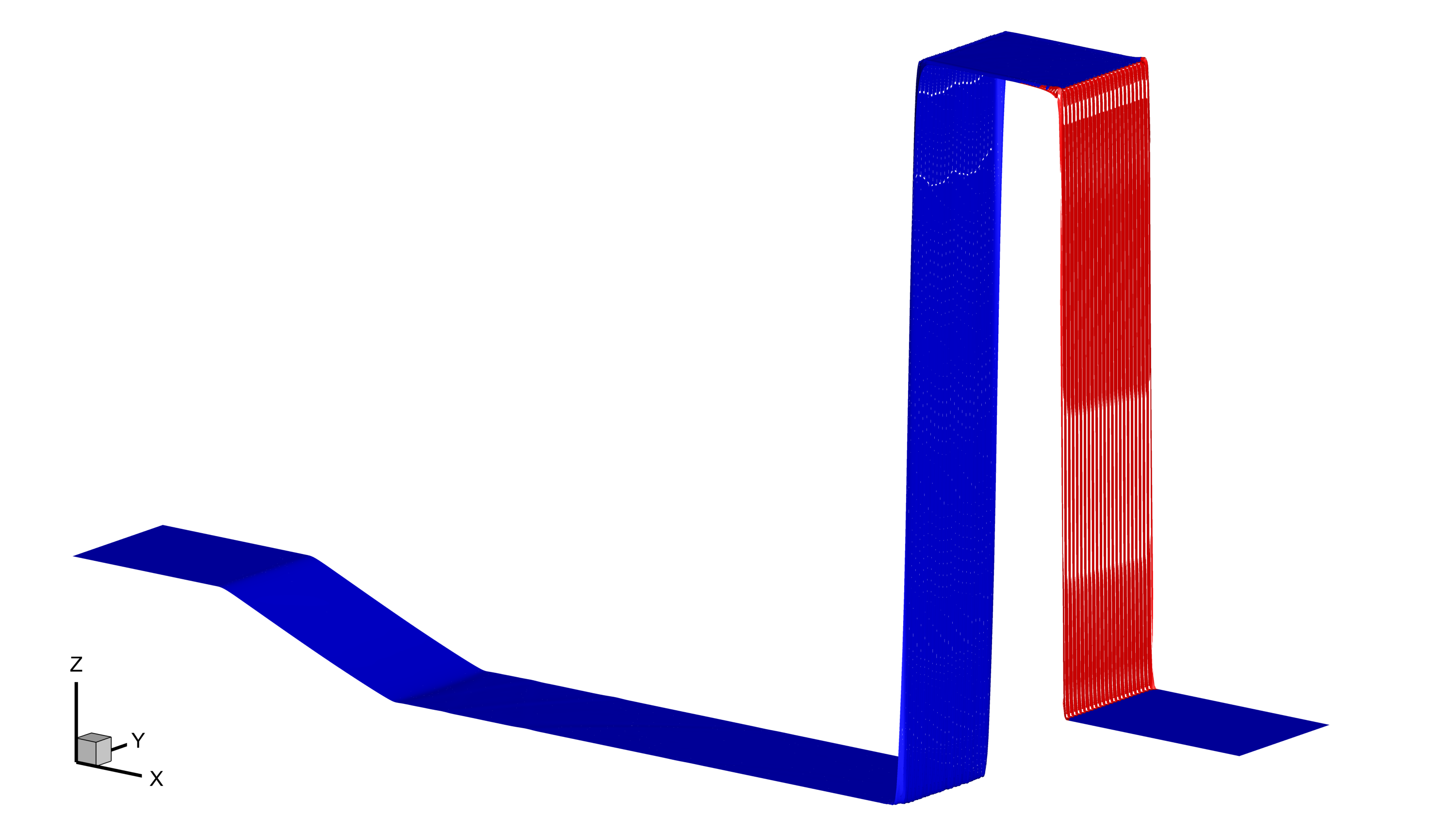}%
	\includegraphics[width=0.37\linewidth]{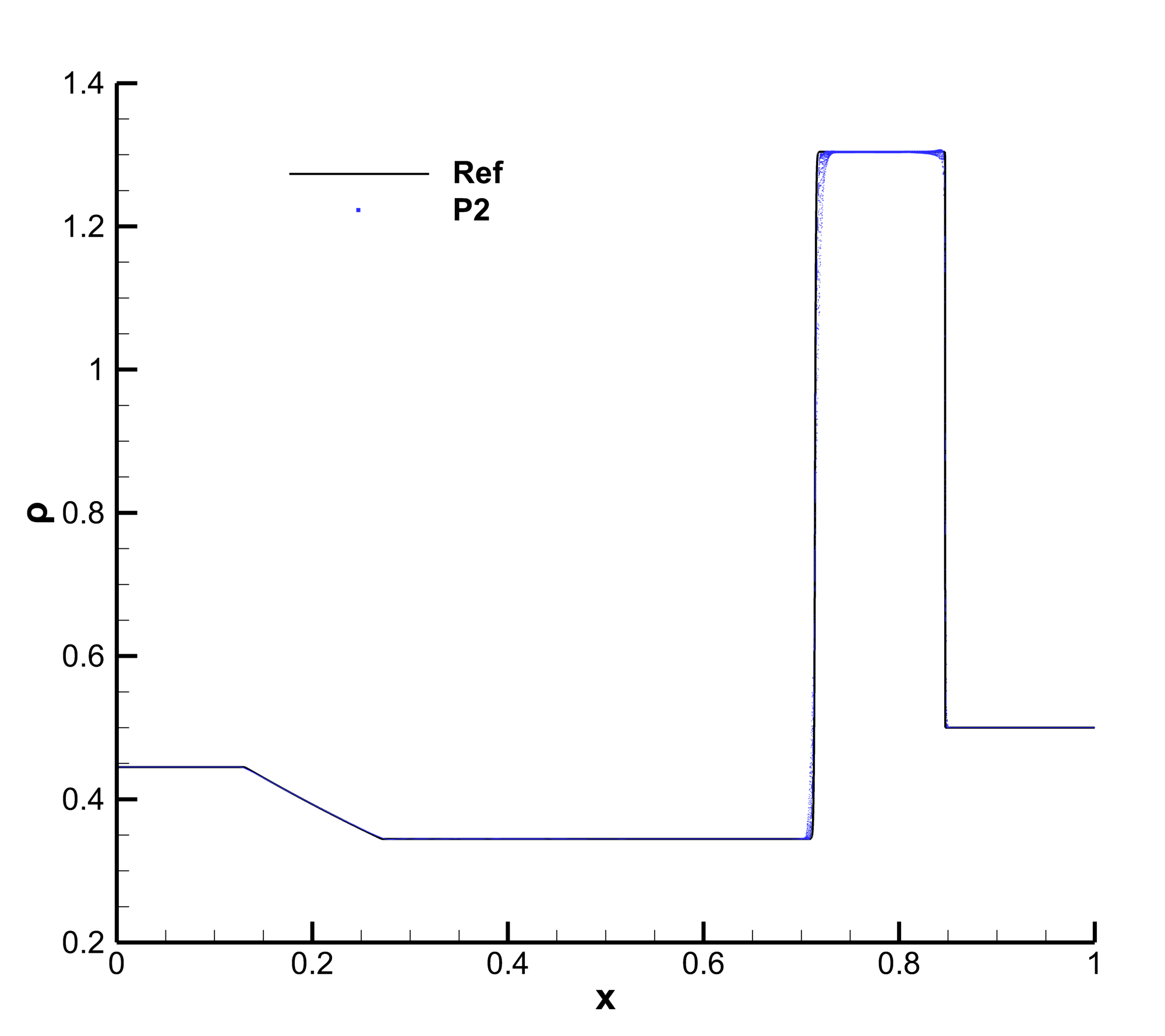}\\
	\includegraphics[width=0.63\linewidth]{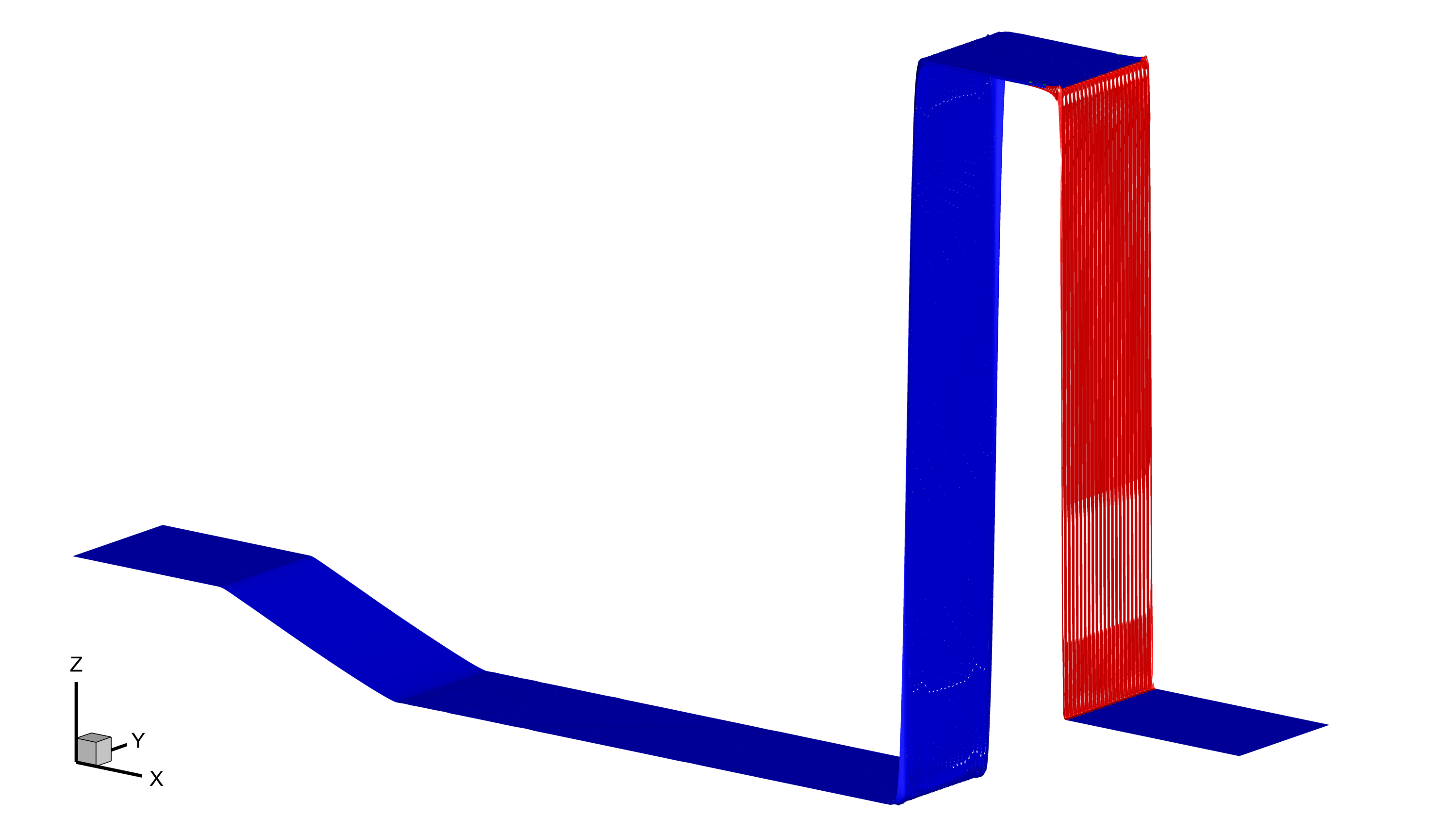}%
	\includegraphics[width=0.37\linewidth]{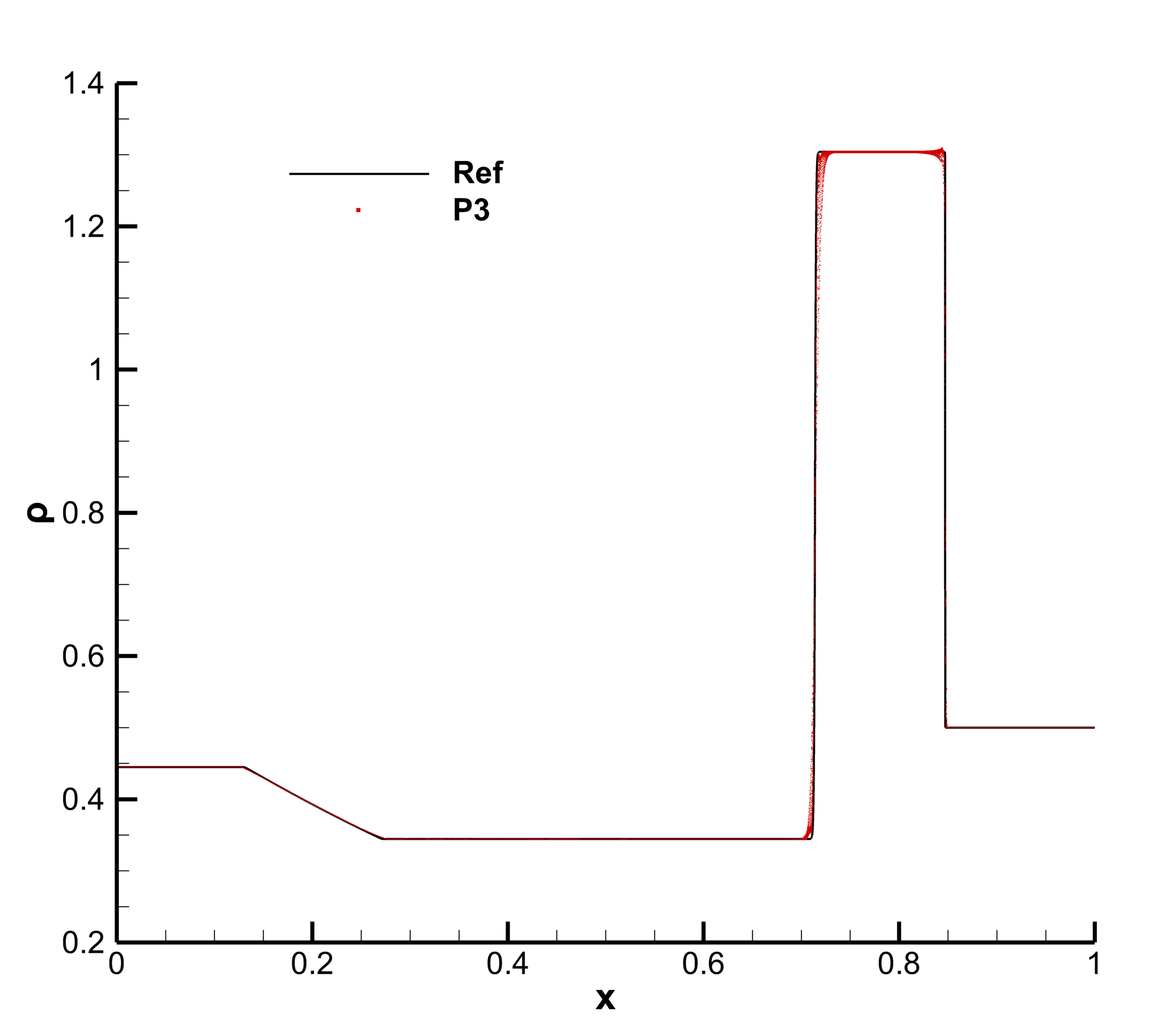}\\
	\caption{Lax shock tube solved on a polygonal tessellation with characteristic mesh size $h=1/430$. %$h=2.3E-3$
		We report the results obtained with our quasi-conservative third order \Pdue (top) and fourth 
		order \Ptre (bottom) schemes. In particular, we show on the left a graph whose $z-$coordinate 
		is given by the density and whose colors refer to the type of scheme employed on each cell: 
		blue for our quasi-conservative scheme, green for the second order limiter in primitive variables 
		and red for the second order limiter applied on the conservative formulation, which is used 
	to discretize the solution where the shock wave is located. 
	Then, on the right, we report a scatter plot of our numerical results for the density profile 
	compared with a reference solution obtained with a one dimensional Euler solver using a fully conservative
	 explicit second order limited Residual Distribution scheme (see e.g.~\cite{RICCHIUTO20105653,AR:17}) on 50000 points.}
	\label{fig.lax}
\end{figure}

We continue our set of Riemann problems by solving the Lax shock tube, originally due to~\cite{lax2}, 
which in addition to one shock wave also develops a contact wave and a rarefaction fun.
The initial conditions that originate the three waves are 
\begin{equation}
	(\rho, u, v, p)(\x) = 
	\begin{cases} 
		(0.5  , 0, 0,  0.571    )   &  \ \text{ if } \  x > 0.5, \\
		(0.445,  0.698, 0, 3.528 )  &   \ \text{ if } \  x \le 0.5, \\	
	\end{cases}
\end{equation}
and we solve this test problem on the domain $\Omega = [0,1]\times[0,0.1]$ up to the final time $t_f =0.14$.
We report our numerical results compared with a reference solution in Figure~\ref{fig.lax}, 
where we can notice that the conservative second order limiter is activated exactly where the shock wave is located 
which allows to perfectly track the shock position without any error on the computation of the shock speed. 

\subsubsection{Circular Sod explosion problem}
\label{ssec.sod}

\begin{figure}[!b]
	\centering
	\includegraphics[width=0.57\linewidth,trim=1 1 1 1,clip]{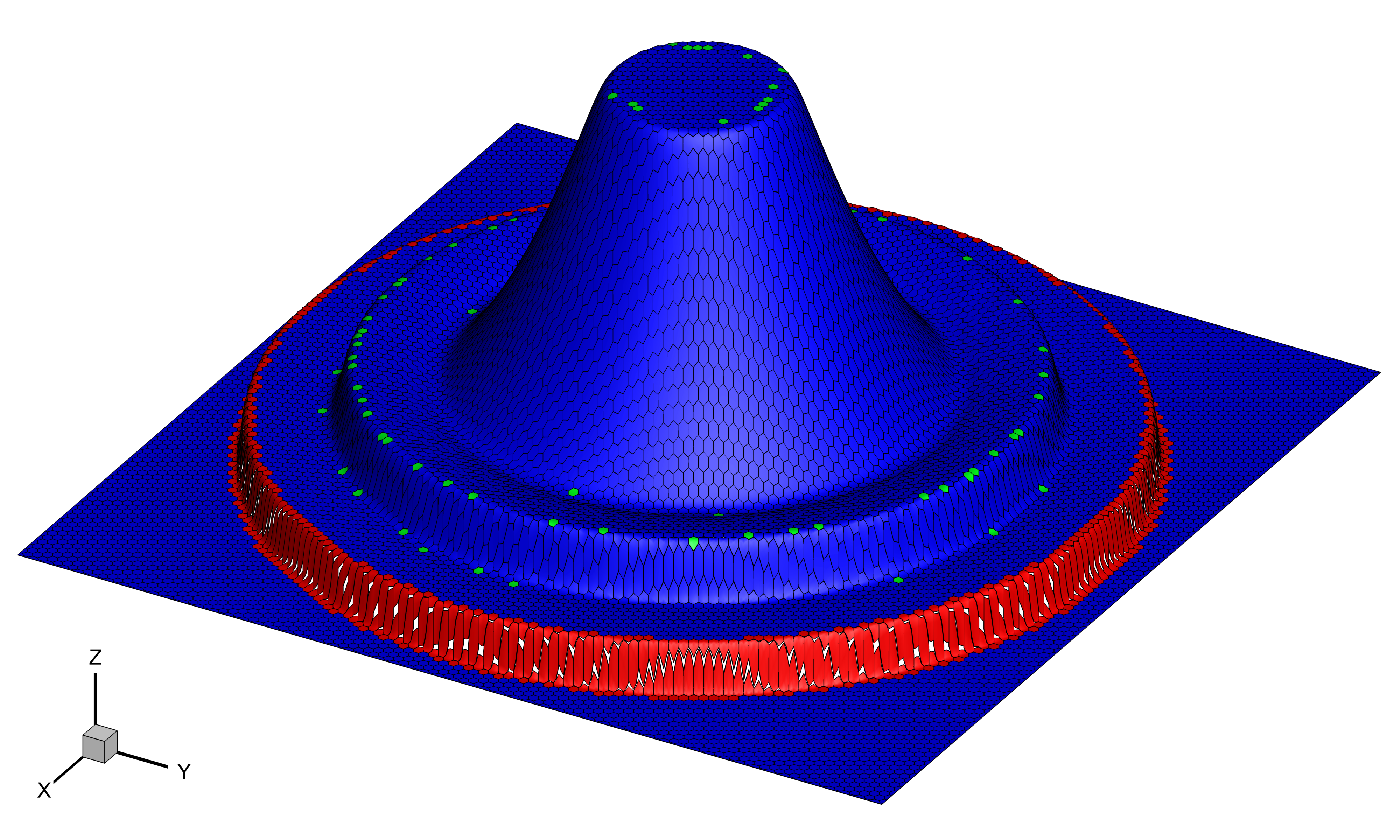}\quad
	\includegraphics[width=0.40\linewidth,trim=1 1 1 1,clip]{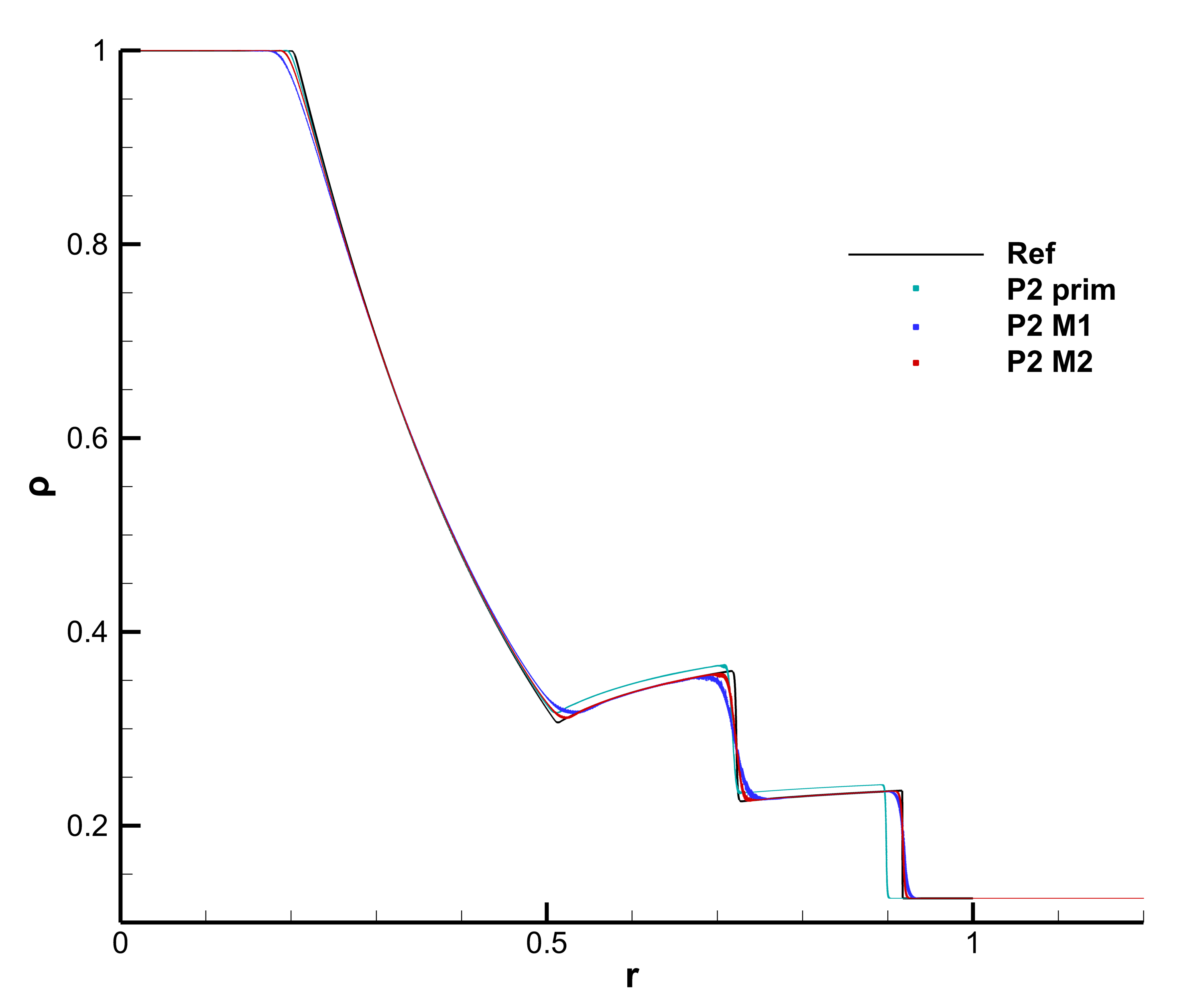}\\
	\includegraphics[width=0.57\linewidth,trim=1 1 1 1,clip]{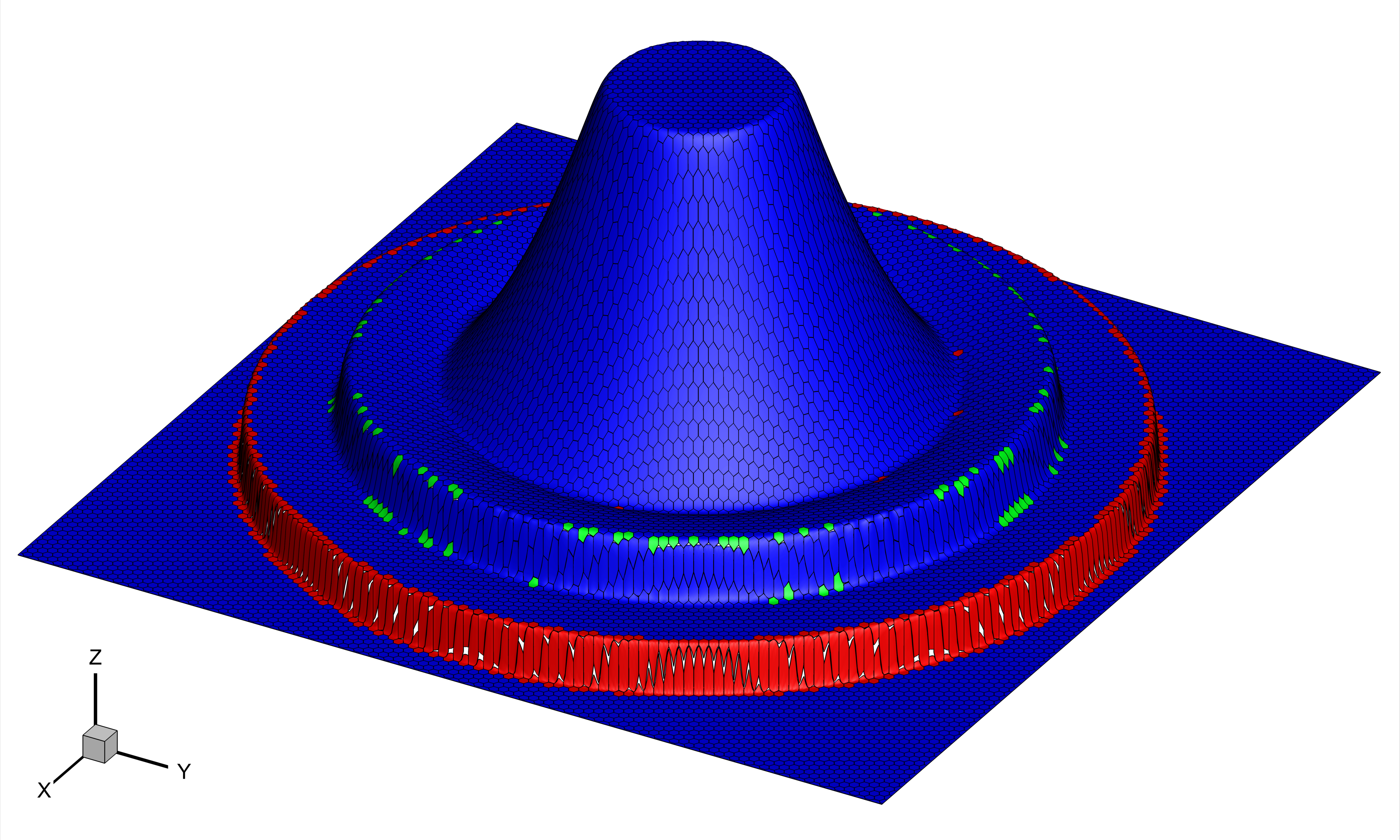}\quad
	\includegraphics[width=0.40\linewidth,trim=1 1 1 1,clip]{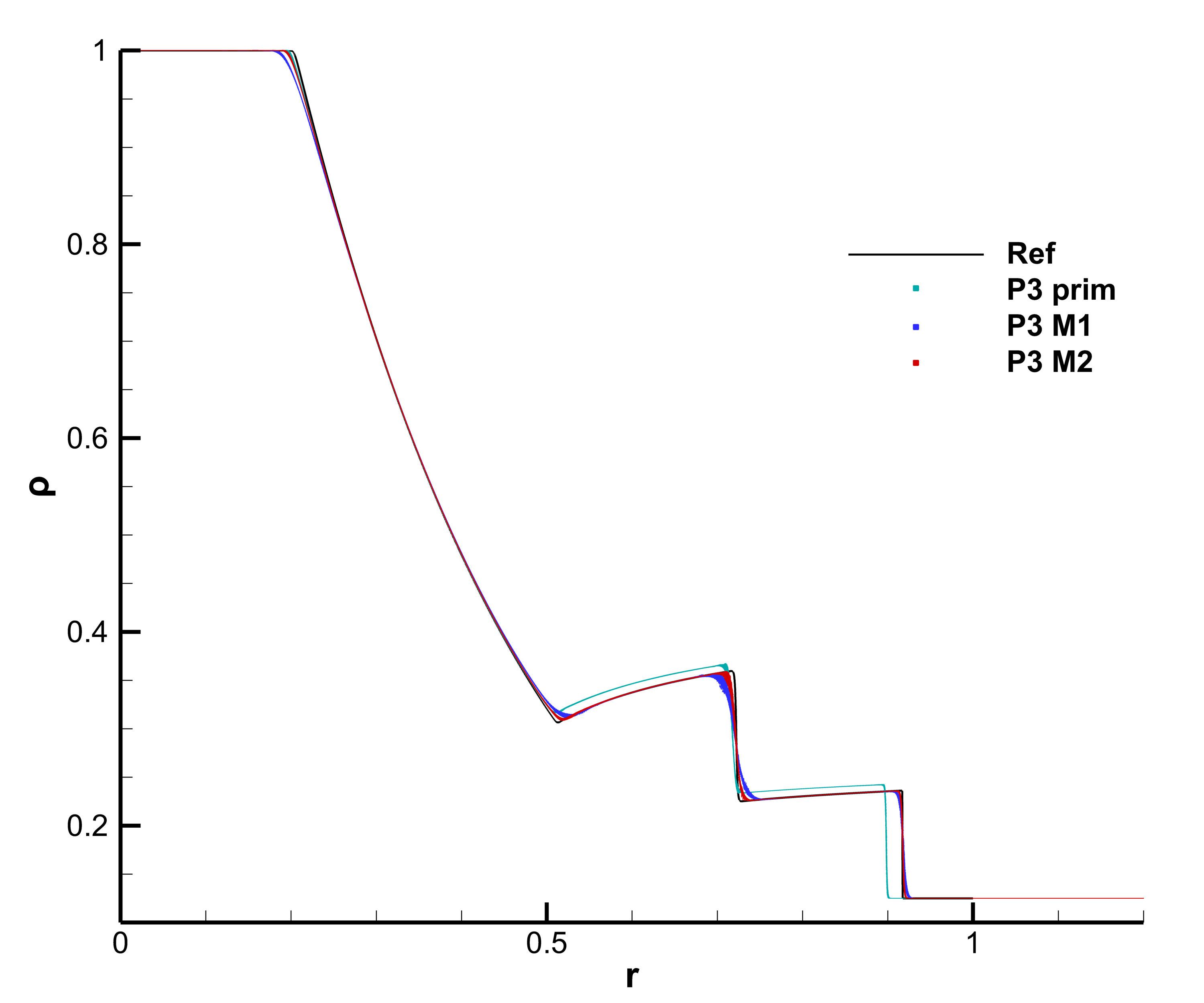}\\
	\caption{Circular Sod explosion. We show the numerical results obtained with our quasi-conservative 
		third order \Pdue (top row) and fourth order \Ptre (bottom row) discontinuous Galerkin schemes. 
		On the left, we depict in blue the cells where our DG scheme has been applied to the non-conservative formulation 
		and in red or green those where our second order subcell FV limiter has been activated. 
		In particular, on the red cells which correspond to \textit{shock-triggered troubled cells} 
		the FV schemes is applied on the conservative formulation of the Euler equations.
		On the right, we report the scatter plot of the numerical solution obtained 
		on a coarse mesh $M_1$ with characteristic mesh size of $h=1/50$ % 1.84E{-2}
		and a finer one $M_2$ with characteristic mesh size of $h=1/215$, %4.64E{-3} 
		compared with a reference solution (black), 
		and the \textit{incorrect} solution ({green}) that one could have obtained 
		by applying the DG scheme for the non-conservative formulation everywhere on the domain, including the shock regions.
	}
	\label{fig.sod}
\end{figure}

This circular explosion problem can be seen as a multidimensional extension of the classical Sod test case.
Here, we consider as computational domain a square of dimension $[-1;1]\times[-1;1]$, 
and the initial condition is composed of two  different states, separated by a discontinuity at radius $r_d=0.5$ 
\begin{equation}
	(\rho, u, v, p)(\x) = 
	\begin{cases} 
		(1.0,  0, 0, 1 )        &  \ \text{ if } \  r \le r_d, \\	
		(0.125, 0, 0,  0.1    ) &  \ \text{ if } \  r > r_d. \\		
	\end{cases}
\end{equation}
The final time is chosen to be $t_f=0.25$, so that the shock wave does not cross 
the external boundary of the domain, where wall boundary conditions are imposed. 
We have obtained a reference solution thanks to the rotational symmetry of the problem 
which reduces to a one-dimensional problem with geometric source terms,
which we have solved by using a classical second order TVD scheme on a very fine mesh. 
In Figure~\ref{fig.sod}, we report the numerical results obtained with our third and fourth order quasi-conservative schemes, 
compared both with a reference solution and the wrong results that one could obtain 
by working with the primitive formulation of the Euler equations everywhere on the domain.
Indeed, with our approach on the shock regions (depicted in red in Figure~\ref{fig.sod}), 
instead of using the original scheme, we apply our second order limiter scheme 
acting on the conservative formulation, 
and it is evident that this local switch is enough to restore the ability of the scheme to correctly
capture the propagation speed of shockwaves.

\subsubsection{Woodward-Colella blast wave}

The last Riemann problem that we consider is the so-called Woodward and Colella blast wave due to~\cite{woodward1984numerical}.
This test case is a standard low energy benchmark problem and it is characterized by the interaction of two strong shocks 
that reach the boundaries and are further reflected so to produce additional mutual interactions. 
The initial data is built as follow 
\begin{equation}
	(\rho, u, v, p)(\x) = 
	\begin{cases} 
		(1, 0, 0, 1000      ) &  \ \text{ if } \  0 \le x \le 0.1, \\	
		(1, 0, 0, 10^{-2}   ) &  \ \text{ if } \  0.1 < x \le 0.9,		 \\		
		(1, 0, 0, 10^2      ) &  \ \text{ if } \  0.9 < x \le 1.0,		 \\		
	\end{cases}
\end{equation}
and the problem is solved on the rectangular domain $\Omega = [0,1]\times[0,0.1]$
	
\begin{figure}[!b]
	\centering
	\includegraphics[width=0.7\linewidth]{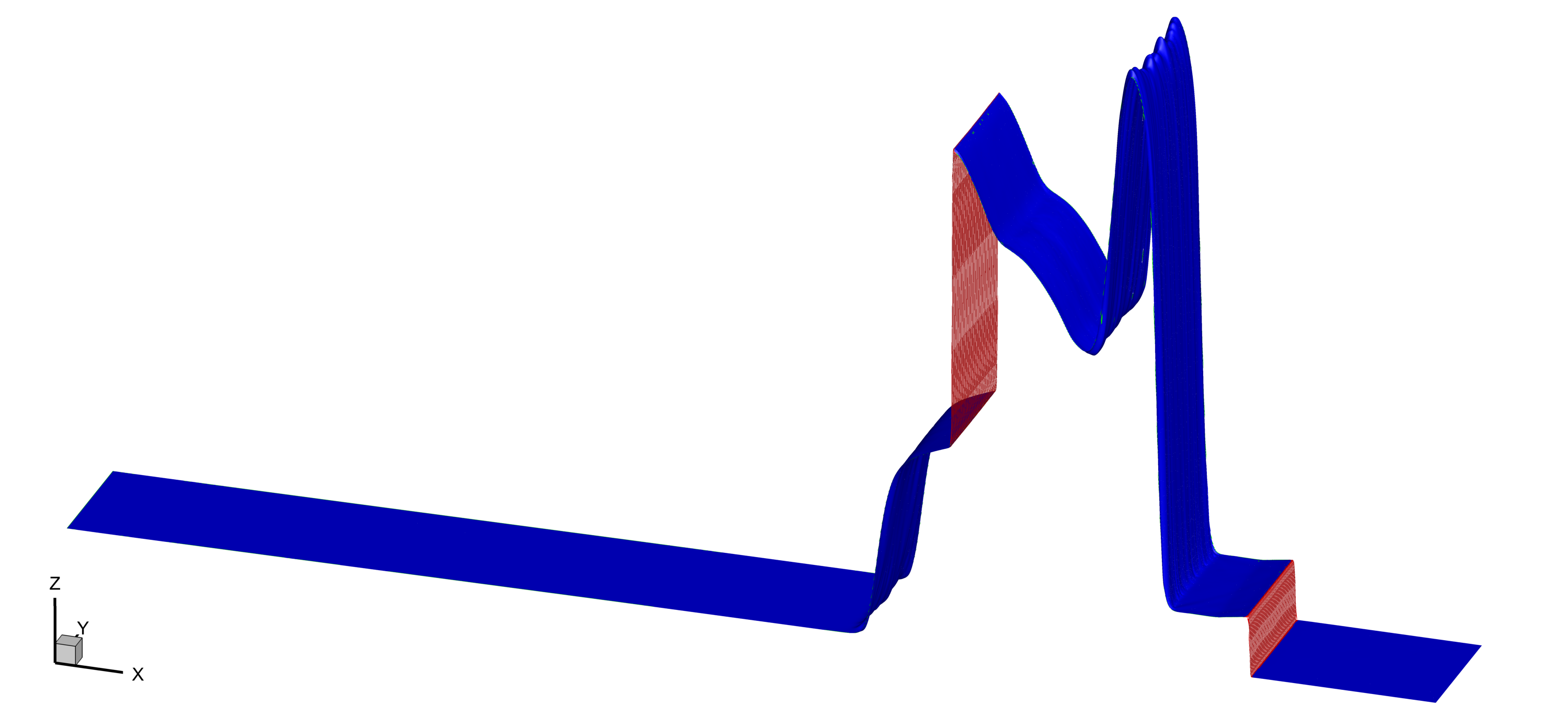}%
	\includegraphics[width=0.3\linewidth]{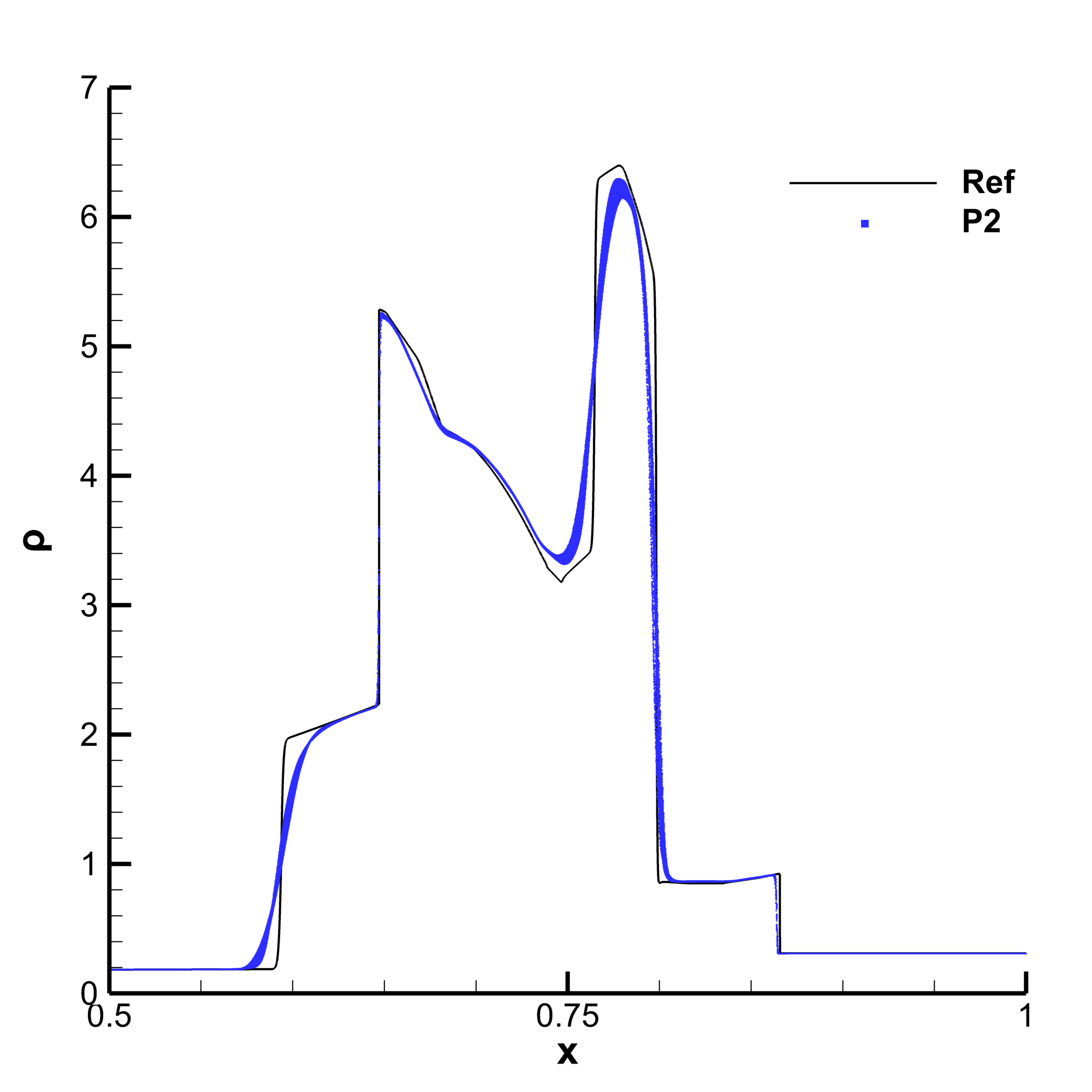}\\
	\includegraphics[width=0.7\linewidth]{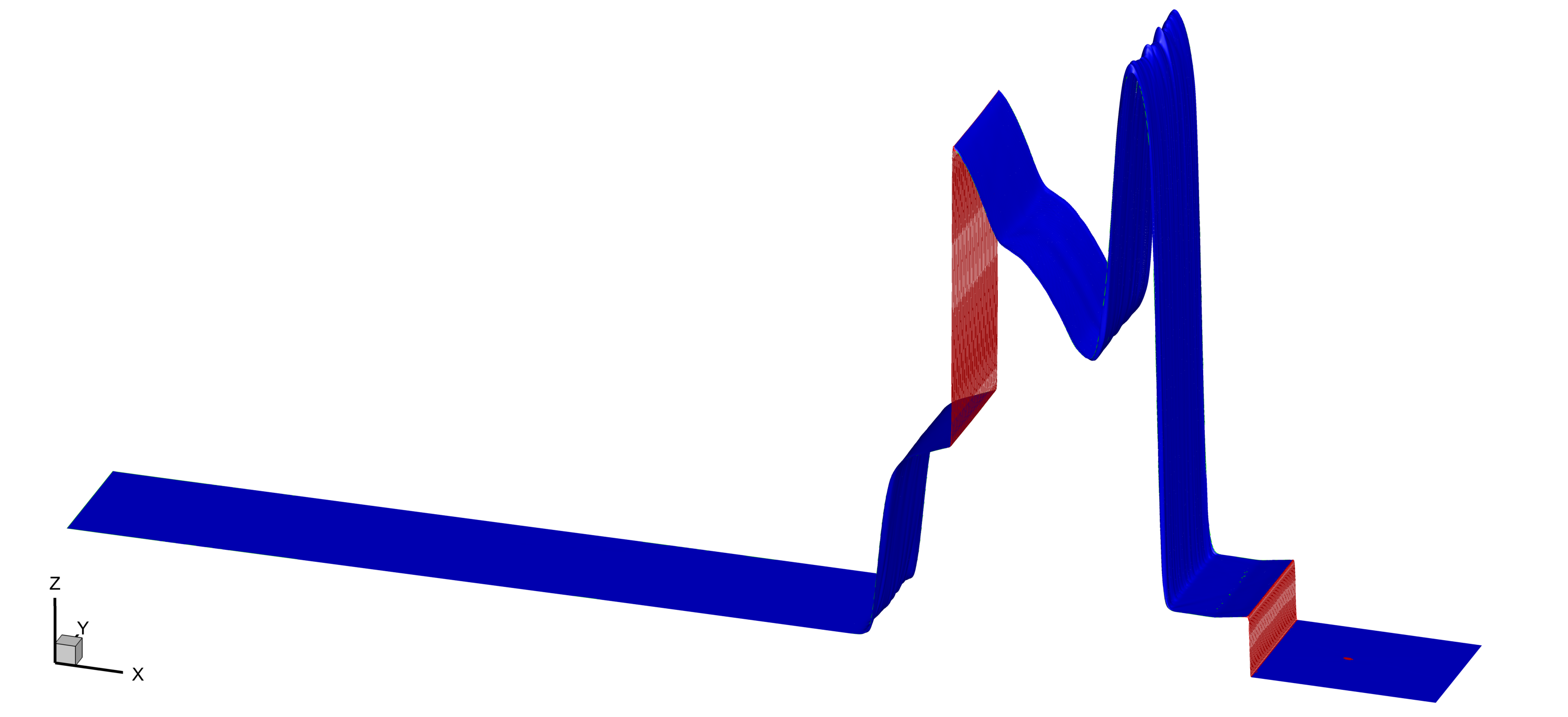}%
	\includegraphics[width=0.3\linewidth]{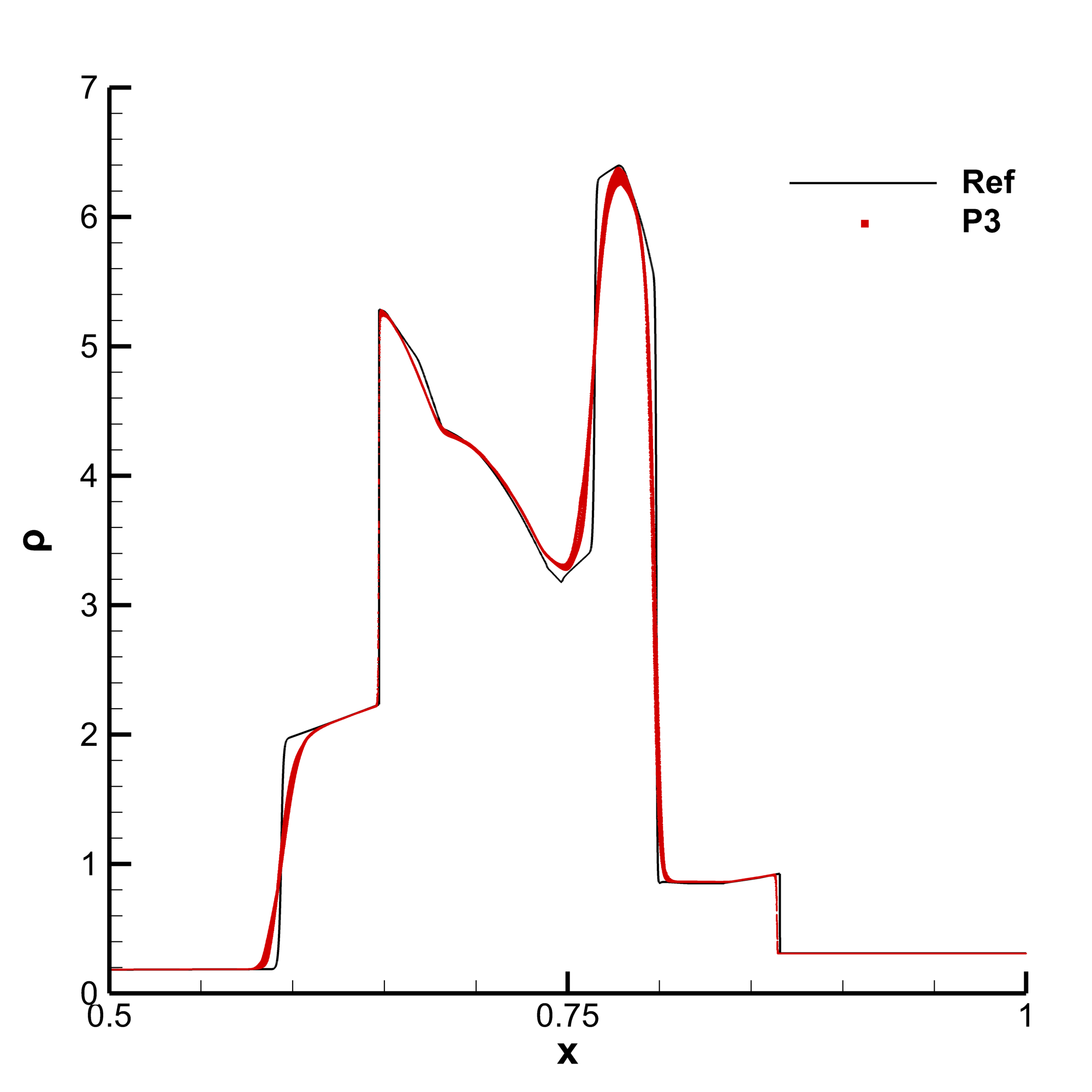}\\
	\caption{Woodward-Colella blast waves at the final time $t= 0.038$ solved on a polygonal tessellation 
	with characteristic mesh size $h=1/1080$. %9.28E-4
	We report the results obtained with our third order \Pdue (top row) and fourth order \Ptre (bottom row) 
	discontinuous Galerkin schemes, 
	which solve the Euler equations formulated in primitive variables everywhere 
	a part for the shock regions, highlighted in red on the left, where a second order FV scheme has been 
	applied to the conservative formulation of the PDE at the subgrid level.
	In particular, on the right we report a zoom on the region $x\in[0.5,1]$ where 
	we compare the scatter plot of our numerical results with a reference solution obtained 
	with a fully conservative non-linear Residual Distribution scheme~\cite{RICCHIUTO20105653,AR:17} 
	applied to the 1D  Euler equations on 50000 points. 
}
	\label{fig.woodwardcolella}
\end{figure}

We report our numerical results in Figure~\eqref{fig.woodwardcolella}: 
in particular, we depict in red the \textit{shock-triggered troubled cells} and we compare 
the scatter plot of the solution obtained with our quasi-conservative schemes with a reference solution obtained with
a fully conservative one dimensional explicit second order limited Residual Distribution scheme (see e.g.~\cite{RICCHIUTO20105653,AR:17}) on 50000 points.
In addition to the correct detection of all the generated waves and the correct location of shocks, 
 convergence can be noticed by comparing visually the results obtained with our third and fourth order schemes.

\subsubsection{Forward Facing Step}

% big verision
\begin{figure}[!b]
	\centering
	\includegraphics[width=0.88\linewidth]{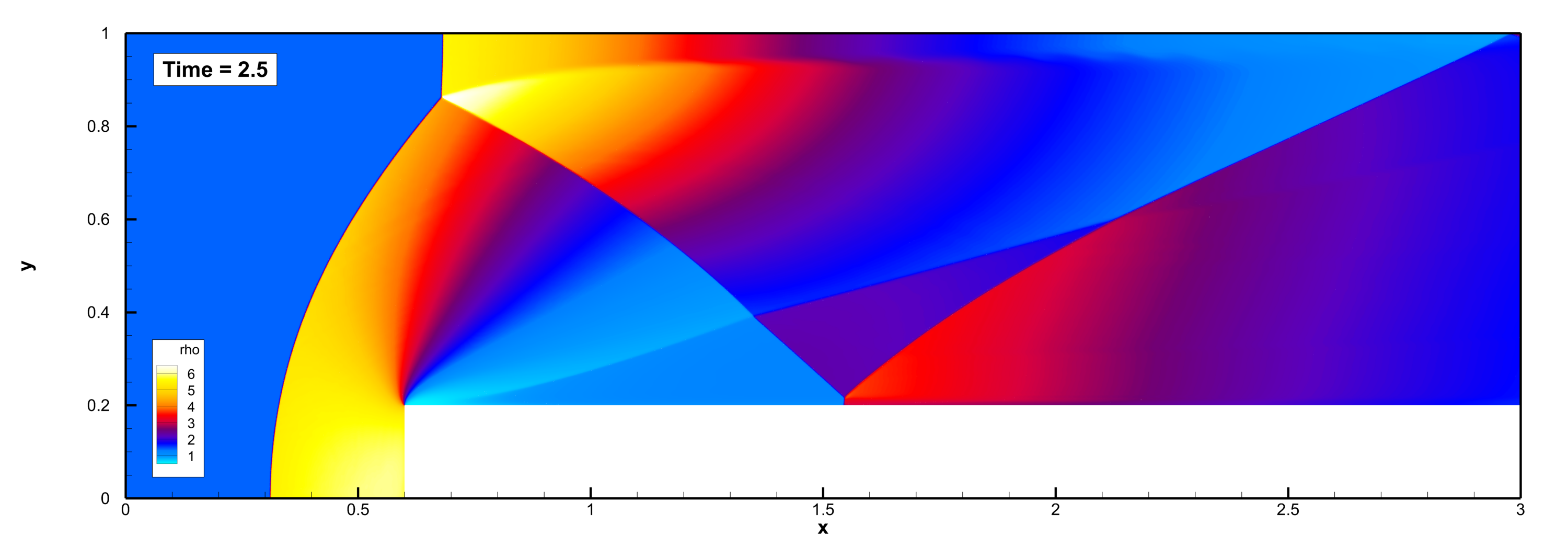}\\[-4pt]
	\includegraphics[width=0.88\linewidth]{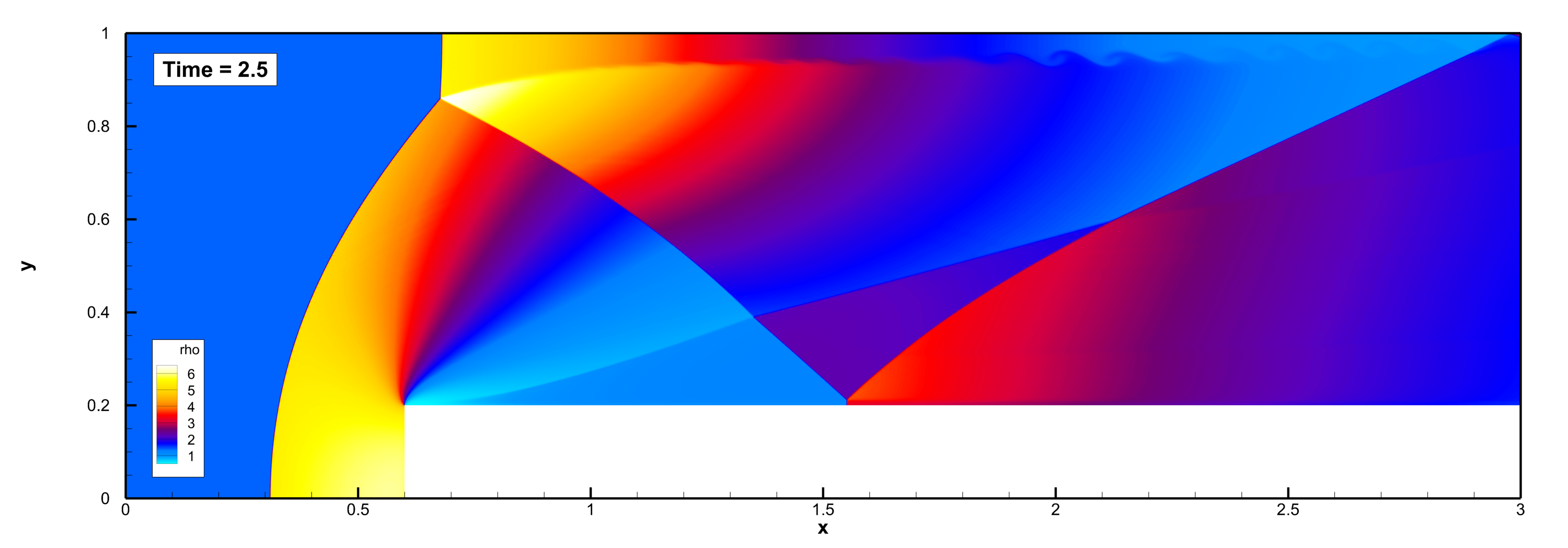}\\[-4pt]
	\includegraphics[width=0.88\linewidth]{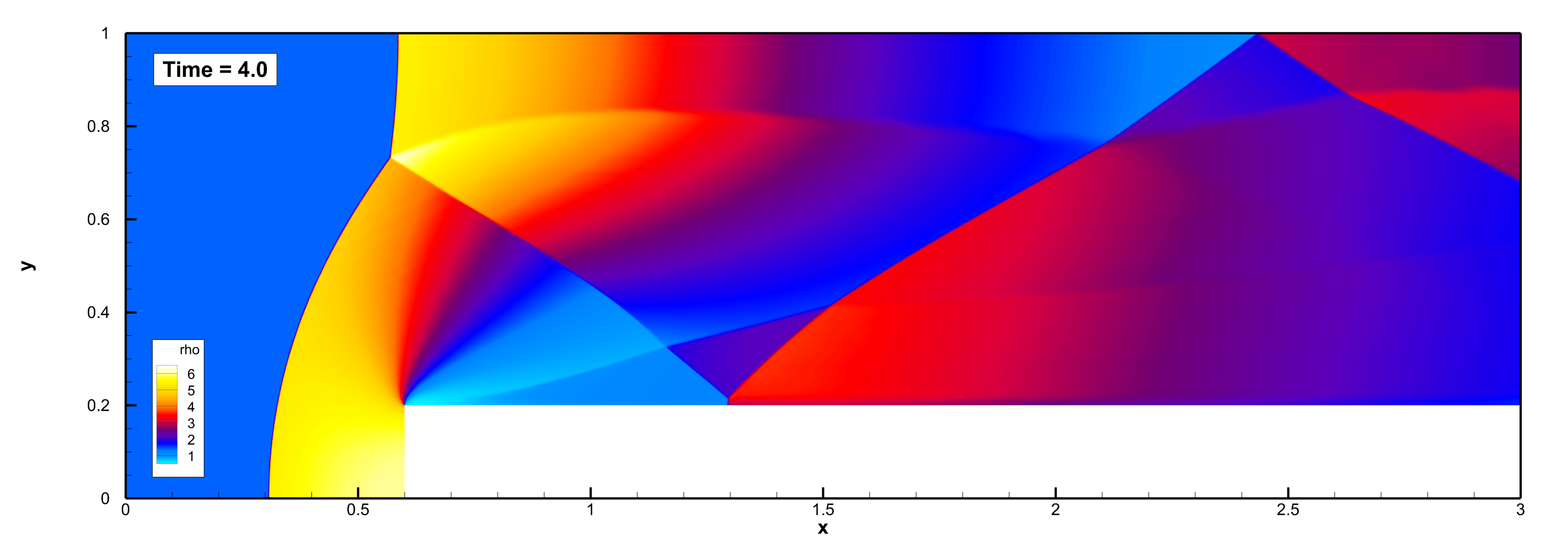}\\[-4pt]
	\includegraphics[width=0.88\linewidth]{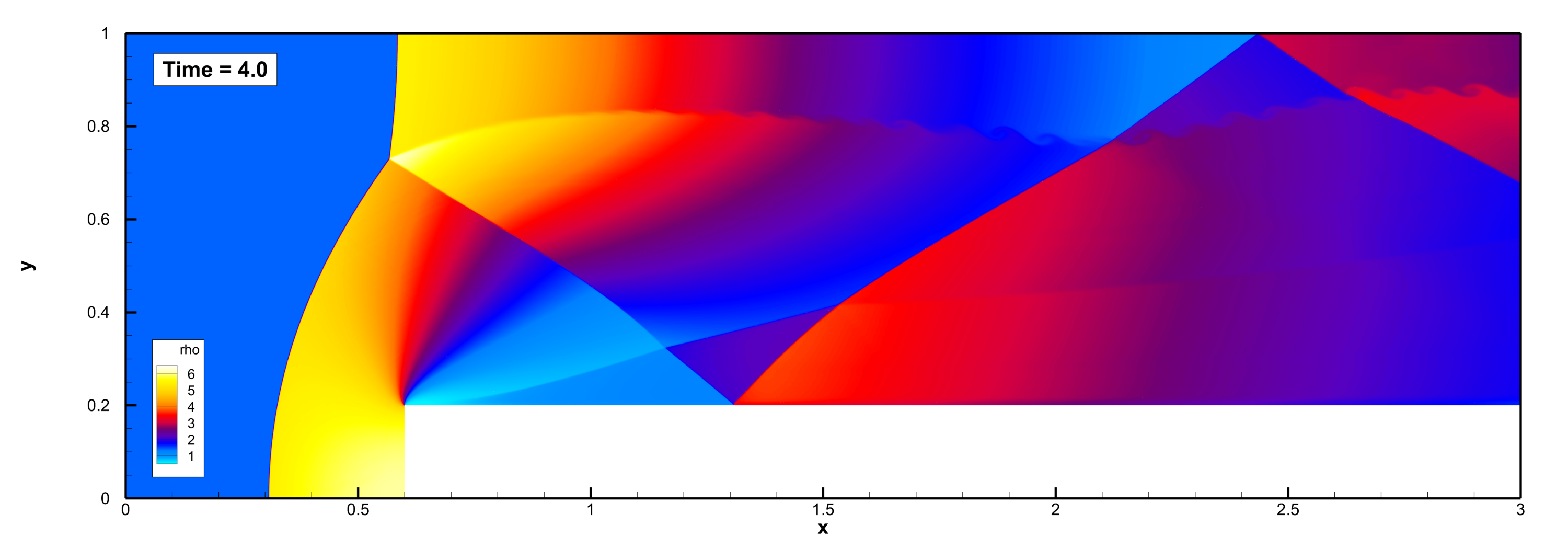}
	\caption{Numerical results for the forward facing step benchmark at time $t=2.5$ and $t=4$ obtained 
	with our quasi-concervative third order \Pdue (first and third image) and fourth order \Ptre (second 
	and fourth image) DG schemes on a polygonal tessellation with characteristic mesh size of $h=1/80$.} %1.3E{-2}
	\label{fig.FFS}
\end{figure}

We now consider the forward facing step problem which is a truly two-dimensional test case to be 
studied on a non trivial domain composed as follow
$\Omega = [0,7]\times[0,3] - [0.6,7]\times[0,0.2]$. 
This challenging benchmark has been introduced in~\cite{emery1968evaluation} and further studied by 
Woodward and Colella in~\cite{woodward1984numerical}.
The problem begins with uniform Mach 3 flow in a wind tunnel containing a step. 
The initial density of the tunnel is $rho_0 = \gamma = 1.4$ and the initial pressure is $p_0 = 1$.
For what concerns boundary conditions the exit right boundary has no effect on the flow, because the exit velocity is always supersonic.
On the other side of the domain the gas is continually fed in from the left-hand boundary by setting 
the left ghost state always equal to the initial conditions. 
Finally, wall boundary conditions are set on the bottom and top part of the domain. 
We report the numerical results obtained with our third order \Pdue and fourth order \Ptre quasi-conservative 
DG scheme in Figure~\ref{fig.FFS} at time $t=2.5$ and $t= 4.0$.

\subsubsection{Double Mach Reflection}

\begin{figure}[!b]
	\centering
	\includegraphics[width=0.711\linewidth,height=0.22\textheight]{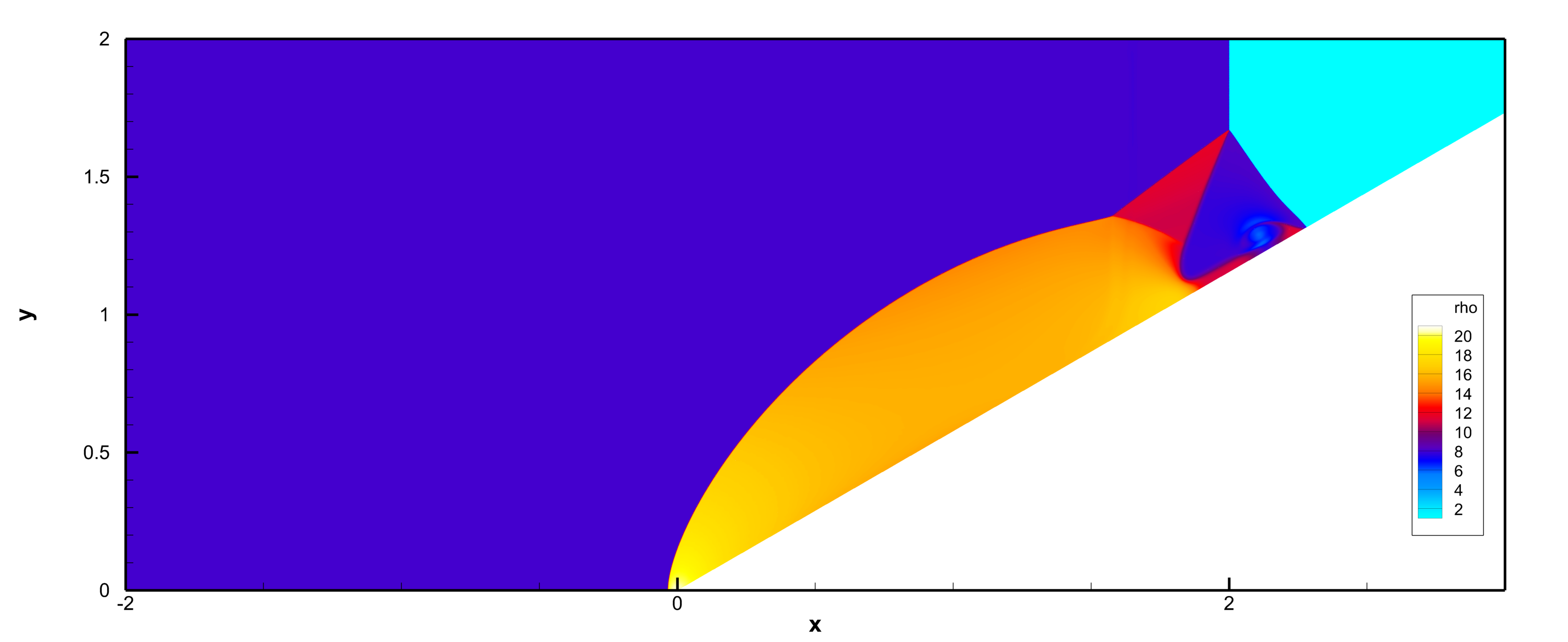}%
	\includegraphics[width=0.289\linewidth,height=0.22\textheight]{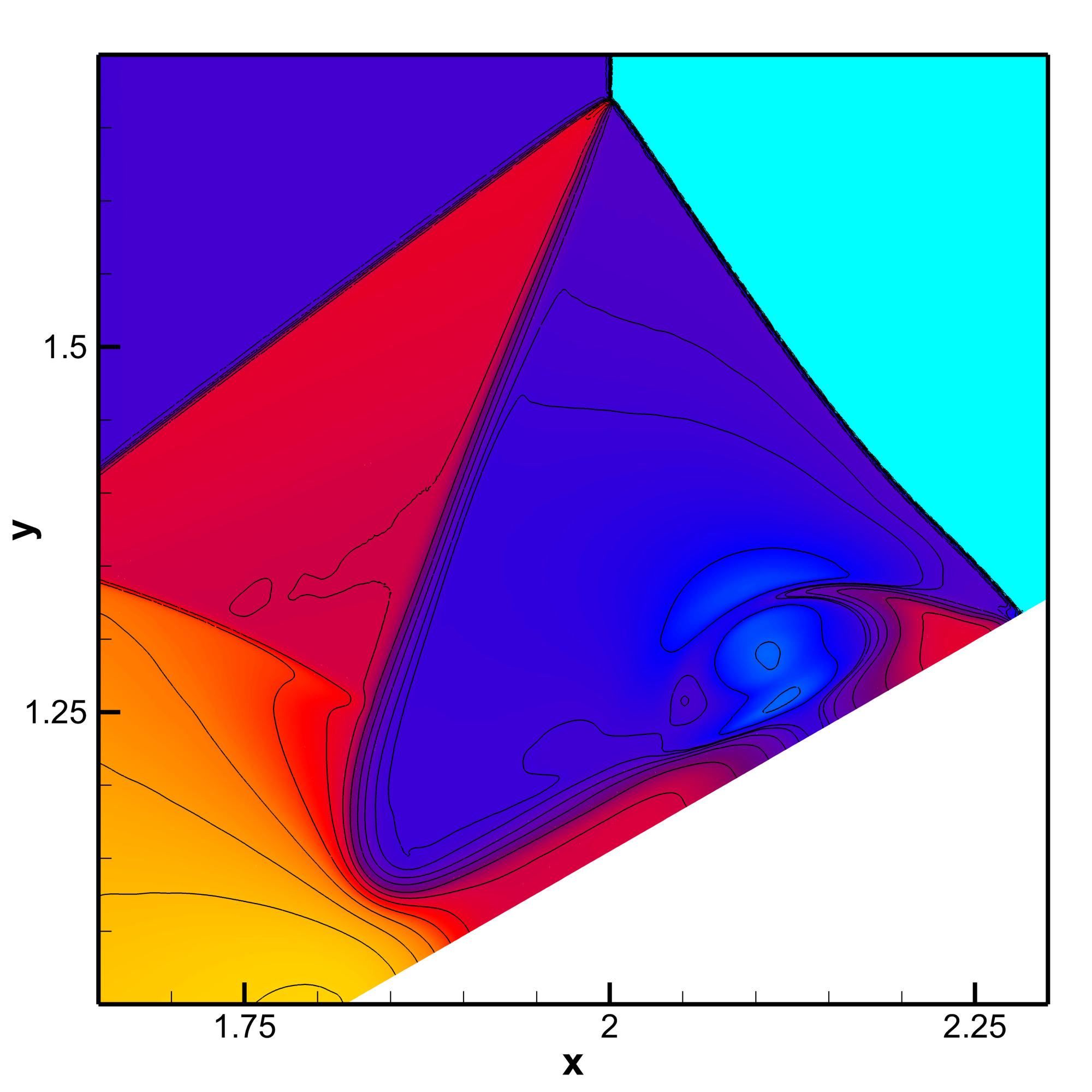}\\
	\includegraphics[width=0.711\linewidth,height=0.22\textheight]{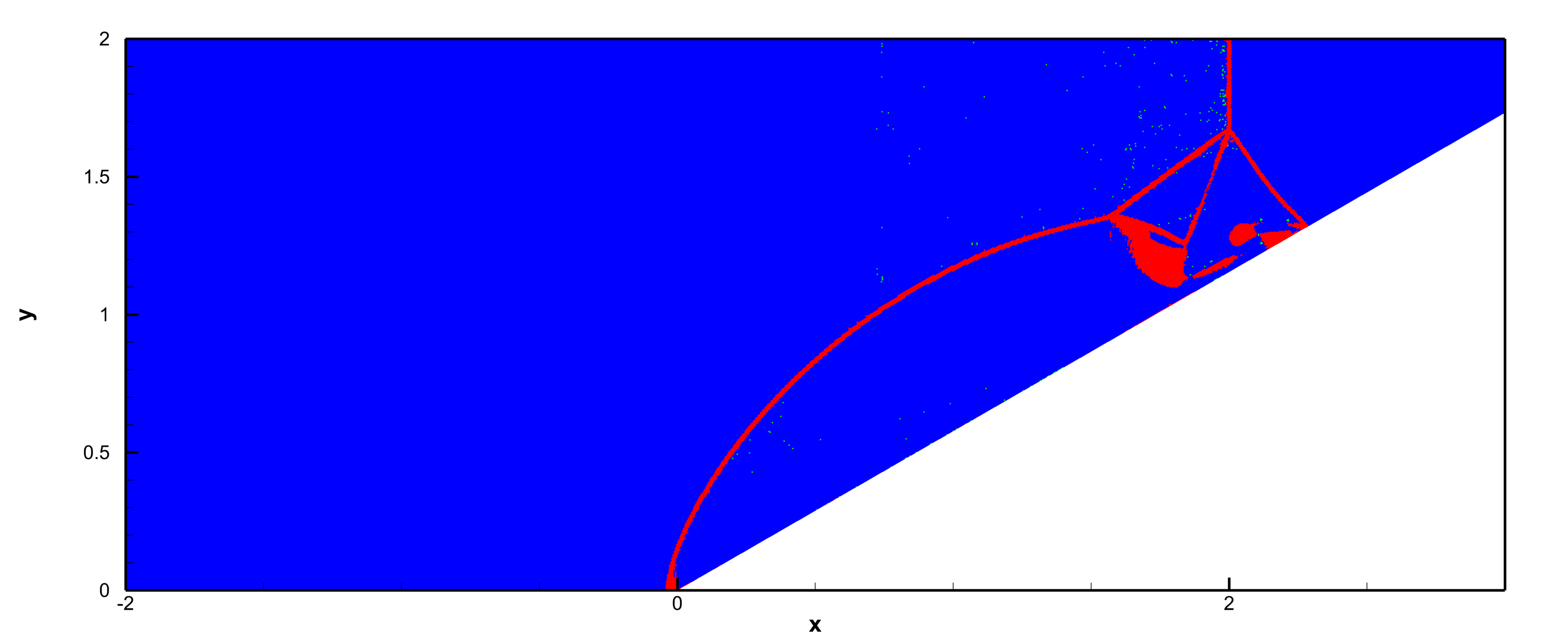}%
	\includegraphics[width=0.289\linewidth,height=0.22\textheight]{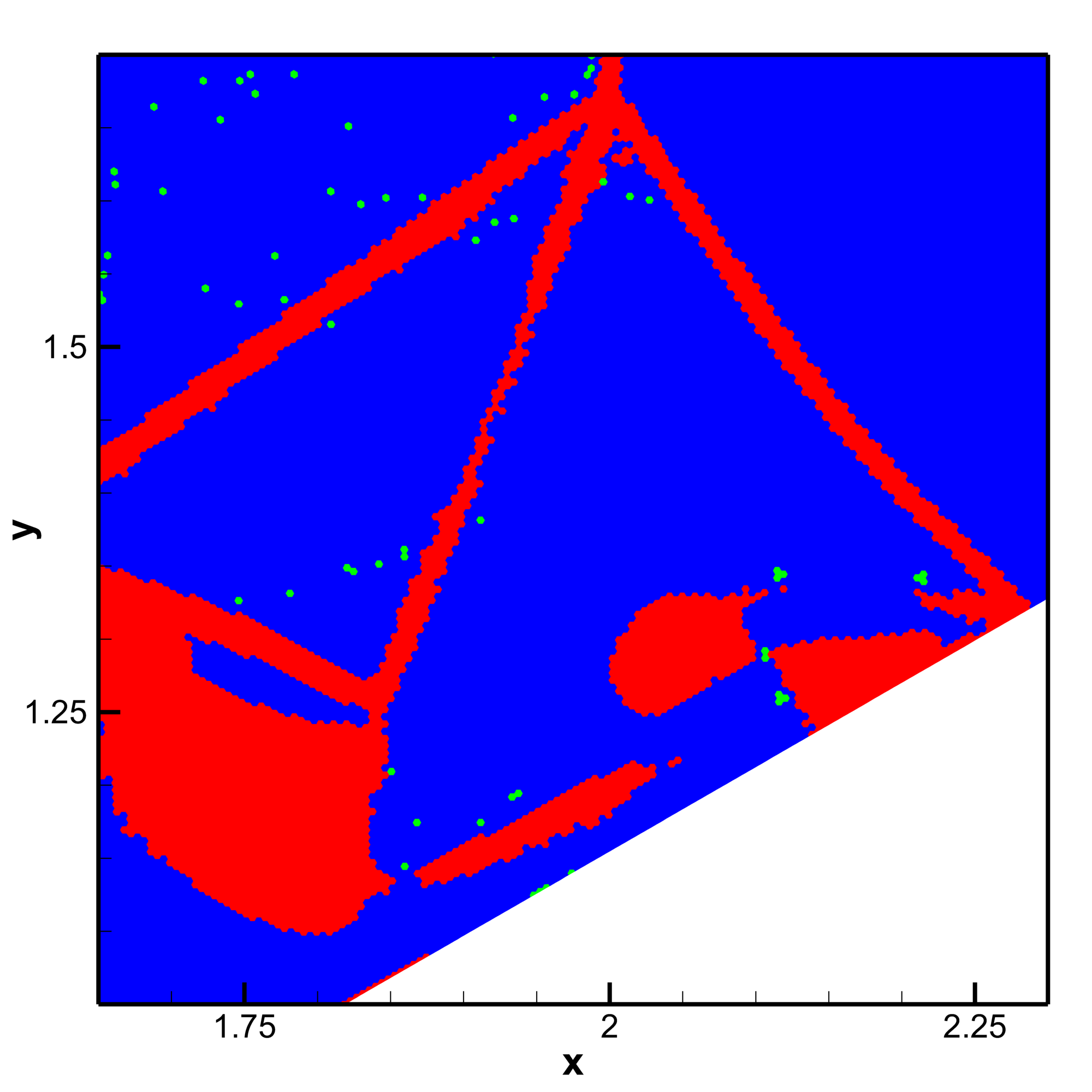}\\
	\caption{Numerical results for the Double Mach Reflection at time $t_f=0.2$ obtained with our 
	quasi-concervative third order \Pdue discontinuous Galerkin scheme on a polygonal tessellation 
	with characteristic mesh size of $h=1/220$: %4.6E{-3}
		density (top) and limiter type (bottom).}
	\label{fig.DMR}
\end{figure}

We close our set of validating benchmarks with the double Mach reflection problem 
which also has been originally proposed by Woodward and Colella in~\cite{woodward1984numerical},
where it has been studied on a rectangular domain thanks to a rotation of the initial conditions. 
Here, since we employ unstructured meshes, the problem can be directly run in physical coordinates as done in~\cite{dumbser2007arbitrary}
by considering as computational domain the polygon obtained by connecting the points $(-2,0)-(0,0)-(3,\sqrt{3})-(3,2)-(-2,2)$.
We consider a very strong shock wave that is moving along the x-direction at Mach $M=10$
and impacting against the ramp of angle $\pi/6$.
The initial condition is given by 
\begin{equation}
	(\rho, u, v, p)(\x) = 
	\begin{cases} 
		(8, 8.25, 0, 116.5  )  &\ \text{ if } \  x \le 0, \\	
		(\gamma, 0, 0, 1       ) & \ \text{ if } \  x > 0,		 \\		
	\end{cases}
\end{equation}
with $\gamma =1.4$ and the final computational time is $t_f= 0.2$.
The strong shock wave that hits the ramp causes the development of other two shock waves, 
one traveling towards the right and the other one propagating towards the top boundary of the domain.
The results obtained with our quasi-conservative third order \Pdue DG scheme on a polygonal tessellation
with characteristic mesh size  $h=1/220$, %$4.6E-3 
are depicted in Figure~\ref{fig.DMR}. 
The small-scale structures produced by the roll-up of the shear layers behind the shock wave are clearly 
visible from the density contour lines 
for which also a zoom is provided on the right part of the Figure. 
Moreover, we depict in red the \textit{shock-triggered troubled cells} which are treated by the a posteriori subcell limiter 
with the conservative correction, and in green the other limited cells treated instead in primitive variables. 
The use of the conservative formulation just on the shock regions allows to correctly reproduce the dynamics of this benchmark, 
while the majority of the domain can still be simulated in primitive variables.

\subsection{Euler equations with {space-time-dependent} EOS parameters}

\subsubsection{Quirk and Karni’s test case}

\begin{figure}[!b]
	\centering
	\includegraphics[width=0.9\linewidth,trim=5 5 5 5,clip]{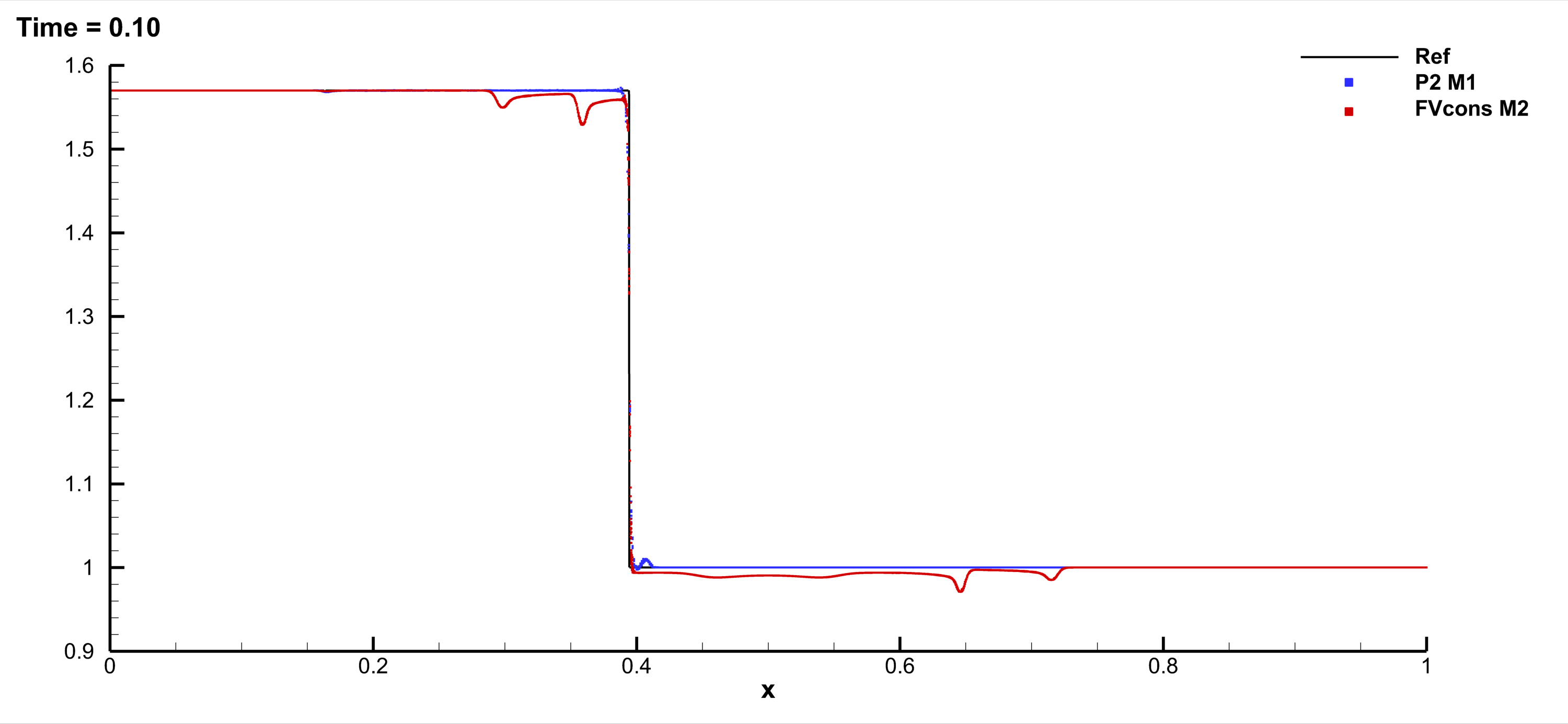}\\%
	\includegraphics[width=0.9\linewidth,trim=5 5 5 5,clip]{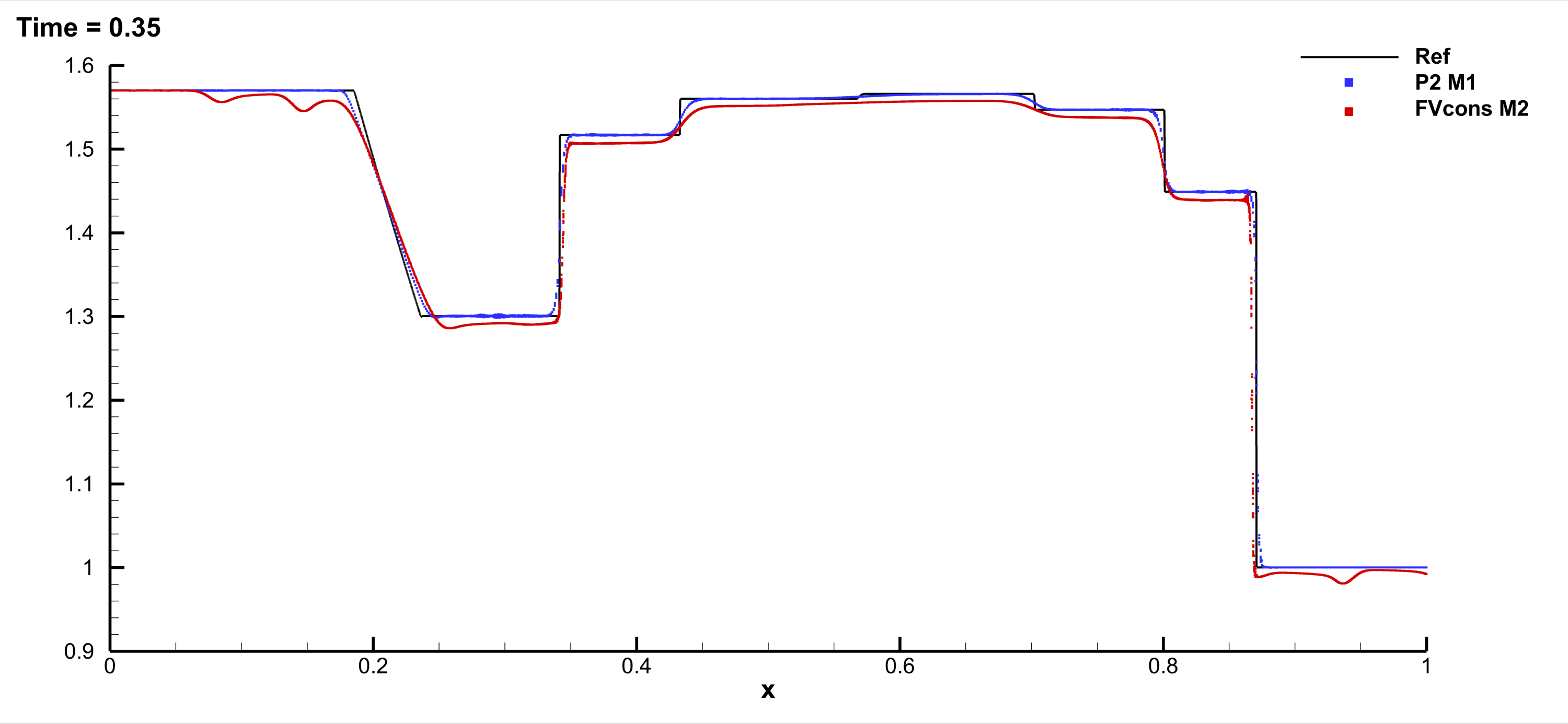}%
	\caption{Pressure profile for the essentially 1D Quirk and Karni’s test case at time $t=0.1$, i.e. just 
	before the shock hits the 1D-bubble, and at the final time $t=0.35$. 
		We compare the numerical results obtained with our third order \Pdue quasi-conservative DG scheme on a 
		polygonal tessellation with characteristic mesh size $h=1/270$ %=3.64E-3
		(blue) with those obtained with a standard conservative second order FV scheme on a fine triangulation 
		with $h=1/2690$ %3.72E-4
		(red); we also provide a reference solution (black solid line). The solution obtained with the classical 
		conservative scheme, even if run on a very fine mesh, develops spurios pressure oscillations already before 
		the shock wave has hit the bubble which are then propagated around the domain and corrupt the solution.
	}
	\label{fig.BubbleHelium1D}
\end{figure}

We start with an essentially one-dimensional test case first proposed by Quirk and Karni in~\cite{quirk1996dynamics} 
and further investigated by Abgrall in~\cite{abgrall1996prevent}.
This 1D test case mimics the shock helium-bubble interaction, that we also study in the two-dimensional setting in 
the next Section, and shares with it the difficulties of multi-material simulations.
Indeed a classical simulation performed by using the conservative formulation~\eqref{eq.multimaterialCons} produces 
spurious oscillations that corrupt the final solution as clearly shown by the red graphs in Figure~\ref{fig.BubbleHelium1D}.

The initial setting for this problem is the following
\begin{equation}
	(\rho, u, v, p,c)(\x) = 
	\begin{cases} 
		( 1.3765,  0.3948, 0, 1.57, 1)  & \ \text{ if } \  x \le 0.2, \\[2pt]	
		( 0.138,  0,  0, 1,  \frac{5}{3\gamma} ) &  \ \text{ if } \   0.4 \le x \le 0.6,		 \\[2pt]		
		( 1,   0, 0, 1,  1,  1)  &  \ \text{{ otherwise}}, 		
	\end{cases}
\end{equation}
on a rectangular domain $\Omega=[0,1]\times[0, 0.05]$ with Dirichlet boundary conditions at the left 
and right borders and wall boundary conditions on the bottom and top edges.
We have reported the numerical results obtained with our \Pdue quasi-conservative approach in Figure~\ref{fig.BubbleHelium1D} which perfectly 
agree with the reference solutions we have obtained by solving Kapila's model~\cite{kapila2001} for one-velocity 
one-pressure compressible mutiphase flow (also known as five-equation Baer--Nunziato system~\cite{baernunziato}) 
on a very fine mesh of 20000 uniform cells by means of a second order path-conservative ADER-FV method 
with TVD reconstruction, 
specifically the one presented and validated in~\cite{hstglm, sisurfacetension}.

\subsubsection{Shock helium-bubble interaction}
\label{ssec.TwoPhase_2dbubble}

\begin{figure}[!bp]
	\centering
	\includegraphics[width=0.5\linewidth,trim=3 10 3 10,clip]{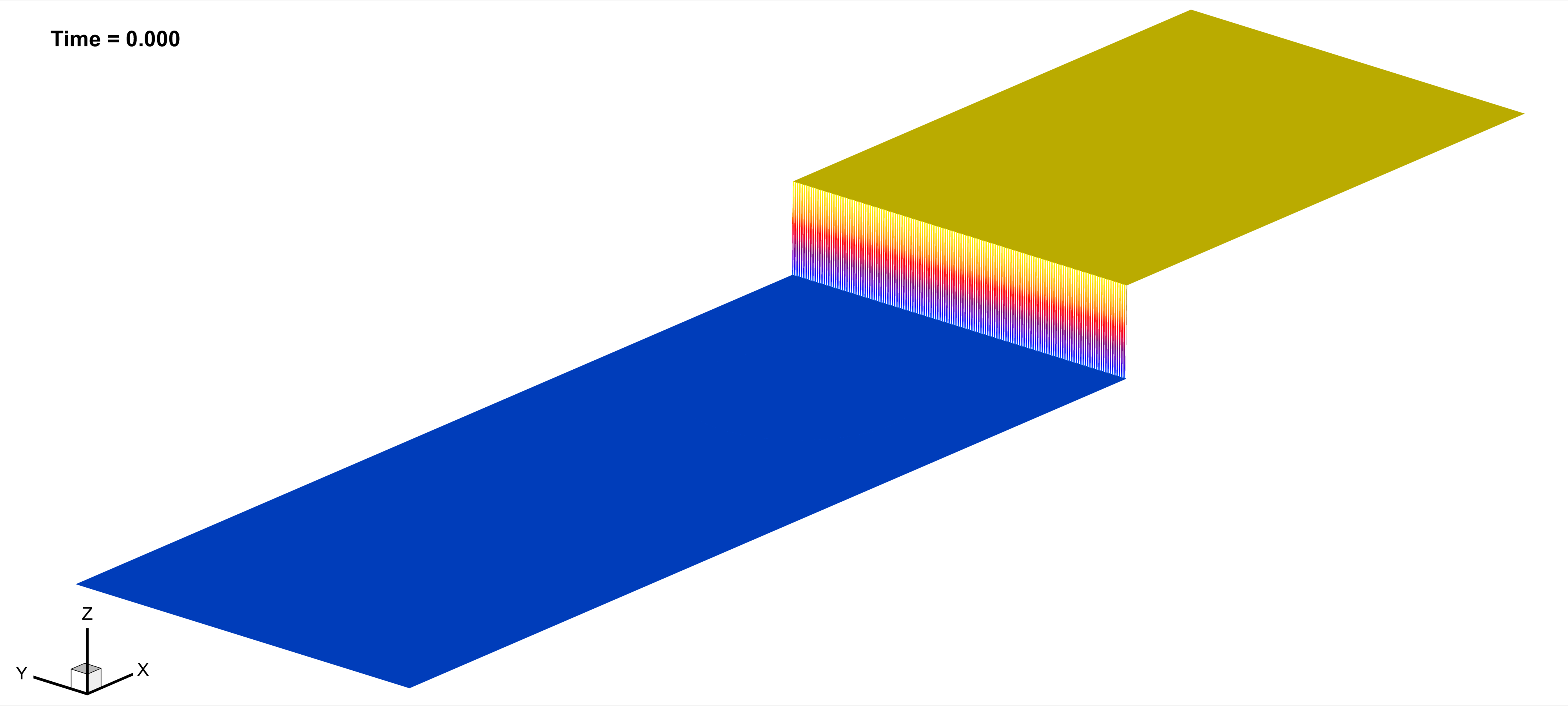}%
	\includegraphics[width=0.5\linewidth,trim=3 10 3 10,clip]{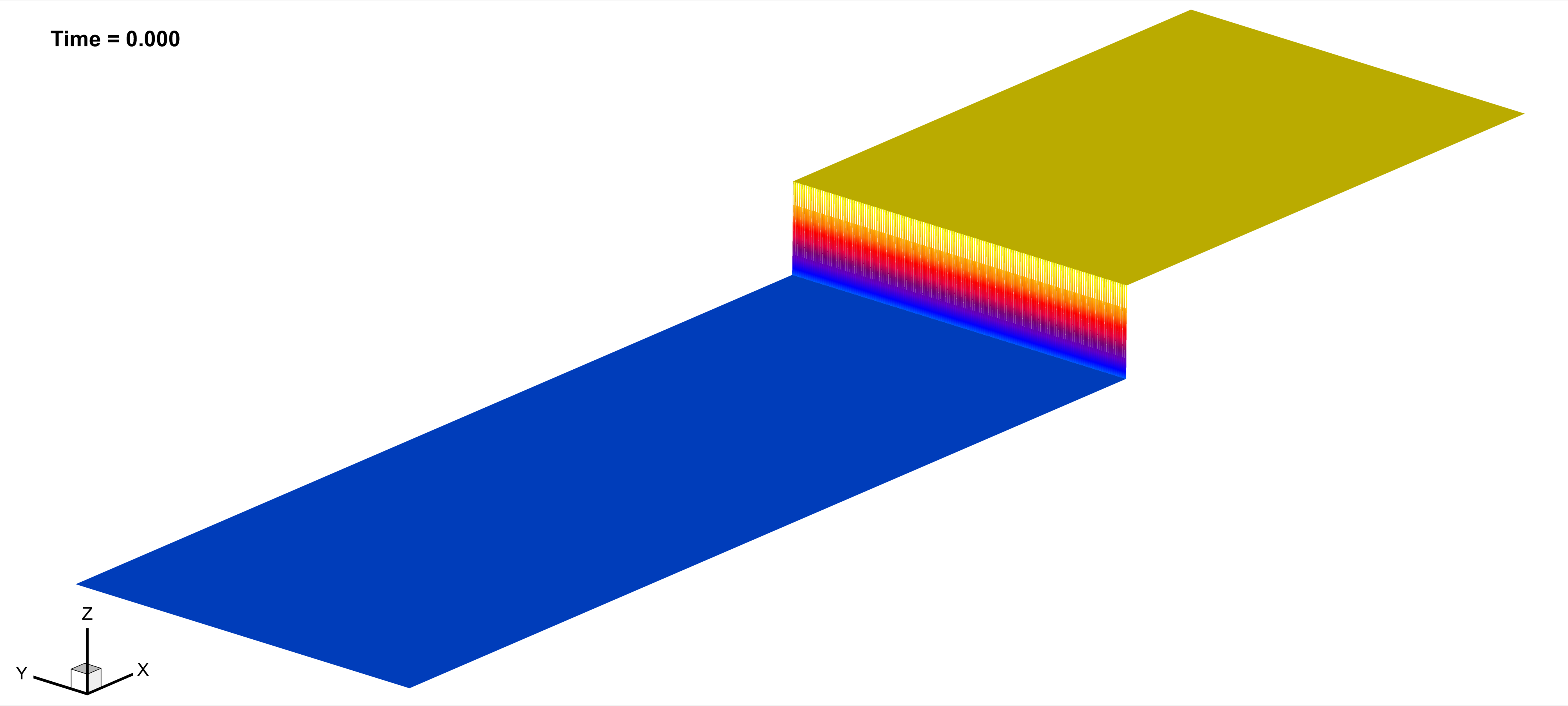}\\[22pt]
	\includegraphics[width=0.5\linewidth,trim=3 10 3 10,clip]{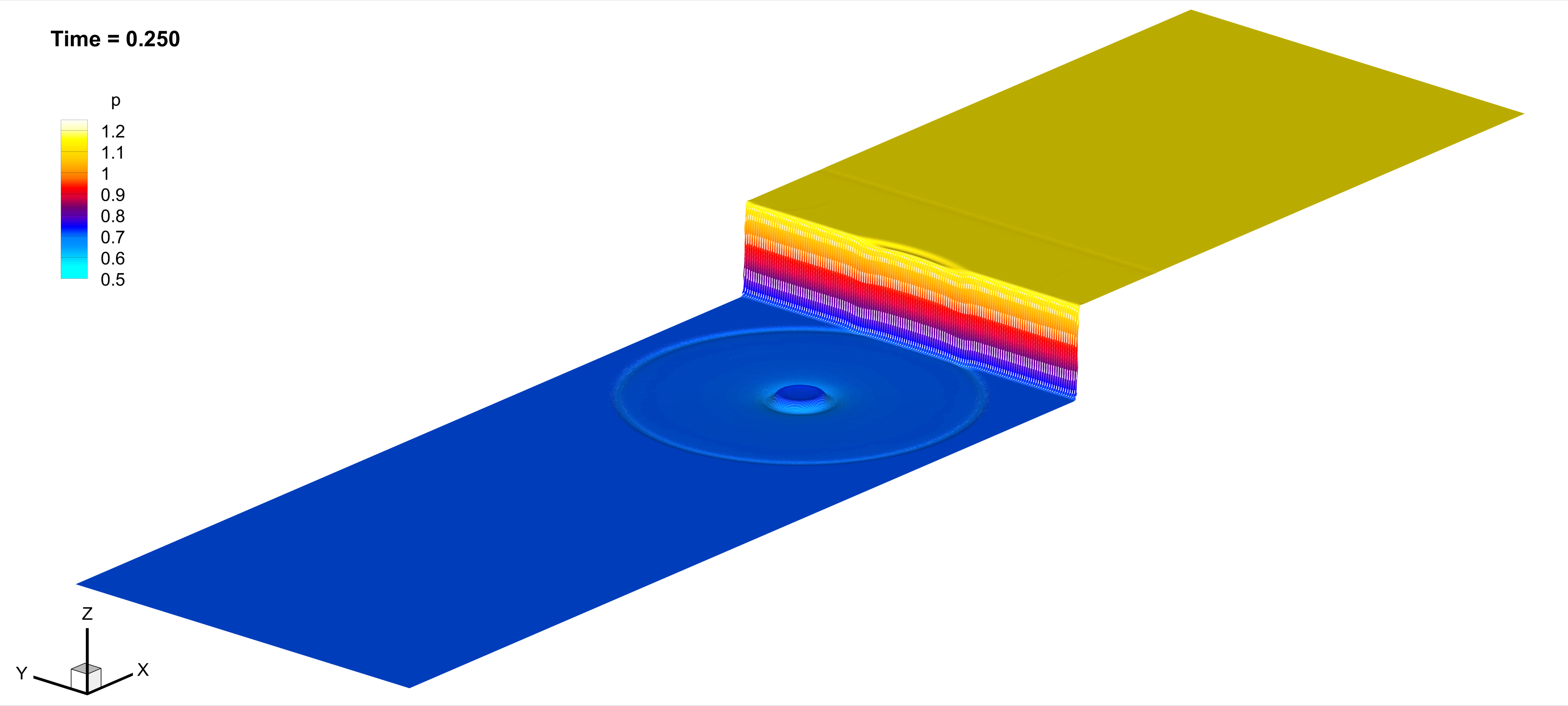}%
	\includegraphics[width=0.5\linewidth,trim=3 10 3 10,clip]{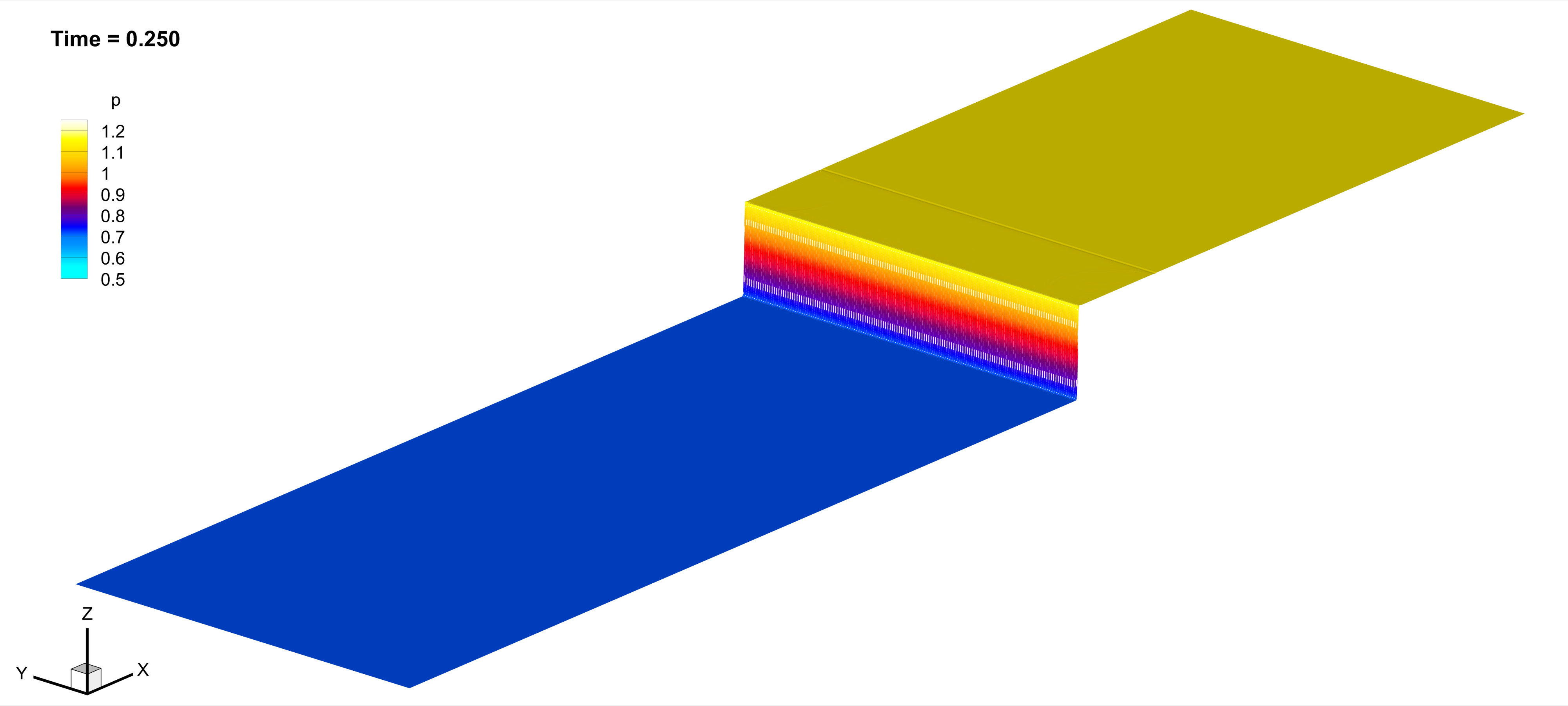}\\[22pt]
	\includegraphics[width=0.5\linewidth,trim=3 10 3 10,clip]{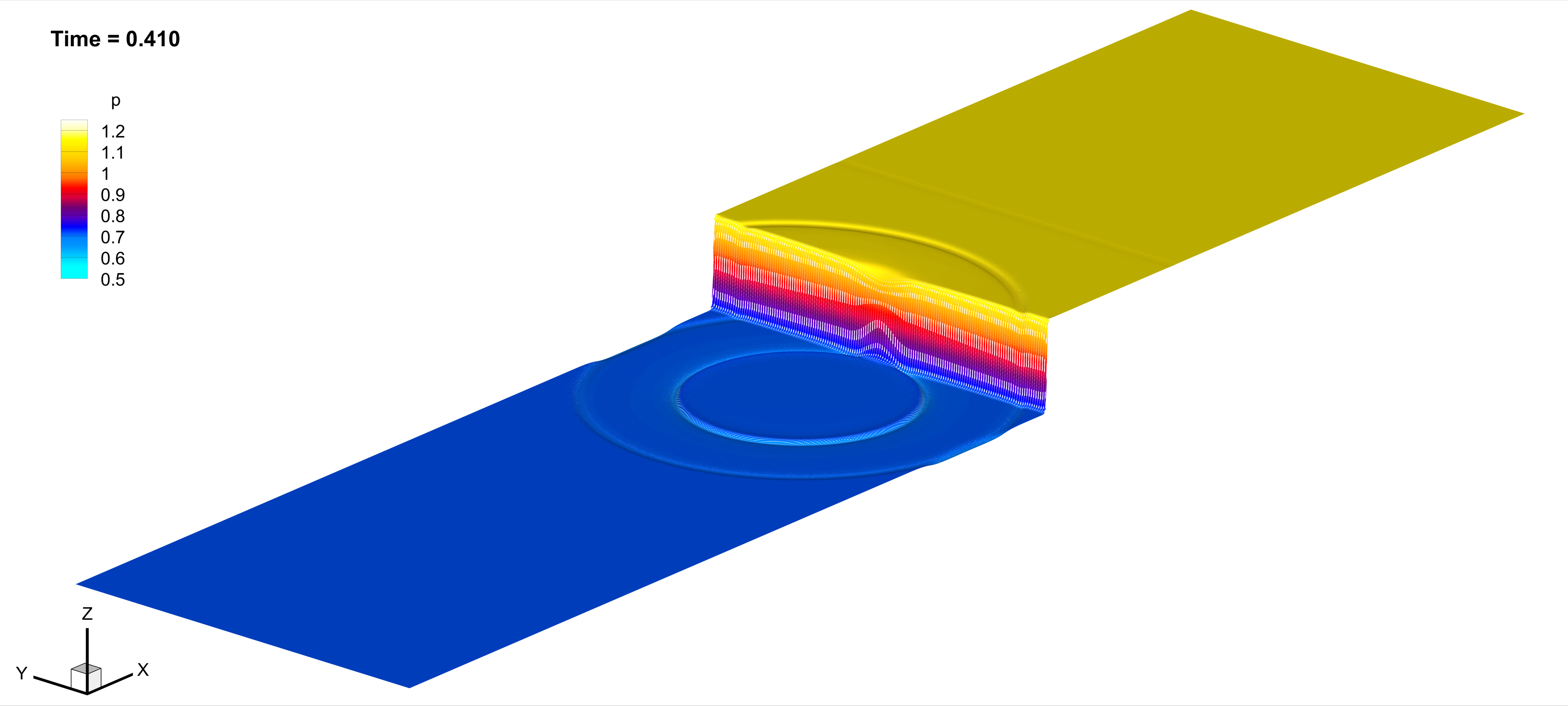}%
	\includegraphics[width=0.5\linewidth,trim=3 10 3 10,clip]{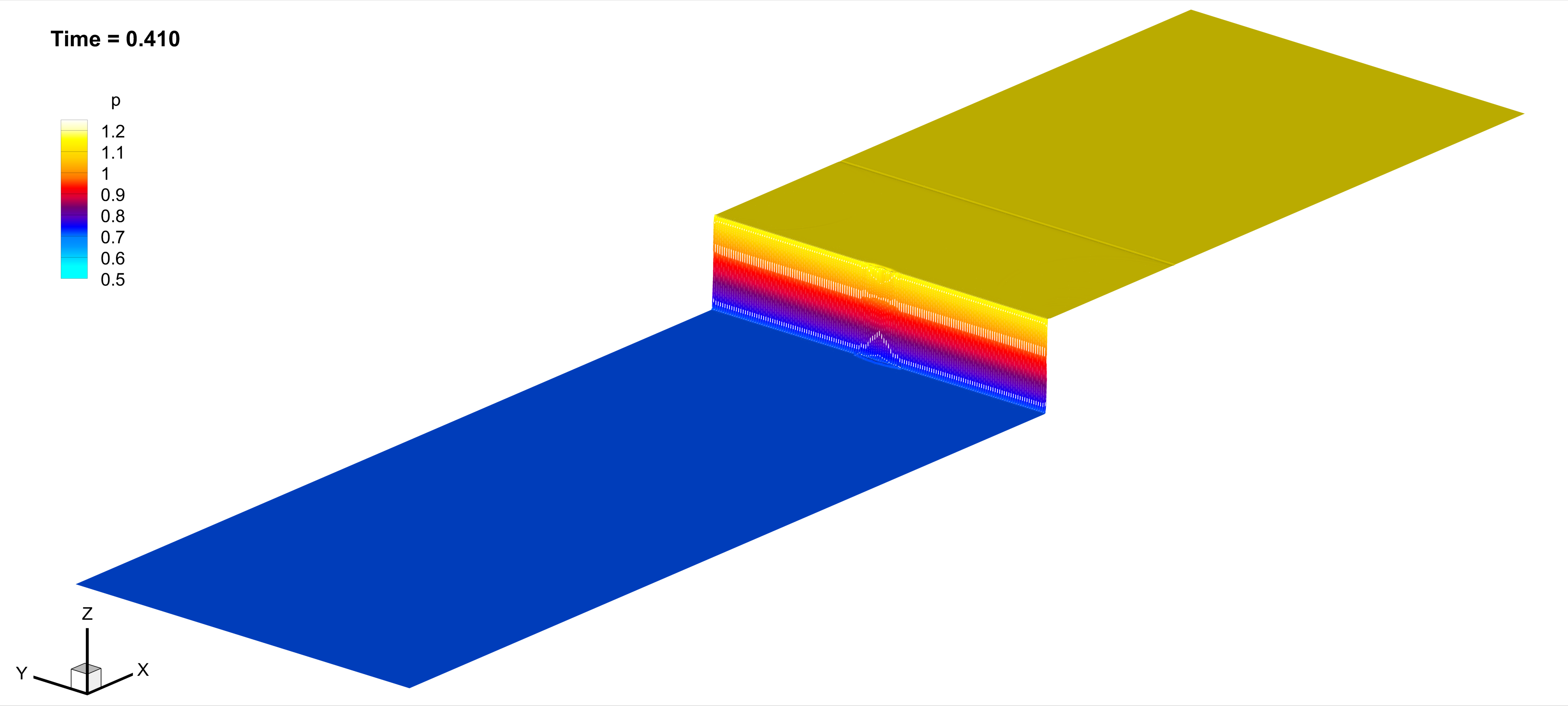}\\[22pt]
	\includegraphics[width=0.5\linewidth,trim=3 10 3 10,clip]{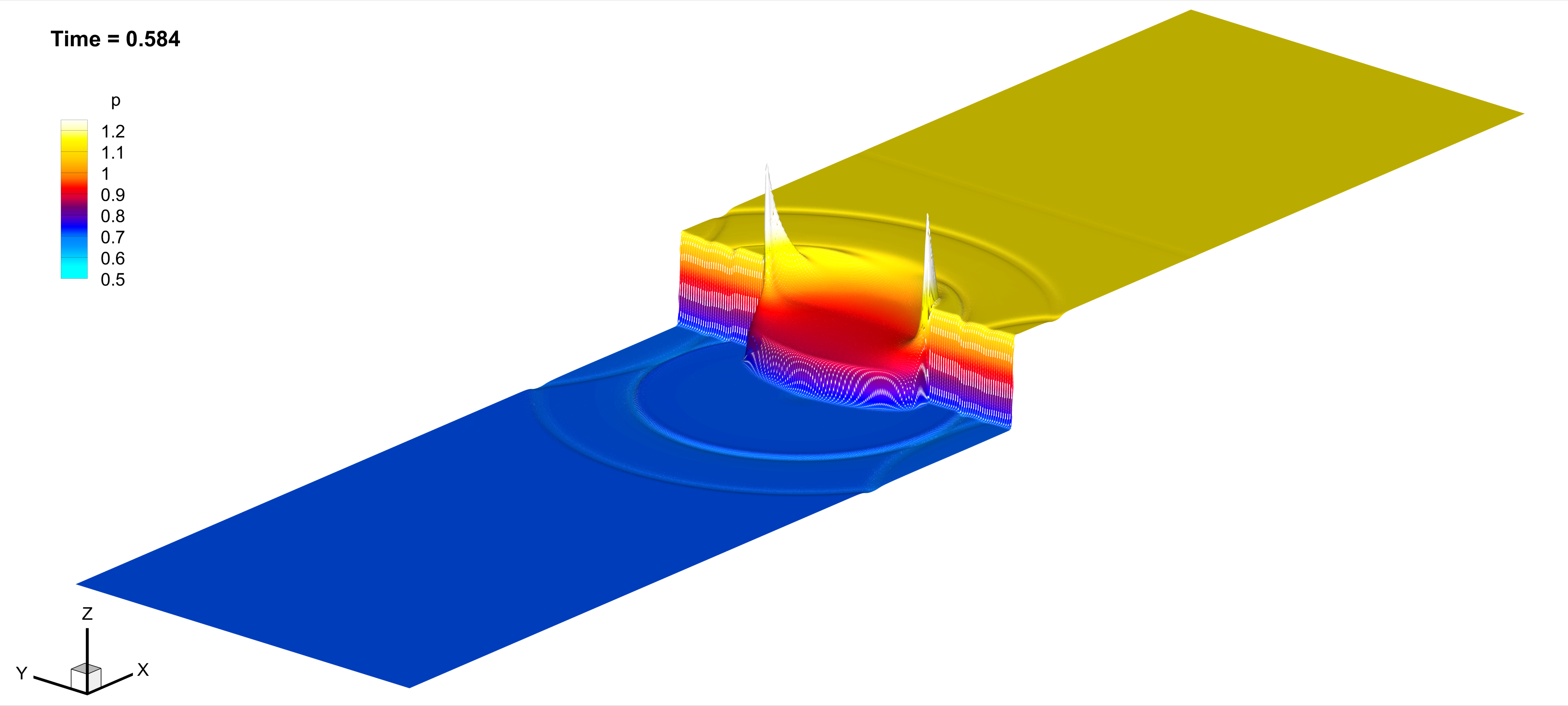}%
	\includegraphics[width=0.5\linewidth,trim=3 10 3 10,clip]{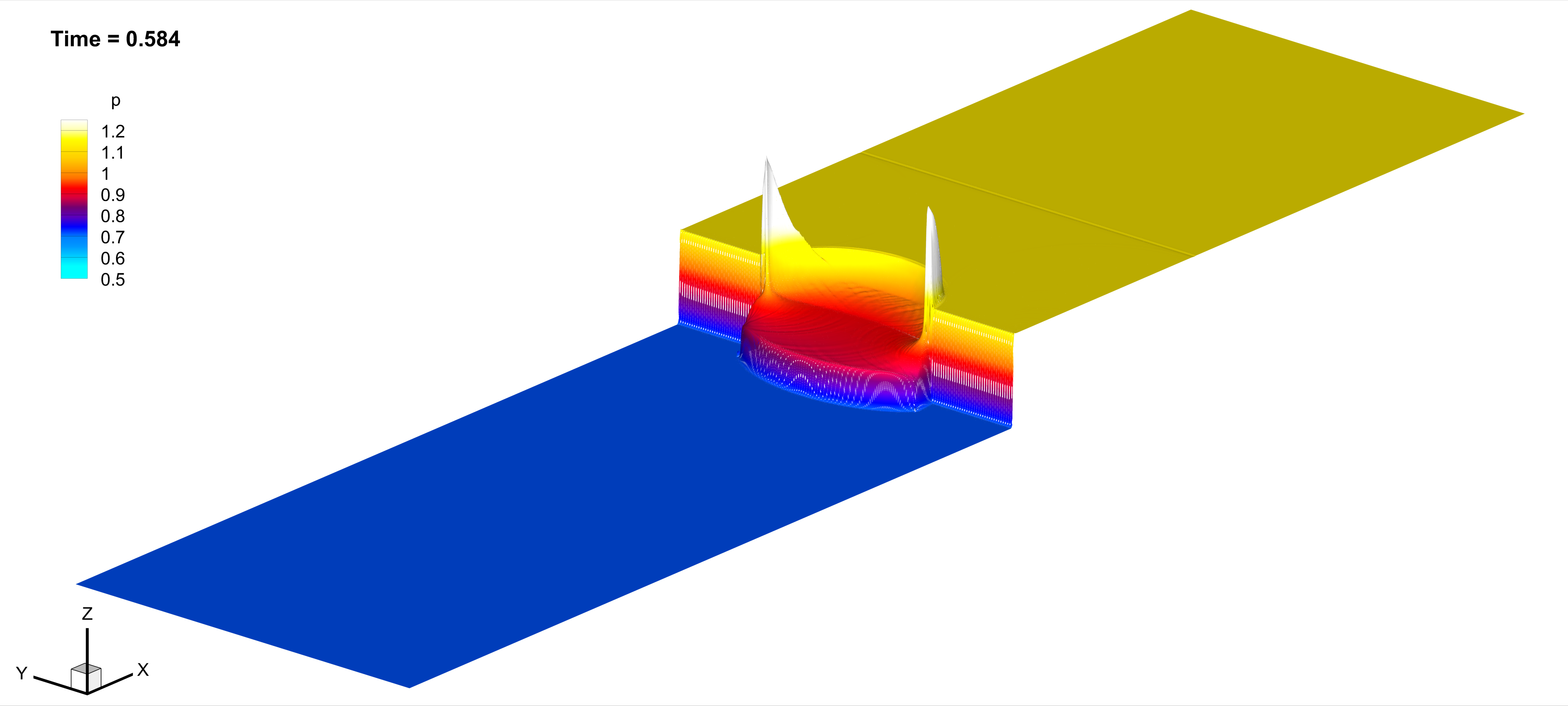}\\[22pt]
	\caption{Pressure profile for the shock helium-bubble interaction. In this first figure, we compare the 
	numerical results obtained with our third order \Pdue quasi-conservative DG scheme on a polygonal tessellation 
		with characteristic mesh size $h=1/150$  %$h=6.49E-3$ 
		(right) 
		with those obtained with a robust second order conservative ENO finite volume scheme 
		on a very fine triangular mesh with $h=1/750$ %$h=1.33E-3$ 
		(left). 
		We remark that the second order FV scheme used here for comparison reasons, 
		produces \textit{spurious} oscillations which propagate away from the interface
		between the air and the bubble, in the form of circular pressure waves.
		Instead, our quasi-conservative approach, by working in primitive variables in particular at the bubble interface (see also next figure), 
		completely avoids any spurious oscillations. Since both schemes employ at the limiter level the same ENO 
		reconstruction, it is clear that the oscillations are induced by the conservative fluxes at material interfaces
		and not related to Gibbs-like oscillations.
	}
	\label{fig.BubbleHelium-Shock_primvscons}
\end{figure}
\begin{figure}[!bp]
	\centering
	\includegraphics[width=0.5\linewidth,trim=3 6 3 6,clip]{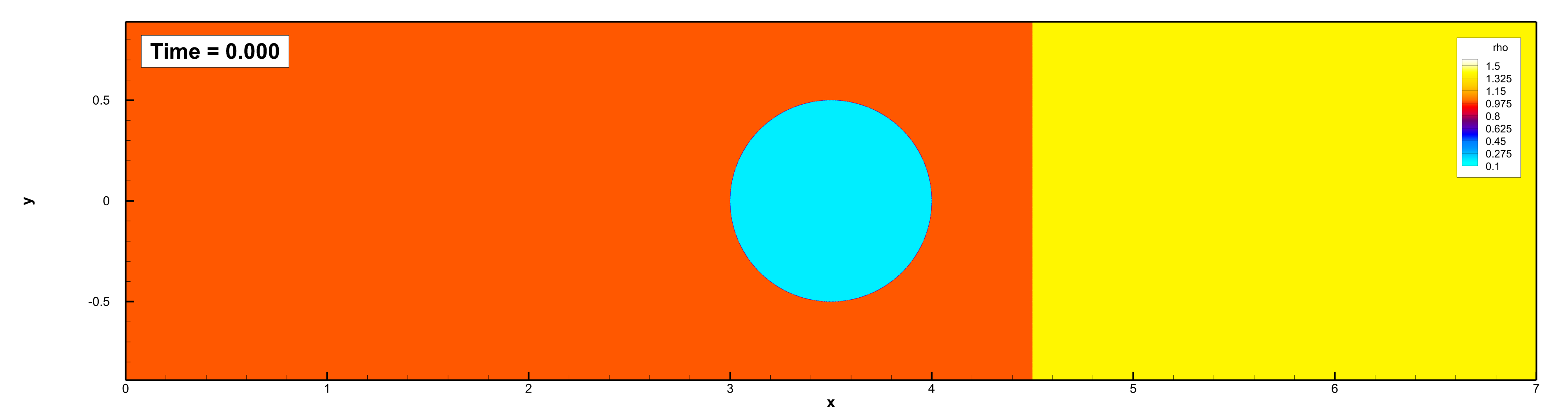}%
	\includegraphics[width=0.5\linewidth,trim=3 6 3 6,clip]{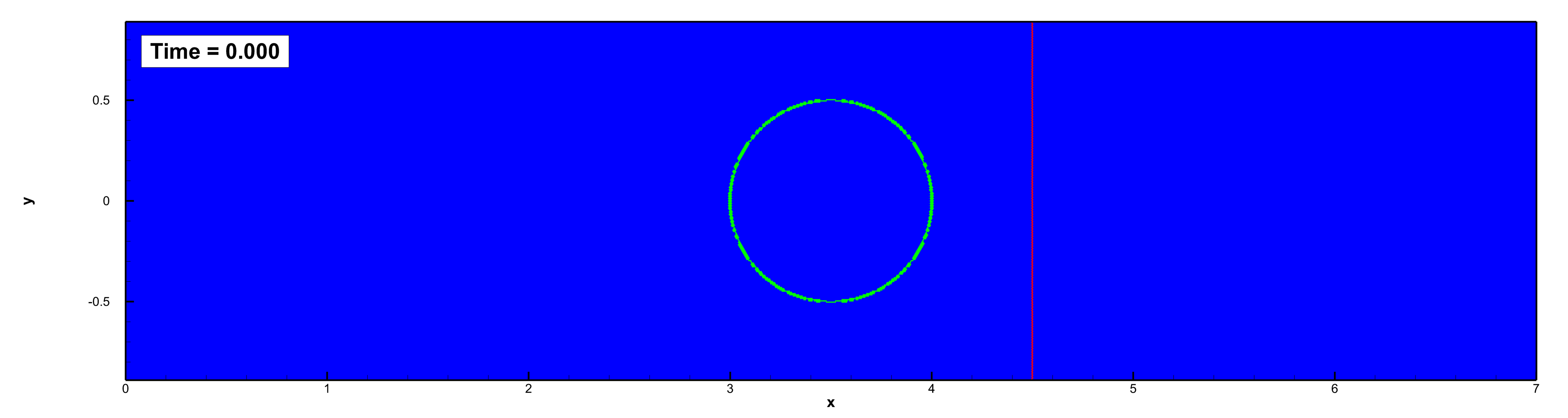}\\[-3pt]
	\includegraphics[width=0.5\linewidth,trim=3 6 3 6,clip]{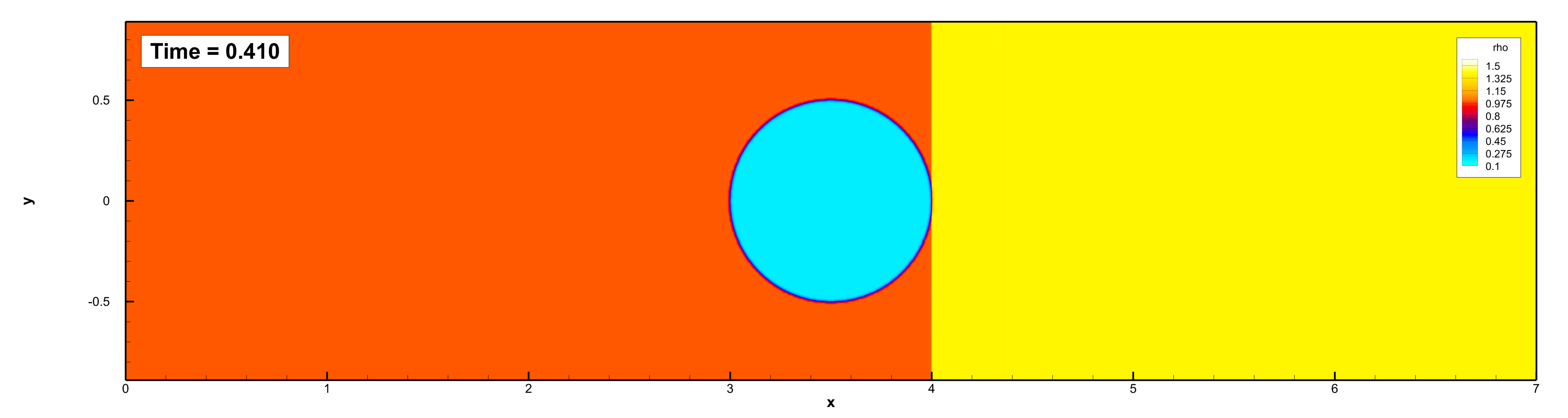}%
	\includegraphics[width=0.5\linewidth,trim=3 6 3 6,clip]{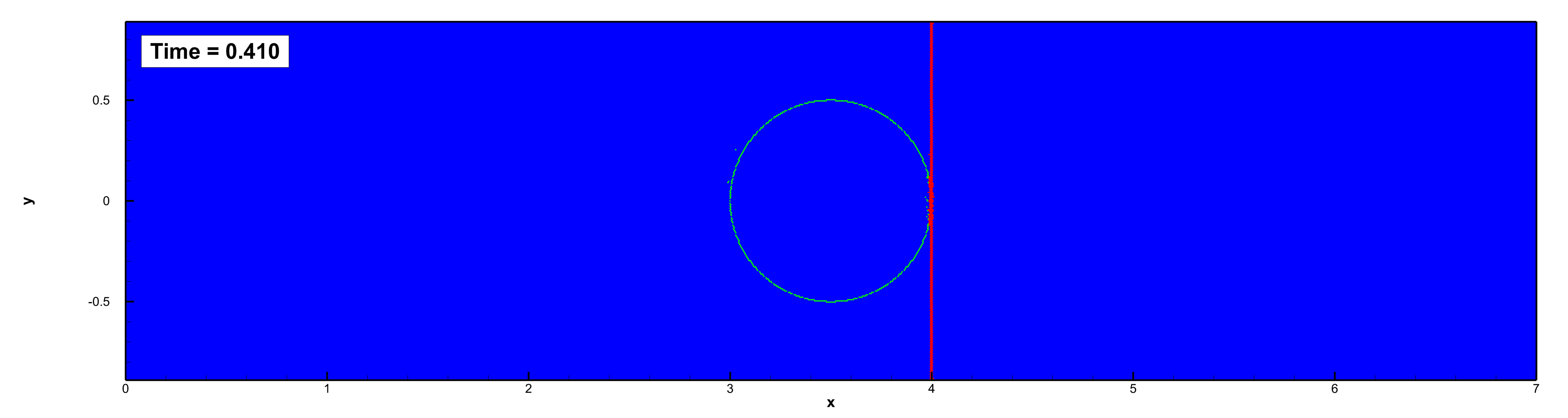}\\[-3pt]
	\includegraphics[width=0.5\linewidth,trim=3 6 3 6,clip]{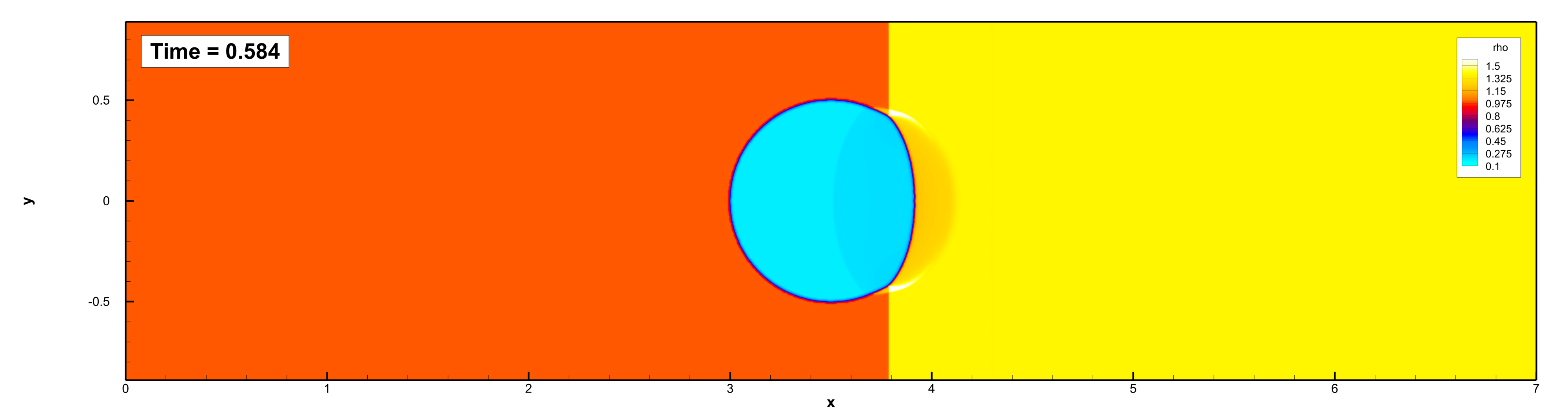}%
	\includegraphics[width=0.5\linewidth,trim=3 6 3 6,clip]{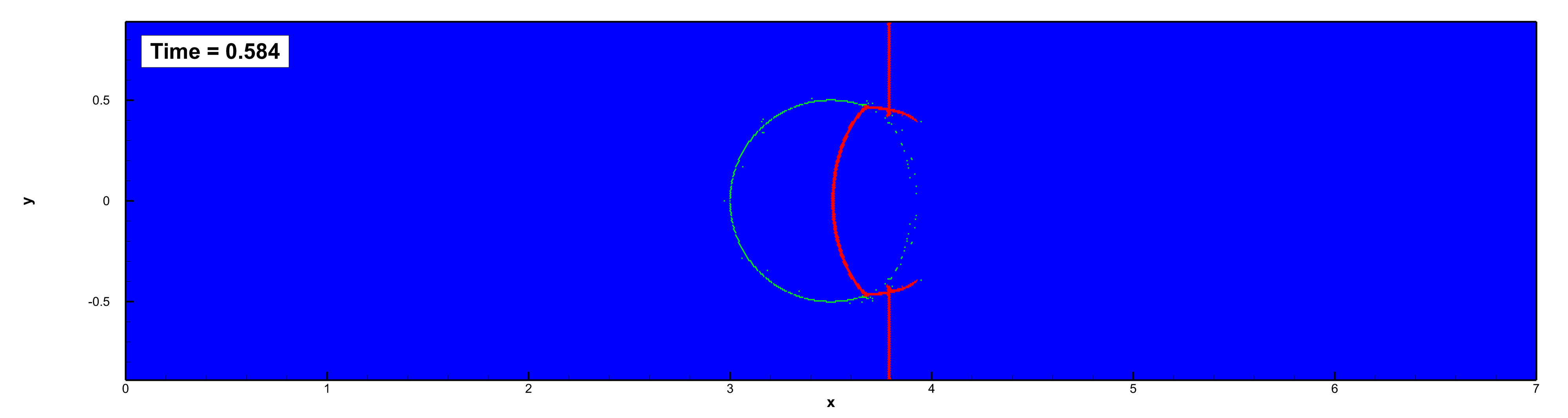}\\[-3pt]
	\includegraphics[width=0.5\linewidth,trim=3 6 3 6,clip]{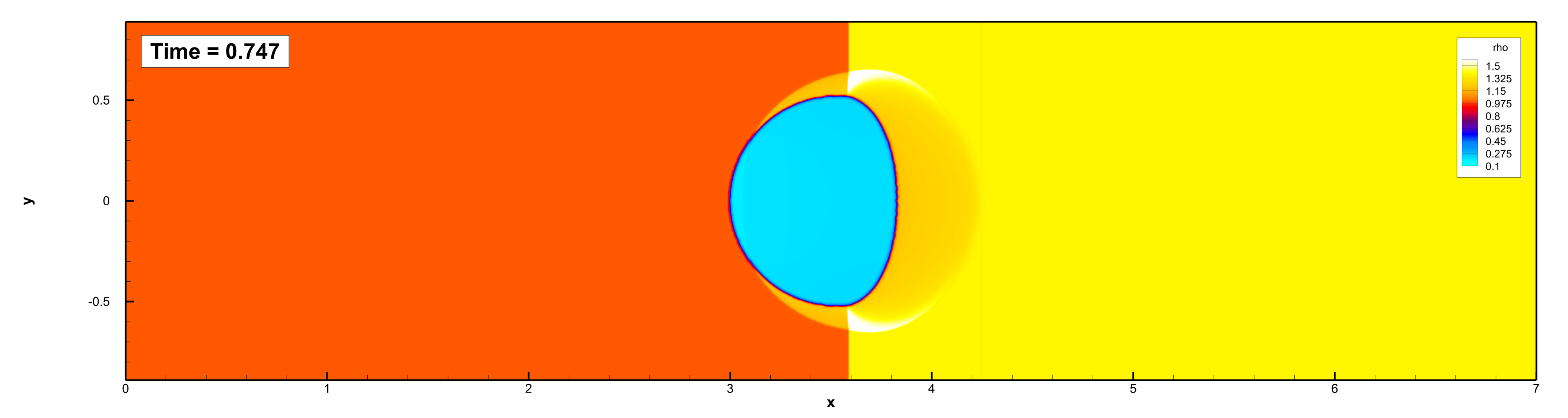}%
	\includegraphics[width=0.5\linewidth,trim=3 6 3 6,clip]{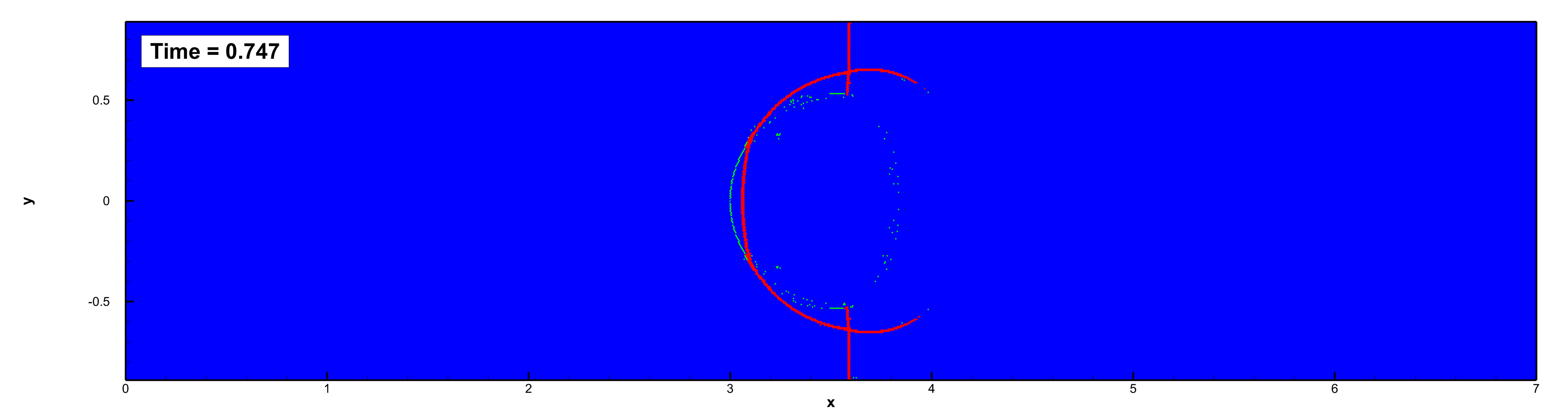}\\[-3pt]
	\includegraphics[width=0.5\linewidth,trim=3 6 3 6,clip]{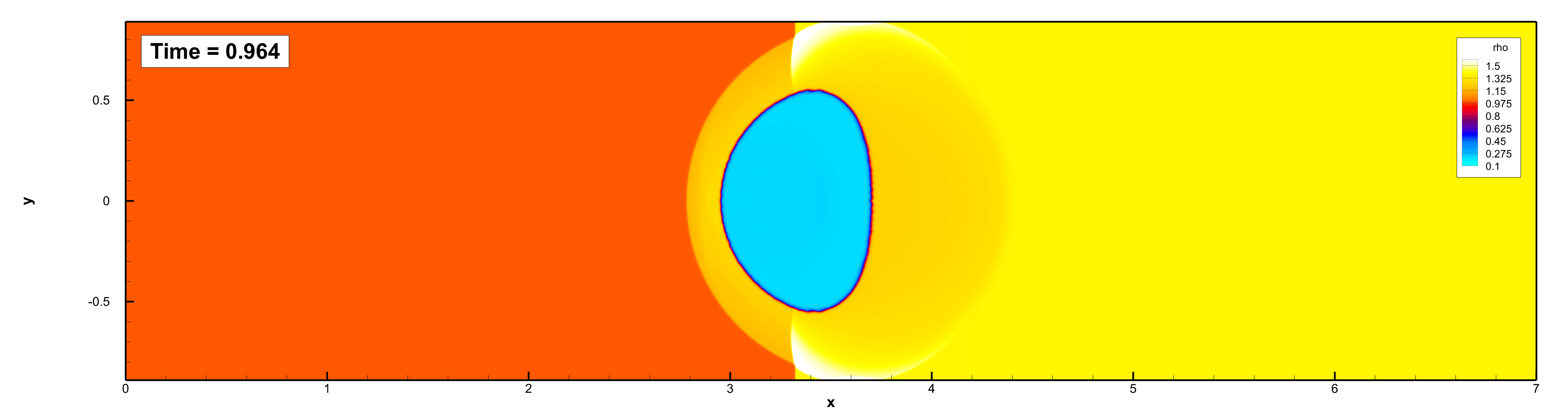}%
	\includegraphics[width=0.5\linewidth,trim=3 6 3 6,clip]{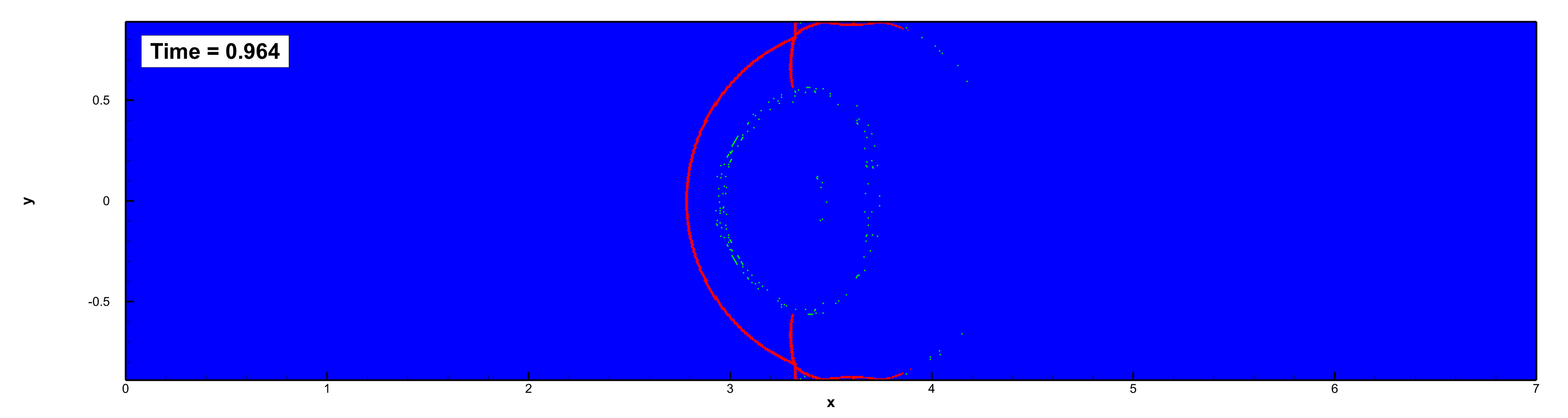}\\[-3pt]
	\includegraphics[width=0.5\linewidth,trim=3 6 3 6,clip]{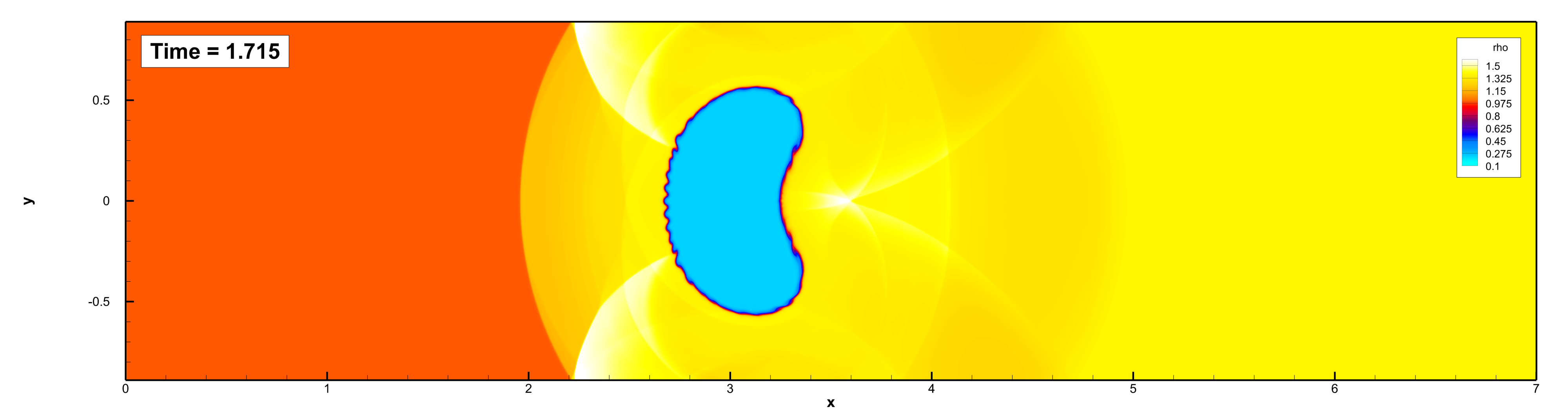}%
	\includegraphics[width=0.5\linewidth,trim=3 6 3 6,clip]{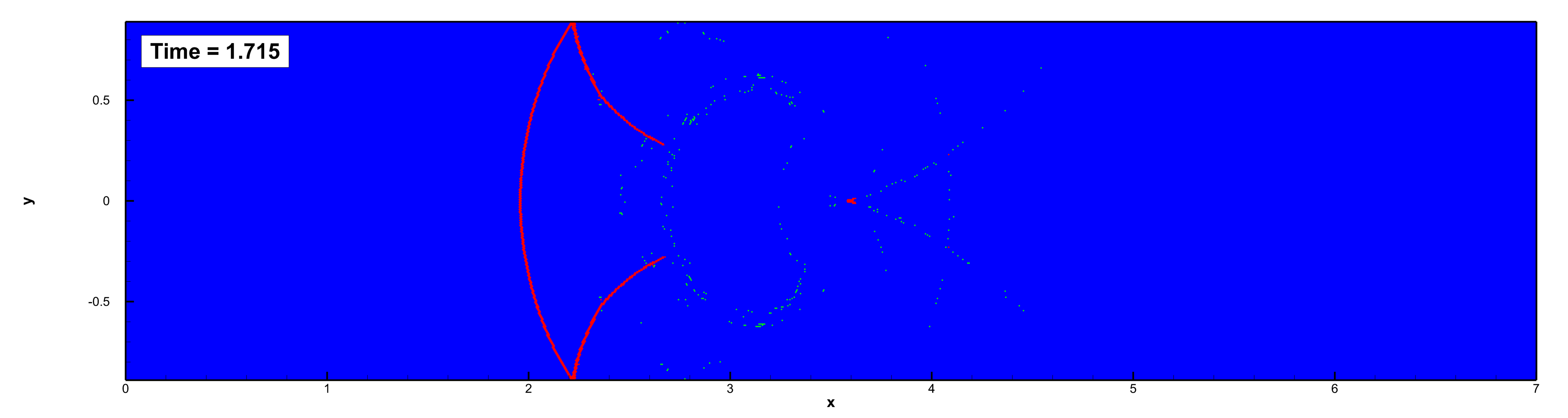}\\[-3pt]
	\includegraphics[width=0.5\linewidth,trim=3 6 3 6,clip]{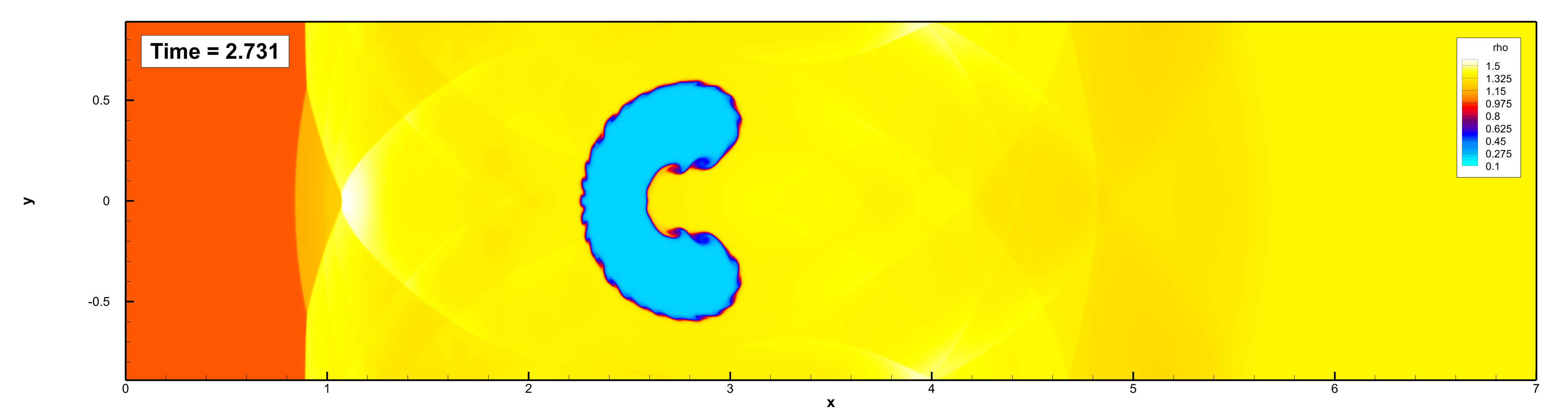}%
	\includegraphics[width=0.5\linewidth,trim=3 6 3 6,clip]{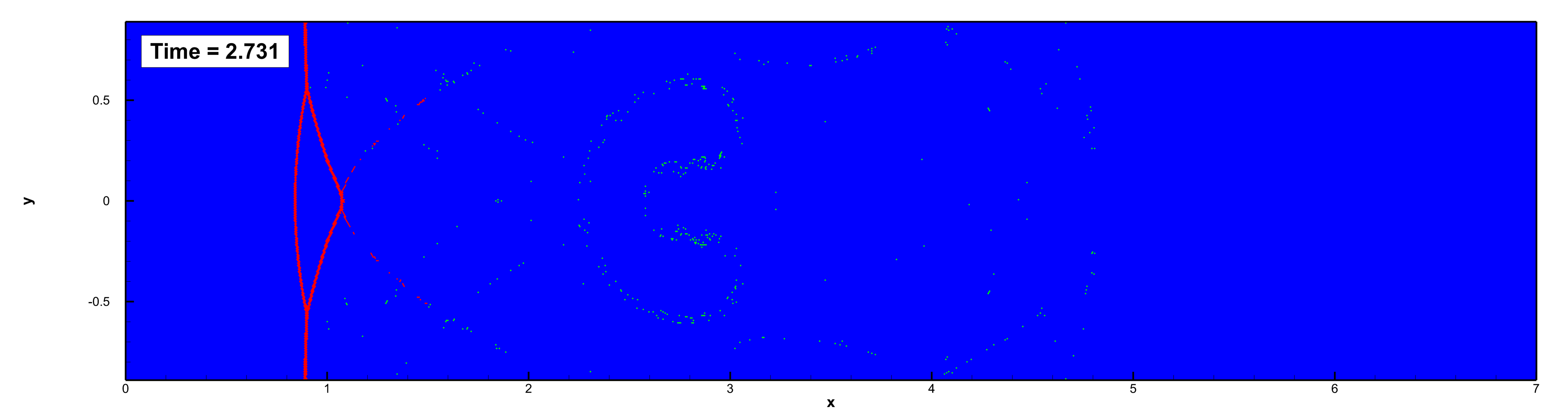}\\[-3pt]
	\includegraphics[width=0.5\linewidth,trim=3 6 3 6,clip]{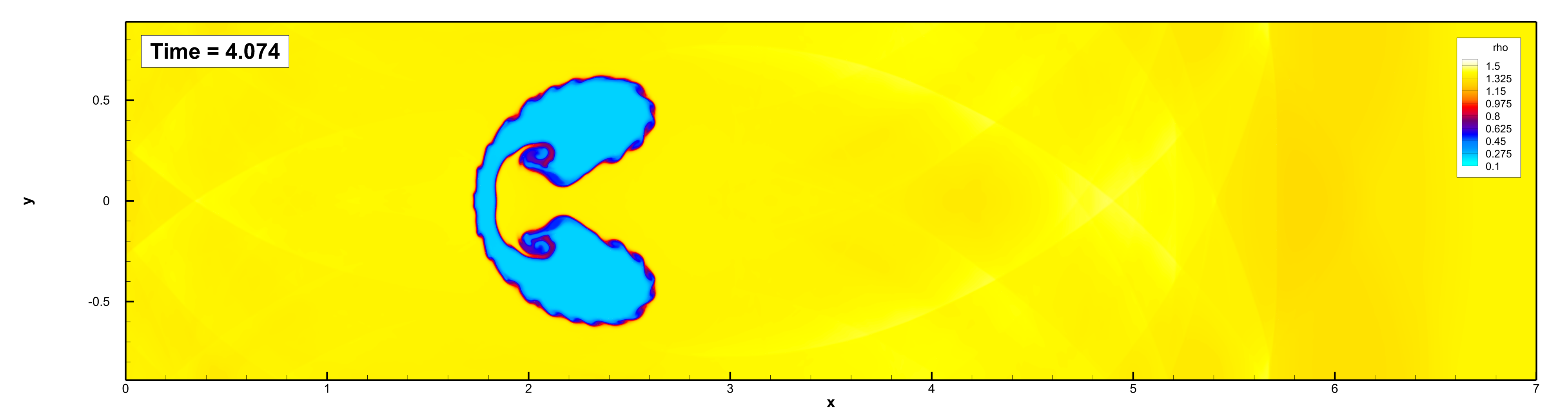}%
	\includegraphics[width=0.5\linewidth,trim=3 6 3 6,clip]{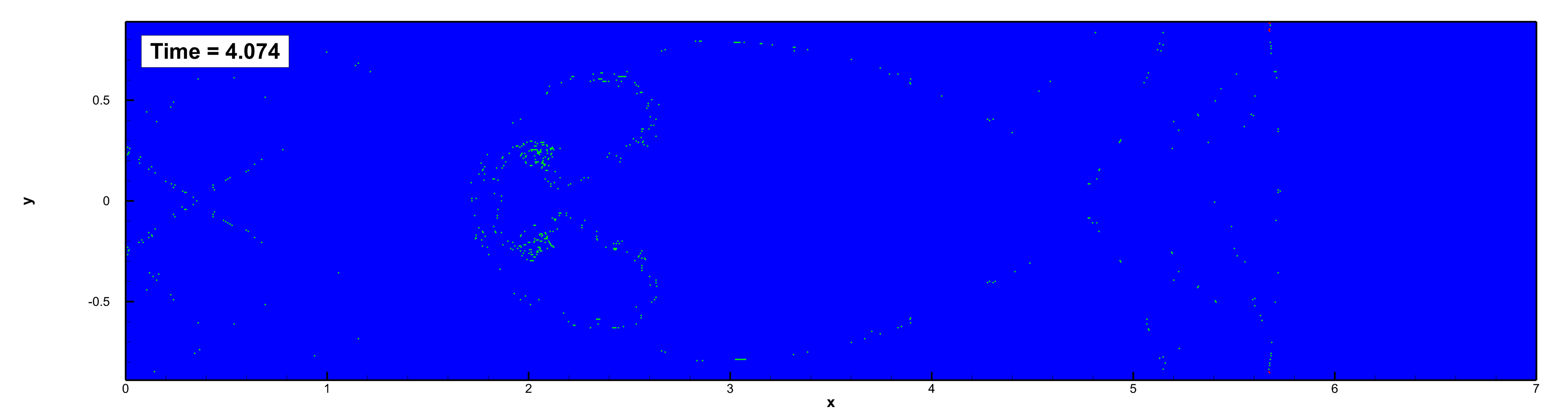}\\[-3pt]
	\includegraphics[width=0.5\linewidth,trim=3 6 3 6,clip]{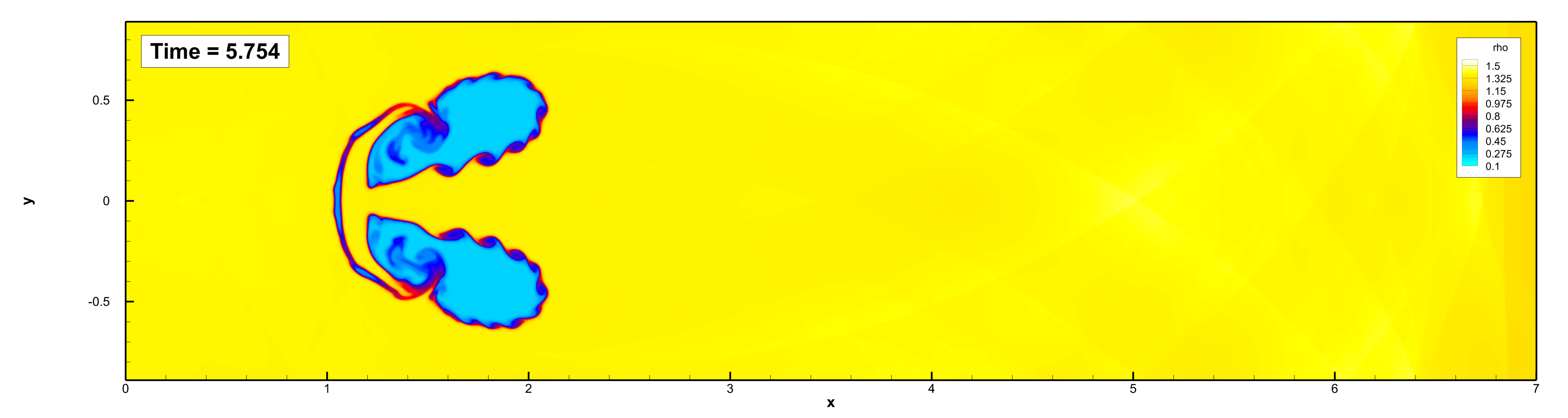}%
	\includegraphics[width=0.5\linewidth,trim=3 6 3 6,clip]{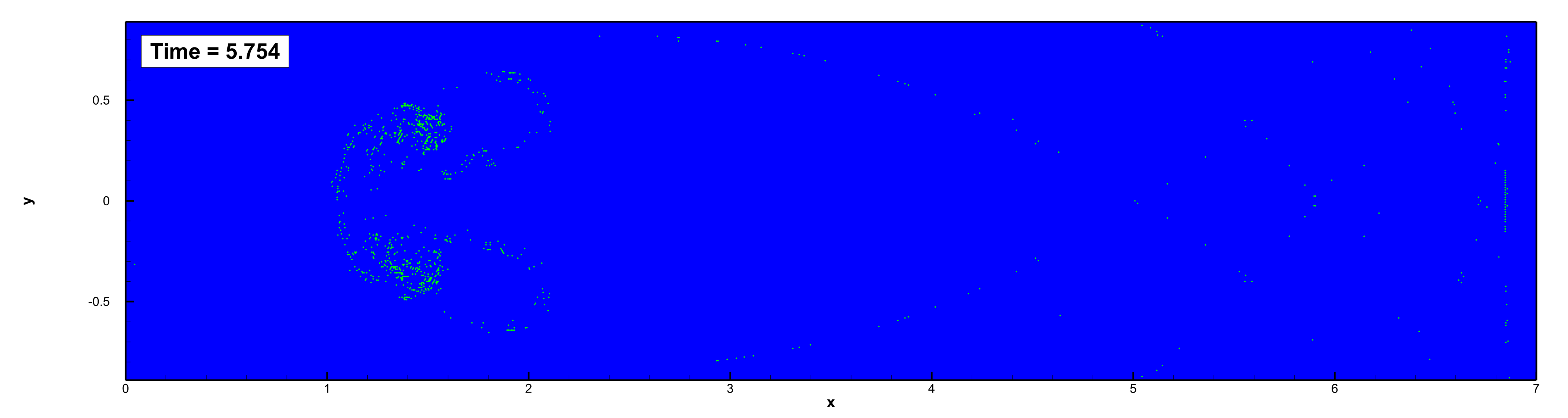}
	\caption{Shock helium-bubble interaction: on the left, we report the density profile obtained with our \Pdue quasi-conservative scheme, 
		and on the right, the distribution of limiter activations and the type of detected troubled cells.
		We emphasize that the limiter activates both at the discontinuity produced by the bubble-air interface and the shock wave;
		however it works with the conservative formulation of the PDE~\eqref{eq.multimaterialCons} only on 
		the \textit{shock-triggered troubled cells} (depicted in red),
		thus without negatively affecting the oscillation-free solution given by the primitive formulation at 
		the material interface.
		The figure shows the initial condition (first line) and the dimensionless time $T^*=0.41$ ($t^*=75$~[$\mu$s]), 
		when the shock first impacts onto the bubble (second line). 
		The other output times correspond to the physical times $t=t^* + [32,62,102,240,427,674,983]$~[$\mu$s].}
	\label{fig.BubbleHelium-Shock_densityandlimiter}
\end{figure}
\begin{figure}[!b]
	\centering
	\includegraphics[width=0.33\linewidth,trim=3 6 3 6,clip]{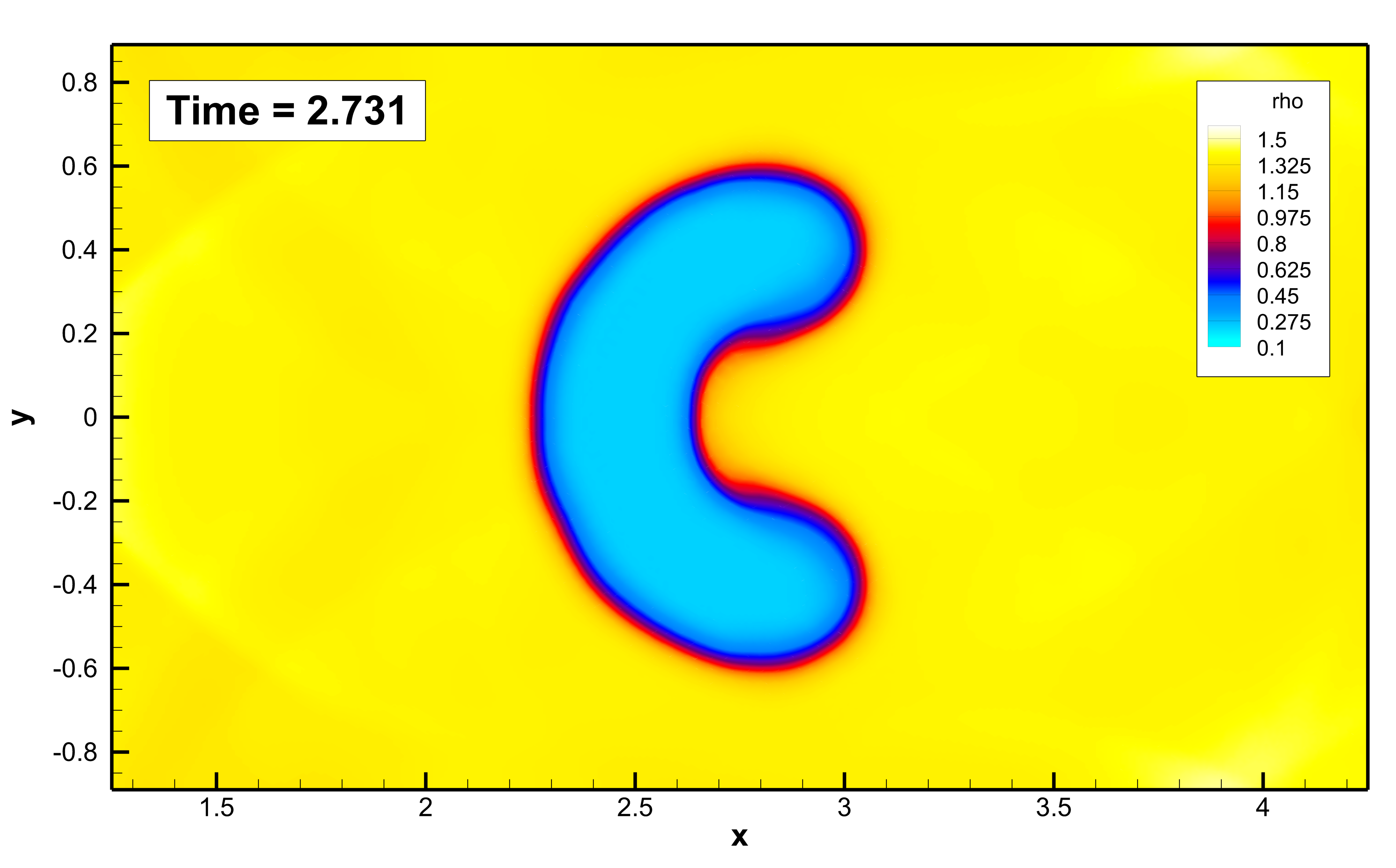}%
	\includegraphics[width=0.33\linewidth,trim=3 6 3 6,clip]{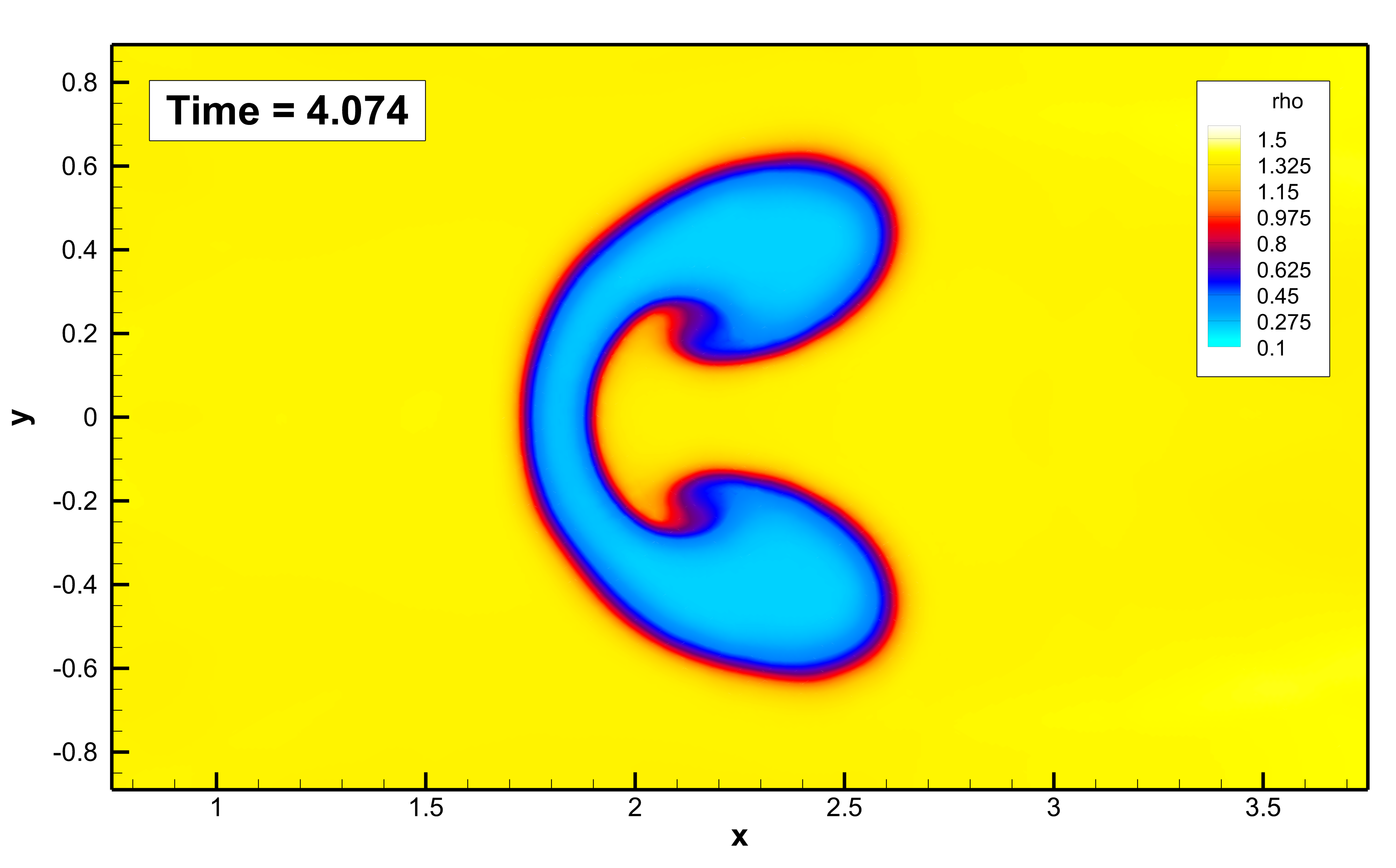}%
	\includegraphics[width=0.33\linewidth,trim=3 6 3 6,clip]{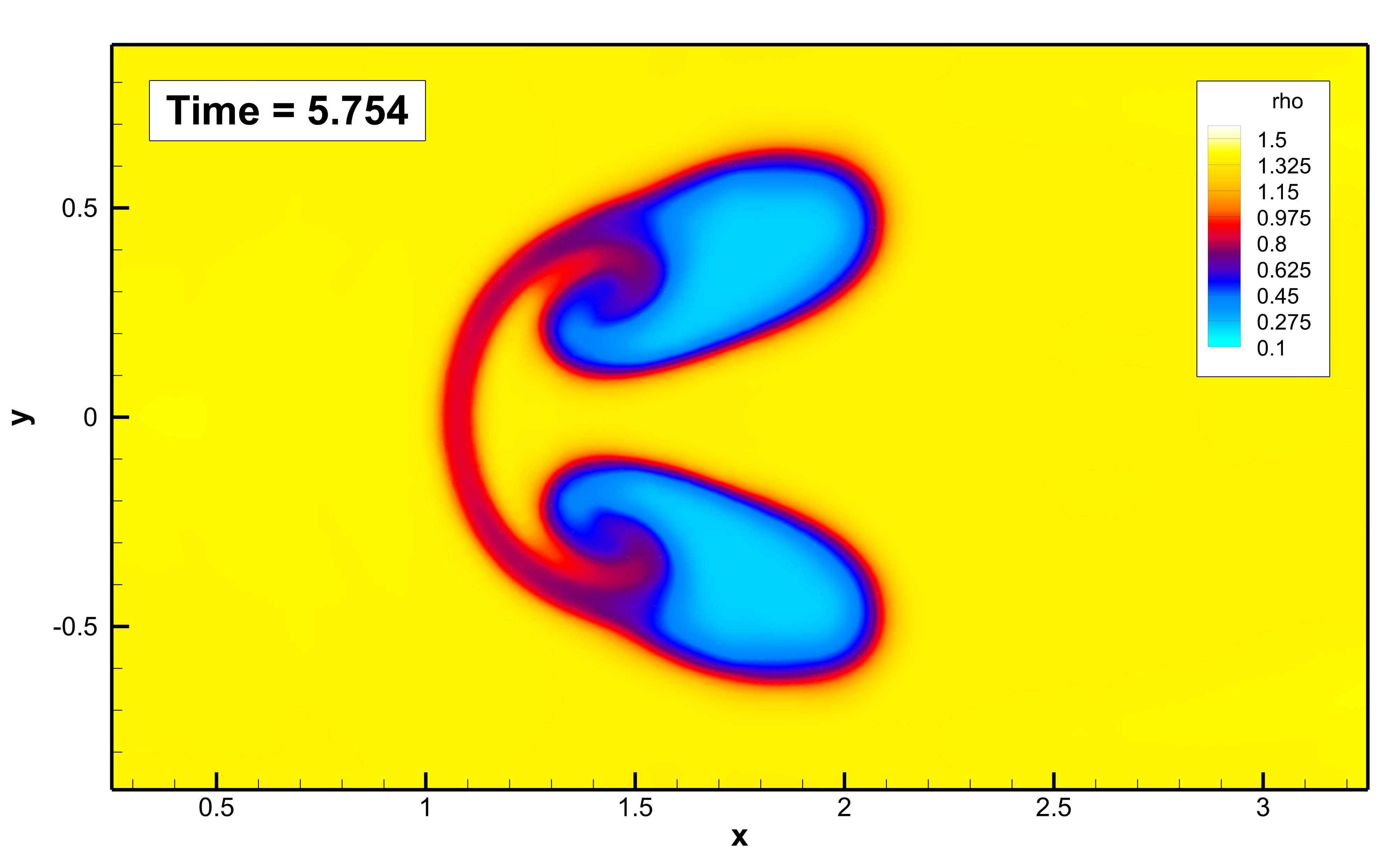}\\[-5pt]
	\includegraphics[width=0.33\linewidth,trim=3 6 3 6,clip]{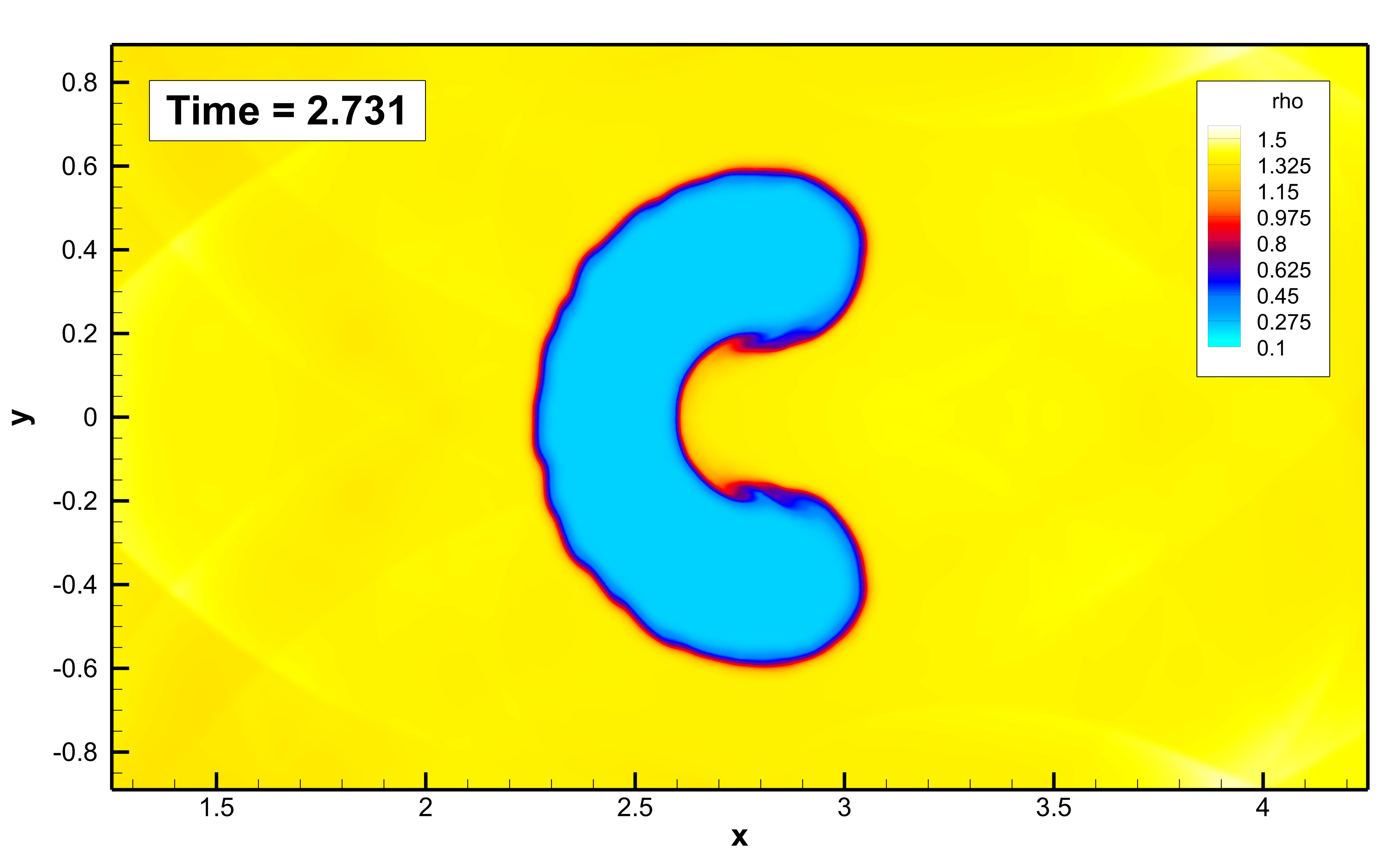}%
	\includegraphics[width=0.33\linewidth,trim=3 6 3 6,clip]{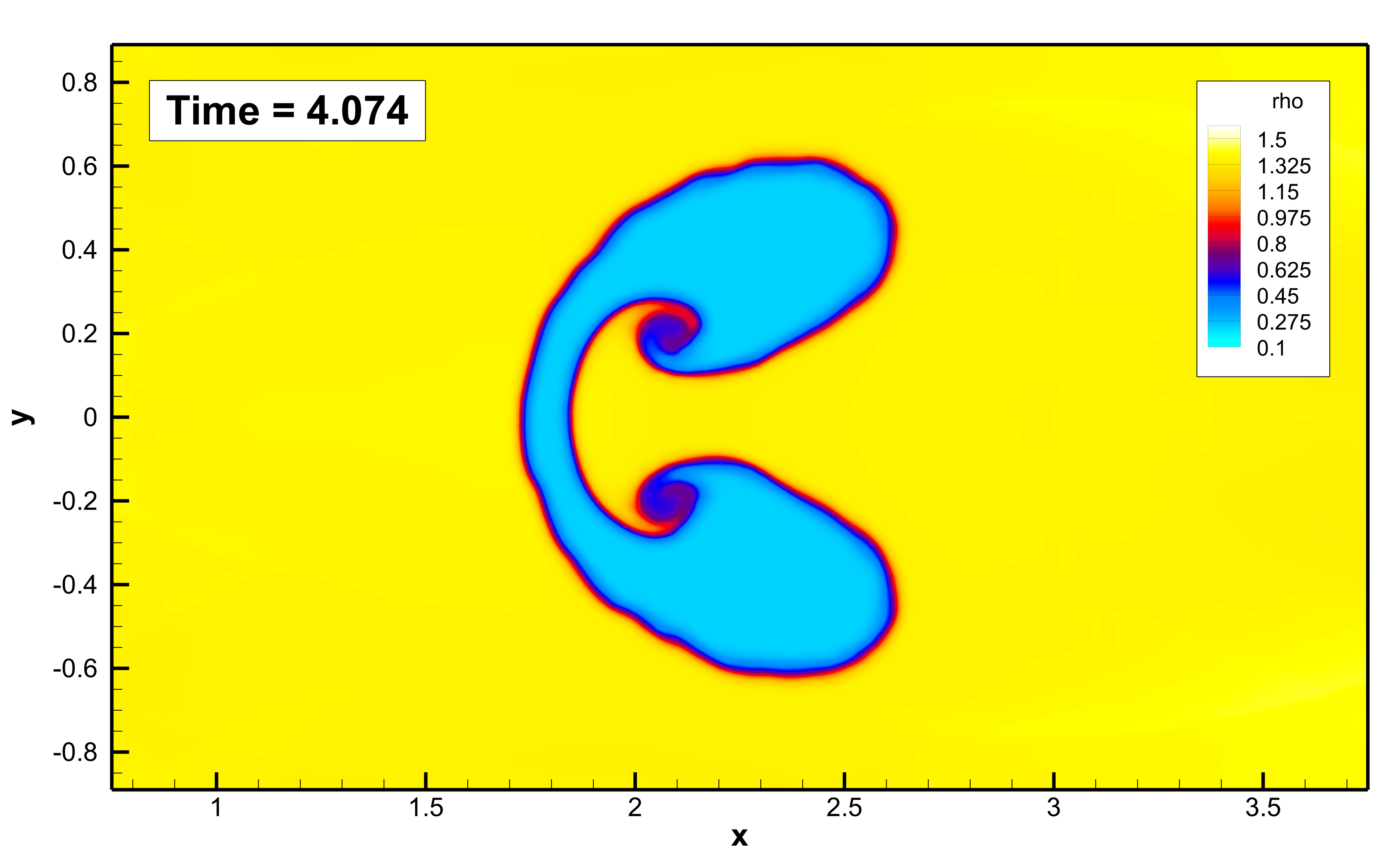}%
	\includegraphics[width=0.33\linewidth,trim=3 6 3 6,clip]{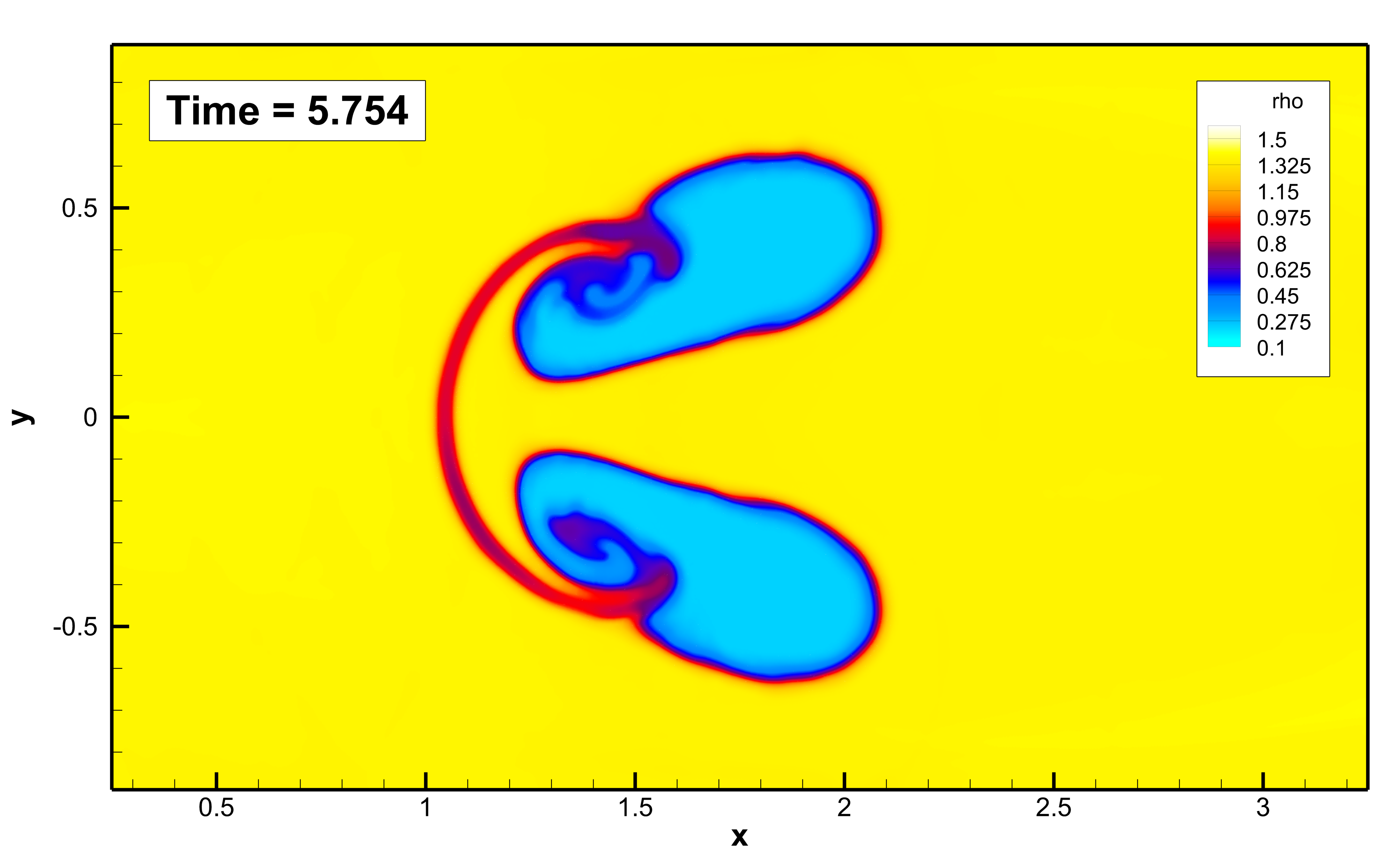}\\[-5pt]
	\includegraphics[width=0.33\linewidth,trim=3 6 3 6,clip]{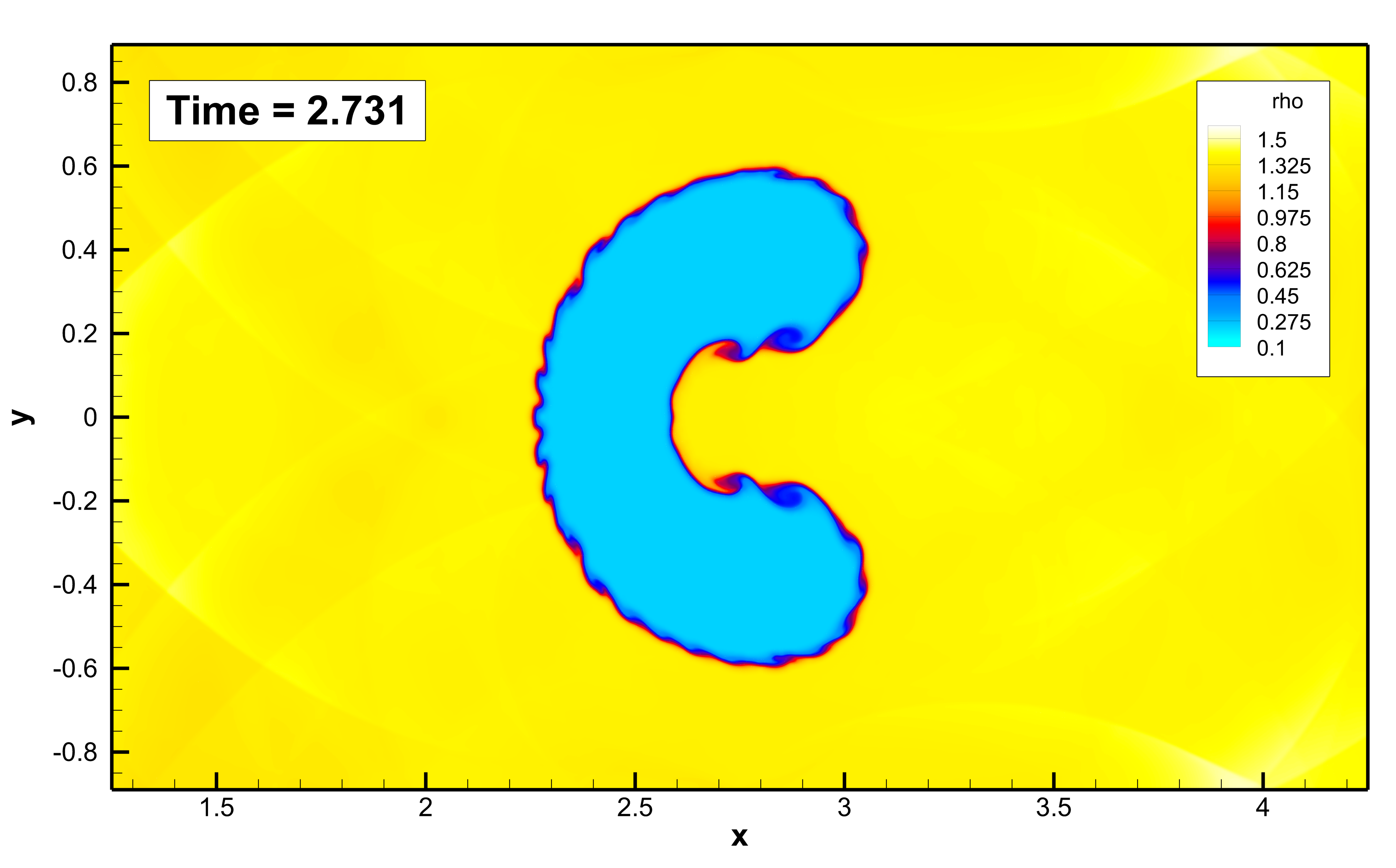}%
	\includegraphics[width=0.33\linewidth,trim=3 6 3 6,clip]{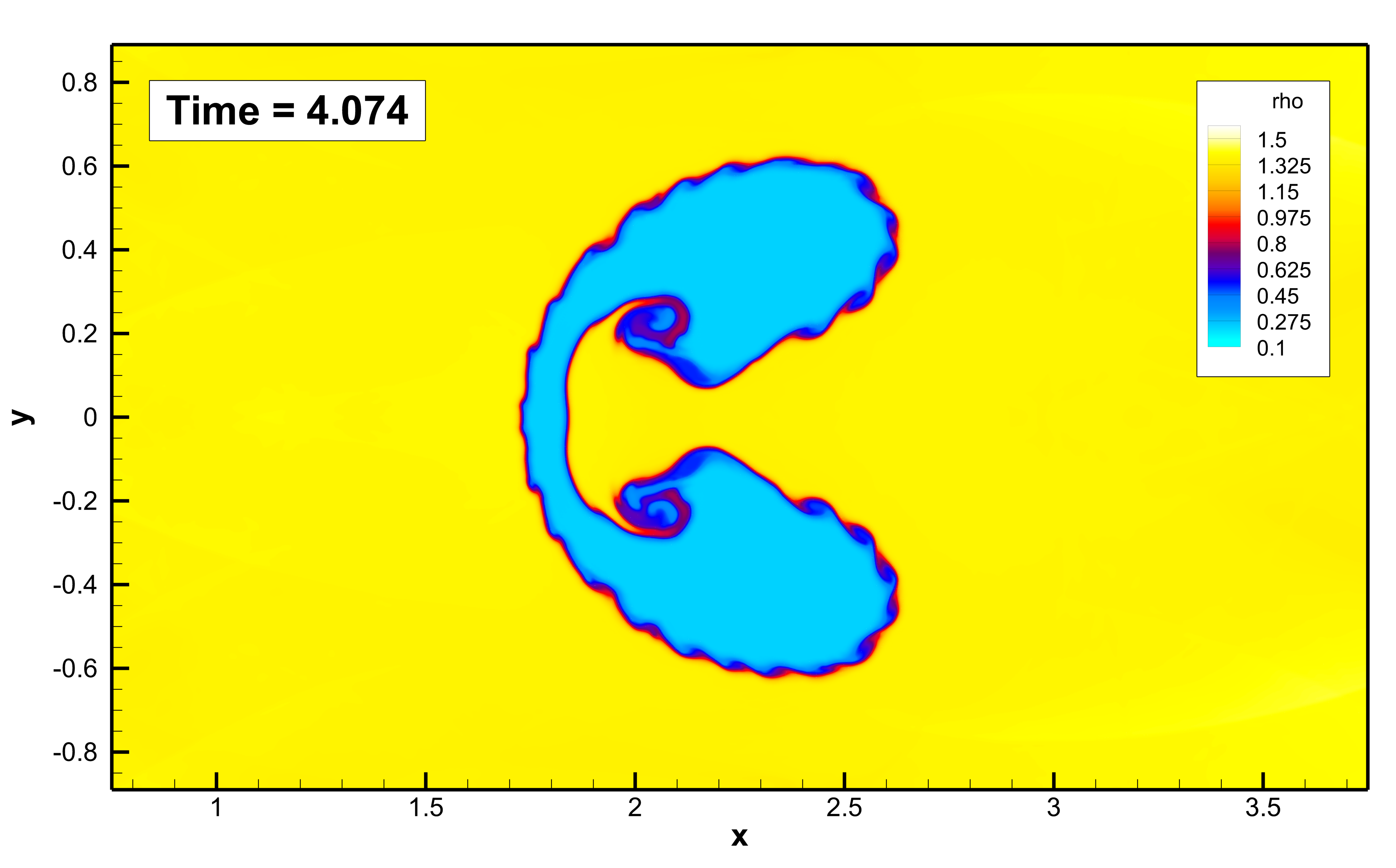}%
	\includegraphics[width=0.33\linewidth,trim=3 6 3 6,clip]{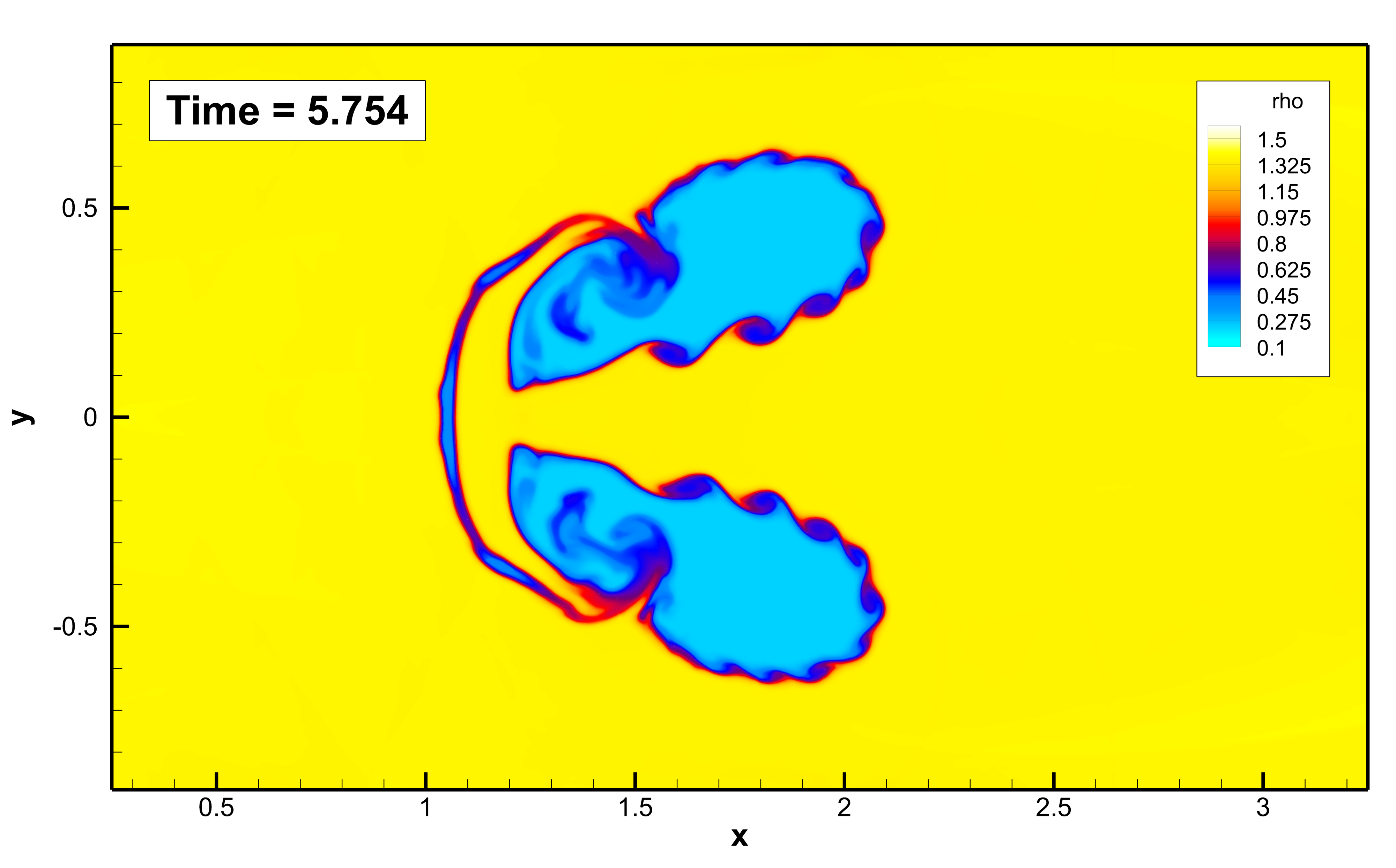}\\[-5pt]
	\includegraphics[width=0.33\linewidth,trim=3 6 3 6,clip]{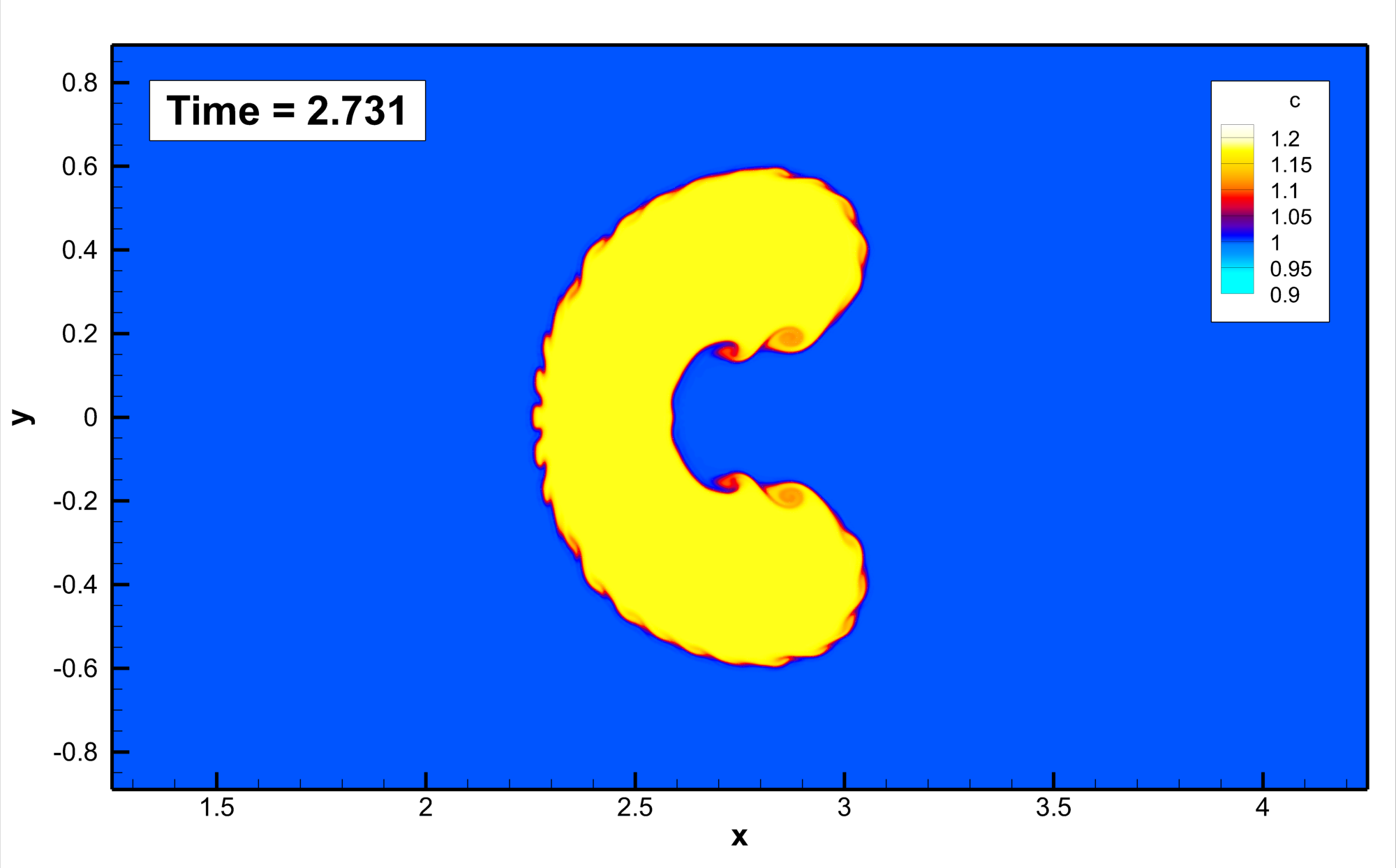}%
	\includegraphics[width=0.33\linewidth,trim=3 6 3 6,clip]{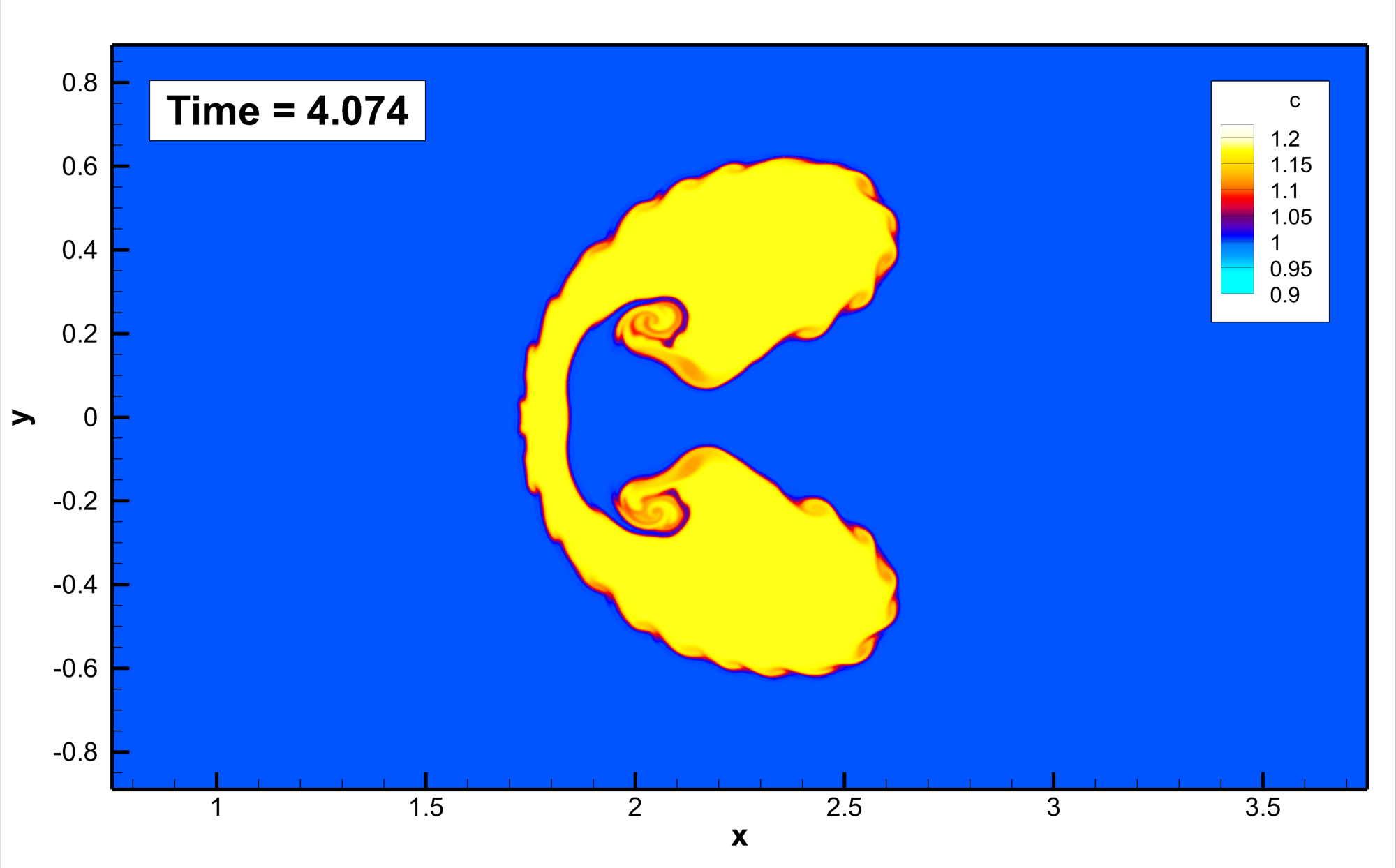}%
	\includegraphics[width=0.33\linewidth,trim=3 6 3 6,clip]{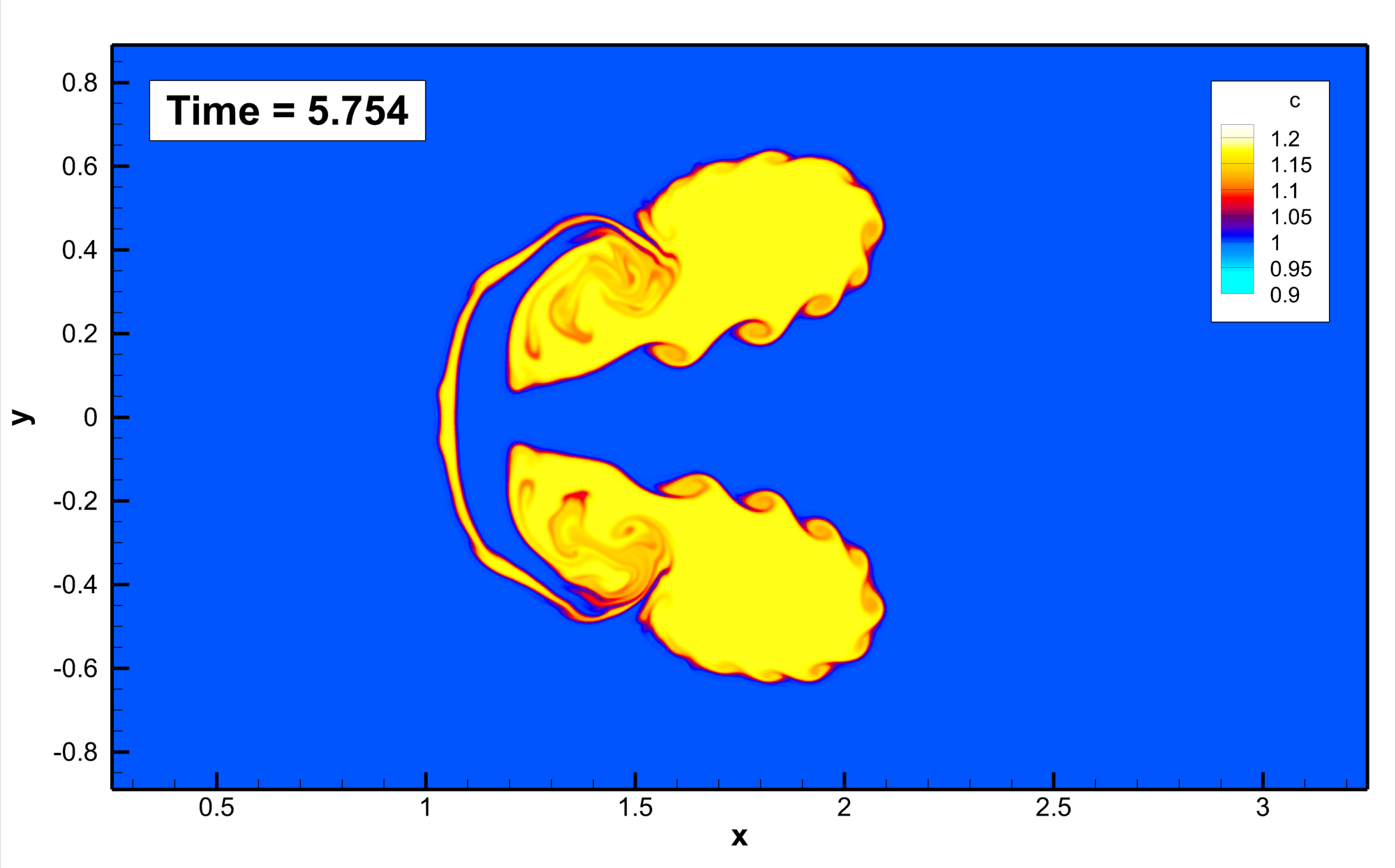}%
	\caption{Mesh convergence of the main flow structures for the shock helium-bubble interaction. In this figure we 
	report the density profile (first three lines) and the values of {$c$} (last line) at three output 
	times $t=2.731$, $t=4.074$, $t=5.754$. In particular, we have used three polygonal tessellations with characteristic mesh size 
	$h=1/30$ (first line), %$h=3.21E-2$
	$h=1/75$ (second line) and %$h=1.30E-2$
	$h=1/150$ (third and fourth line).  %$h=6.49E-3$
	} 	
	\label{fig.BubbleHelium-Shock_meshconvergence}
\end{figure}

Finally, we investigate the two-dimensional dynamics of the shock–bubble interaction. 
We consider a rectangular computational domain of size $[0, 7D]\times[-0.89D, 0.89D]$ filled with air, 
with $D$ being a spatial dimensionless parameter. A shock wave of Mach $M=1.22$ is located at $x=4.5D$ 
and is traveling from right to left, and it will hit a helium bubble of radius $r_b=0.05$ and 
center $(3.5D, \, 0)$. The time coordinate $t$ can be also made non-dimensional using the shock velocity 
and the radius of the bubble, thus obtaining the dimensionless time $T=t/(r_b/c \, M)$, 
where $c_a=\sqrt{\gamma p_a/\rho_a}$ is the sound speed of the air. The pre-shock conditions in the air are 
\begin{equation}
	(\rho, u, v, p, c)(\x) = \left(1, 0, 0, \frac{1}{\gamma}, 1\right), 
\end{equation}
with $\gamma = 1.4$, and the post-shock conditions are
\begin{equation}
	(\rho, u, v, p, c)(\x) = \left(1.3764, -0.3336, 0, \frac{1.5698}{\gamma}, 1\right).
\end{equation}
The initial conditions inside the helium bubble are
\begin{equation}
	(\rho, u, v, p, c)(\x) = \left(0.1819, 0, 0, \frac{1}{\gamma}, \frac{1.648}{\gamma}\right).
\end{equation}
The dimensional quantities adopted for the simulation are computed by considering $D=1$~[m], $c_a=331.6$~[m/s], 
$p_a=101325$~[Pa] and $\rho_a=1.29 \cdot 10^{-3}$~[kg/m$^3$]. 
For what concerns boundary conditions, 
we impose the wall type on the top and bottom side of the domain and Dirichlet on the left and right parts.

This test case corresponds to the experiment of Hass and Sturtevant reported in~\cite{haas1987interaction}. 
It involves shock waves, material interfaces and their interaction
and has been used by several researchers~\cite{karni1992viscous,quirk1996dynamics,abgrall1996prevent,abgrall2003efficient,marquina2003flux,shyue2006wave},
and more recently in~\cite{kundu2021numerical, he2015preventing,henneaux2020extended}, to demonstrate the capability of their
numerical schemes in the simulation of compressible multi-component flows.

The output times are chosen according to the experimental setup in~\cite{haas1987interaction} and the numerical 
results of~\cite{ShockBubble2Phase_IJNMF2011}, hence we consider the series of physical 
times $t=t^* + [32,62,102,240,427,674,983]$~[$\mu$s], which correspond to the dimensionless 
times $T=T^* + [0.584,0.747,0.964,1.714,2.731,4.074,5.754]$. The shock wave hits the bubble at 
time $T^*=0.410$ ($t^*=75$~[$\mu$s]), that is also shown in the results. We report the results obtained 
with our third order \Pdue quasi-conservative DG schemes in Figures~\ref{fig.BubbleHelium-Shock_primvscons},
\ref{fig.BubbleHelium-Shock_densityandlimiter} and~\ref{fig.BubbleHelium-Shock_meshconvergence}. 
In particular, in Figure~\ref{fig.BubbleHelium-Shock_primvscons} we compare our results with those one can obtain 
with a standard conservative second order finite volume scheme on a fine mesh to demonstrate that
when working with the conservative PDE formulation~\eqref{eq.multimaterialCons} the scheme produces spurious pressure oscillations at the material interface, 
while our approach completely avoids these spurious waves.
Then, in Figure~\ref{fig.BubbleHelium-Shock_densityandlimiter} we show the cells 
where the limiter have been activated at various output times, 
and we differentiate between the \textit{shock-triggered troubled cells} (depicted in red),
where we work in conservative variables, by the activations due to the discontinuities of $\rho$ and $c$ at 
the material interface (depicted in green) 
which are treated with our \textit{a posteriori} subcell limiter but still employing the non-conservative formulation~\eqref{eq.multimaterialPrim}.
Finally, in Figure~\ref{fig.BubbleHelium-Shock_meshconvergence} we show the results obtained with successive mesh refinements 
to provide a visual convergence of the main flow structures.

%--------- END OF SECTION -------------------------------------------------	

%--------- SECTION --------------------------------------------------------
\section{Conclusions} \label{sec_concl}	
 
In this paper we have discussed a new approach to efficiently use the non-conservative form of non-linear hyperbolic systems,
while retaining a correct approximation of shocks.  The method is based on a clever combination
of a    non-conservative fully explicit ADER-DG method, with a posteriori  subcell corrections 
which include both a subcell limited update   and  a subcell  correction to remove the conservation defect.
The two corrections are activated by appropriately  defined cell markers, which may flag cells as  either
troubled or shock activated troubled cells (or both).   To provide some theoretical justification to this correction, a  framework
is  proposed to define  the   conservation defect  associated to the fully non-conservative  method.
Conditions on this defect are given for a modified  Lax-Wendroff theorem to hold.
This has allowed to  generalize  to our approach  known  results concerning   the approximation smooth flows, 
contact discontinuities, and material interfaces,
while justifying the idea of removing such defect in shocks.
In practice, the shock detection is performed using a sensor inspired by   Lagrangian  hydrodynamics. 
The  method  has been tested on a variety of problems going from smooth solutions,
to very strong shocks, and shock-contact interactions showing excellent results.

Several extensions and enhancements of this work are   possible.  
The first one is clearly the application  to more complex models, with more  involved  thermodynamics. 
One issue with many such models for multi-material flows is that they often contain non-conservative 
terms. In this case the question arises on what would be the appropriate way to design the conservative correction.
For some systems, enough information exist however to caracterize shocks, of even  build exact Riemann solvers 
(see e.g.~\cite{abgrall_kumar_2014}). This can give some idea.
Another point of enhancement is  the detection of shocks, currently  involving the 
free parameter $m_1$. Some reduction of the sensitivity to this parameter may  come from combining the current
scaled divergence sensor with    smoothness indicator related to the regularity of the approximating polynomials, see e.g.
\cite{KRIVODONOVA2004323}, or using finer calibration based on cells clustering~\cite{doi:10.1137/20M1344081}. 
More involved definitions  may also be explored as   e.g. those used in~\cite{SCIACOVELLI2021105134} for hypersonic flows.

Other extensions are interesting and open technical challenges as e.g. Lagrangian approximations,
to further improve the approximation of contacts, or  implicit explicit time stepping, to improve efficiency 
in presence of stiff  thermodynamics (e.g. large variations of speed of sounds in different materials).

%--------- END OF SECTION -------------------------------------------------	

%%%%%%%%%%%%%%%%%%%%%%%%%%%%%%
\section*{Acknowledgments}
W.~Boscheri and S.~Chiocchetti are members of the INdAM GNCS group in Italy. 
M.~Ricchiuto is member of the CARDAMOM team at the Inria center of the University of Bordeaux. 

E.~Gaburro gratefully acknowledges the support received from the European Union 
with the ERC Starting Grant \textit{ALcHyMiA} (grant agreement No. 101114995).
W.~Boscheri received financial support by Fondazione Cariplo and Fondazione CDP (Italy) under the grant No. 
2022-1895 and by the Italian Ministry of University and Research (MUR) with the PRIN Project 2022 No. 2022N9BM3N.
S.~Chiocchetti acknowledges the support obtained by the Deutsche Forschungsgemeinschaft (DFG) 
via the project \textit{DROPIT}, grant no. GRK 2160/2, and from the European Union’s Horizon Europe 
Research and Innovation Programme under the Marie Skłodowska-Curie Postdoctoral 
Fellowship \textit{MoMeNTUM} (grant agreement No. 101109532).

Views and opinions expressed are however those of the authors only and do not necessarily 
reflect those of the European Union or the European Research Council Executive Agency. 
Neither the European Union nor the granting authority can be held responsible for them.

\appendix

\section{Proof of Proposition \ref{prop.smooth}}%MARIUZSMOOTH proof
\label{app:prop-smooth-proof}

Proposition  \ref{prop.smooth} claims that under the smoothness hypothesis \eqref{eqn.smooth-HP} the 
conservation defect associated to the fully non-conservative formulation  \eqref{eqn.DGscheme}  verifies 
the conditions of Proposition \ref{prop:quasi-cons}.
To demonstrate this estimate, we proceed as follows. First, we note that 
the continuity of the polynomial bases used for the approximation and the mean value theorem allow to write for  
the average of the primitive variables as
$$
\bar \v_i  = \v_h(\x^*_i)
$$
for some point  $\x^*_i\in\omega_i$. {By virtue} of the smoothness assumption \eqref{eqn.smooth-HP} % i like by virtue very much
we thus have in each cell
$$
\v_h(\x)-\bar \v_i  = \mathcal{O}(h).
$$
Using  again \eqref{eqn.smooth-HP} as well as \eqref{eqn.smooth-HP1} we can estimate the terms  in the corrector step as 
$$
\int_{t^n}^{t^{n+1}}  \int_{\omega_i} \phi_k \mathbf{H} d\x dt=  
\mathcal{O} (\Delta t\times h^{2})\;,\quad 
\int_{t^n}^{t^{n+1}}\int_{\partial \omega_i}\phi_k\mathbf{D}\cdot\mathbf{n}dSdt   =  \mathcal{O} (\Delta t\times h\times  h)% +\beta}),
$$
which allows to state that  for all degrees of freedom
$$
\hat \v_{\ell}^{n+1} -\hat \v_{\ell}^{n} \le  \mathcal{O}\left ( 1/\int_{\omega_i} \phi_{\ell}\phi_m d\x \right )    
\; C  \Delta t \, h^{2} =\mathcal{O}(\Delta t) =\mathcal{O}(h),
$$
having used the time step definition \eqref{eq.CFL} in the last equality. 

We now estimate the variation of the local cell averages  combining \eqref{eqn.localConservation3a} and \eqref{eqn.constoprim-diff} as 
$$
\begin{aligned}
\dfrac{\bar{\mathbf{q}}_{i}^{n+1} -\bar{\mathbf{q}}_{i}^{n} }{\Delta t}  =  &\,
\dfrac{1}{|\omega_i| }\dfrac{1}{\Delta t}\left(\,\int\limits_{\omega_i}  \T^{-1}\left (\mathbf{v}_h^{n},\mathbf{v}_h^{n+1} 
\right )\left (\mathbf{v}_h^{n+1} -\mathbf{v}_h^{n}\right )  d\mathbf{x}
\right)\\
= &\,  \dfrac{1}{|\omega_i| }\dfrac{1}{\Delta t} \left(\,\int\limits_{\omega_i}
 \T^{-1}\left ( \tilde\v_h^{n}(\x^*), \tilde\v_h^{n+1}(\x^*) \right ) \left (\tilde \v_h^{n+1}-\tilde  \v_h^{n}\right )d\mathbf{x} 
\right) \\
& + \, \dfrac{1}{|\omega_i| }\dfrac{1}{\Delta t} \left(\,\int\limits_{\omega_i}
\left( \T^{-1}(\mathbf{v}_h^{n},\mathbf{v}_h^{n+1} )-
 \T^{-1}\left( \tilde\v_h^{n}(\x^*), \tilde\v_h^{n+1}(\x^*) \right ) \right) \left (\tilde \v_h^{n+1}-\tilde  \v_h^{n} \right )d\mathbf{x} 
\right).   
\end{aligned}
$$
Further, using  the estimate on $\v_{h}- \bar \v_i$,  the continuity  properties   \eqref{eqn.Jac_prim_cont}
and the definition of $\x_i^*$ we can write 
$$
 \dfrac{\bar{\mathbf{q}}_{i}^{n+1} -\bar{\mathbf{q}}_{i}^{n} }{\Delta t}  = 
(\bar\T^{n+1/2})^{-1}\dfrac{\bar \v_i^{n+1}-\bar \v_i^{n}}{\Delta t} +  C_0 h
$$ 
for some bounded constant $C_0$ depending on the mesh geometry and on the solution derivatives, and  
having set  $(\bar\T^{n+1/2})^{-1}\wb{:=}\T^{-1}(\bar\v_i^{n},\bar\v_i^{n+1} )$. We can now use the update of the non-conservative corrector to see that
\begin{equation}\label{eqn.smooth-proof1}
 \dfrac{\bar{\mathbf{q}}_{i}^{n+1} -\bar{\mathbf{q}}_{i}^{n} }{\Delta t}  = 
- \dfrac{1}{\Delta t} \dfrac{1}{|\omega_i |}\int_{t^n}^{t^{n+1}}\int_{\omega_i} \left (\bar\T^{n+1/2} \right )^{-1} \H(\tilde\v_h)d\x\,dt
- \dfrac{1}{\Delta t} \dfrac{1}{|\omega_i |}\int_{t^n}^{t^{n+1}}\int_{\partial\omega_i} 
\left (\bar\T^{n+1/2}\right)^{-1}\mathbf{D}\cdot\mathbf{n}\, dS dt+  C_0 h.
\end{equation}
We can now manipulate the definition \eqref{eqn.flux-balance} as follows
$$
\begin{aligned}
\Phi_{i}=  & \sum\limits_{{\ell_{i,\alpha}}} \mathbf{F}_{\n_{{\ell_{i}}}}  - 
\dfrac{1}{\Delta t} \int\limits_{t^n}^{t^{n+1}} \int\limits_{\partial \omega_i} \mathbb{F}(\tilde{\mathbf{v}}_h^-)\cdot\mathbf{n}dS
+\dfrac{1}{\Delta t} \int\limits_{t^n}^{t^{n+1}} \int\limits_{\partial \omega_i} \mathbb{F}(\tilde{\mathbf{v}}_h^-)\cdot\mathbf{n}dS\\
= &
\sum\limits_{{\ell_{i}}} \mathbf{F}_{\n_{{\ell_{i}}}}  - 
\dfrac{1}{\Delta t} \int\limits_{t^n}^{t^{n+1}} \int\limits_{\partial \omega_i} \mathbb{F}(\tilde{\mathbf{v}}_h^-)\cdot\mathbf{n}dS
+\dfrac{1}{\Delta t} \int\limits_{t^n}^{t^{n+1}} \int\limits_{ \omega_i} \nabla\cdot \mathbb{F}(\tilde{\mathbf{v}}_h^-)d\x,
\end{aligned}
$$
having assumed exact integration of the boundary terms so that the last term can be exactly expressed as a volume integral
using the Gauss theorem. Consider the integrated numerical flux 
\begin{equation}\label{eqn.cons.flux.proof}
\mathbf{F}_{\mathbf{n}_{i\ell}} = \dfrac{1}{2} \dfrac{1}{\Delta t} \int\limits_{t^n}^{t^{n+1}} \int\limits_{\partial \omega_i}
\left ( \F(\tilde\v_h^+) + \F(\tilde\v_h^-) \right )\cdot \mathbf{n}_{i\ell}
- \dfrac{1}{2} \dfrac{1}{\Delta t} \int\limits_{t^n}^{t^{n+1}} \int\limits_{\partial \omega_i}  \T^{-1}\left ( \dfrac{ 
\tilde\v_h^++\tilde\v_h^-}{2} \right ) \left |\mathsf{D}_{\mathbf{n}_{i\ell}}\right | ( \tilde\v_h^+-\tilde\v_h^-)
\end{equation}
where 
$$
\begin{aligned}
\dfrac{1}{2} \left |\mathsf{D}_{\mathbf{n}_{i\ell}}\right | : = \, & \dfrac{1}{2}\int \limits_0^1 \left|
\left[ \A_{\V}\left(\Path(\vv_h^-,\vv_h^+,s)\right), \ \B_{\V}\left(\Path(\vv_h^-,\vv_h^+,s)\right) \right] \cdot \n \right|  ds \\
= \,  &  \int \limits_0^1 
\dfrac{1}{2} \left[ \A_{\V}\left(\Path(\vv_h^-,\vv_h^+,s)\right), \ \B_{\V}\left(\Path(\vv_h^-,\vv_h^+,s)\right) \right] \cdot \n  \,ds 
-\mathsf{D}_{\mathbf{n}_{i\ell}}^-,
\end{aligned}
$$
having used the relation $(\cdot)^- = [(\cdot)  - |\cdot|]/2$.    Using this definition in the expression for $\Phi_i$  we can write
$$
\begin{aligned}
\Phi_{i{\alpha}}=& \dfrac{1}{2}\dfrac{1}{\Delta t} \int\limits_{t^n}^{t^{n+1}} \int\limits_{\partial \omega_i}\left(\mathbb{F}(\tilde{\mathbf{v}}_h^+)- 
\mathbb{F}(\tilde{\mathbf{v}}_h^-)\right)\cdot\mathbf{n}dS- \dfrac{1}{2}\dfrac{1}{\Delta t} \int\limits_{t^n}^{t^{n+1}}\int\limits_{\partial \omega_i}
\left (\bar\T^{n+1/2} \right )^{-1} \left |\mathsf{D}_{\mathbf{n}_{i\ell}}\right | ( \tilde\v_h^+-\tilde\v_h^-)dS
+\dfrac{1}{\Delta t} \int\limits_{t^n}^{t^{n+1}} \int\limits_{\omega_i} \nabla\cdot \mathbb{F}(\tilde{\mathbf{v}}_h^-)d\x   \\
- &\dfrac{1}{2}\dfrac{1}{\Delta t} \int\limits_{t^n}^{t^{n+1}}\int\limits_{\partial\omega_i}
\bigg( \T ^{-1}(\tilde \v_h^-) -  \left(\bar\T^{n+1/2}\right)^{-1} \bigg) \left |\mathsf{D}_{\mathbf{n}_{i\ell}} \right | ( \tilde\v_h^+-\tilde\v_h^-)dS
-  \dfrac{1}{2}\dfrac{1}{\Delta t} \int\limits_{t^n}^{t^{n+1}}\int\limits_{\partial\omega_i}
\bigg( \T^{-1}\left(\dfrac{ \tilde\v_h^++\tilde\v_h^-}{2} \right )  - \T ^{-1}(\tilde \v_h^-)\bigg)  \left |\mathsf{D}_{\mathbf{n}_{i\ell}} \right | 
( \tilde\v_h^+-\tilde\v_h^-)dS.
\end{aligned}
$$%
The last two terms can be estimated using the  continuity of the Jacobians \eqref{eqn.Jac_prim_cont}, combined with \eqref{eqn.smooth-HP} as 
well as \eqref{eqn.smooth-HP1},
$$
\begin{aligned}
\dfrac{1}{\Delta t}& \int\limits_{t^n}^{t^{n+1}}\int\limits_{\partial \omega_i}
\bigg( \T ^{-1}(\tilde \v_h^-) -  \left(\bar\T^{n+1/2}\right)^{-1} \bigg) \left |\mathsf{D}_{\mathbf{n}_{i\ell}} \right | ( \tilde\v_h^+-\tilde\v_h^-)dS
= \mathcal{O}(|\partial\omega_i|\times h\times h)%^{1+\beta}) 
=\mathcal{O}(h^{3}), \\%+\beta}) \\
\dfrac{1}{\Delta t}&\int\limits_{t^n}^{t^{n+1}}\int\limits_{\partial\omega_i}
\bigg( \T^{-1}\left(\dfrac{ \tilde\v_h^++\tilde\v_h^-}{2}\right)  - \T ^{-1}(\tilde \v_h^-)\bigg)  \left |\mathsf{D}_{\mathbf{n}_{i\ell}} 
\right | ( \tilde\v_h^+-\tilde\v_h^-)dS 
= \mathcal{O}(|\partial\omega_i|\times h\times h)
%h^{1+\beta}\times h^{1+\beta}) 
=\mathcal{O}(h^{3}), %+2\beta})
\end{aligned}
$$
where we still do not account for possible quadrature errors in the boundary integrals.
We now invoke the   continuity of  the flux $\mathbb{F}$, to use formally a mean value linearization of its Jacobian  matrices such that
$$
\left(\mathbb{F}(\tilde{\mathbf{v}}_h^+)- \mathbb{F}(\tilde{\mathbf{v}}_h^-)\right)\cdot\mathbf{n}= 
\bar J_{\mathbf{n}}(\tilde{\mathbf{v}}_h^-,\tilde{\mathbf{v}}_h^+)(  \tilde\v_h^+-\tilde\v_h^- )
= \left\{ \bar \T^{-1}(\tilde{\mathbf{v}}_h^-,\tilde{\mathbf{v}}_h^+)\left[ \bar \A_{\V}(\tilde{\mathbf{v}}_h^-,\tilde{\mathbf{v}}_h^+), 
\bar \B_{\V}(\tilde{\mathbf{v}}_h^-,\tilde{\mathbf{v}}_h^+) \right]\cdot\mathbf{n}\right\} (  \tilde\v_h^+-\tilde\v_h^- ).
$$
We can again use \eqref{eqn.Jac_prim_cont}, combined with \eqref{eqn.smooth-HP}-\eqref{eqn.smooth-HP1}
to obtain
$$
\left(\mathbb{F}(\tilde{\mathbf{v}}_h^+)- \mathbb{F}(\tilde{\mathbf{v}}_h^-)\right)\cdot\mathbf{n}
=  \left\{ (\bar\T^{n+1/2})^{-1}   \int \limits_0^1 
 \left[ \A_{\V}\left(\Path(\vv_h^-,\vv_h^+,s)\right), \ \B_{\V}\left(\Path(\vv_h^-,\vv_h^+,s)\right) \right] \cdot \n  \,ds \right\}(  
 \tilde\v_h^+-\tilde\v_h^- ) + \mathcal{O}(h^{2}). %+2\beta})
$$
Combining all  the previous expressions we can show that  
$$
\begin{aligned}
\Phi_{i} =&
 \dfrac{1}{2}\dfrac{1}{\Delta t} \int\limits_{t^n}^{t^{n+1}} \int\limits_{\partial \omega_i}\left(\mathbb{F}(\tilde{\mathbf{v}}_h^+)- 
 \mathbb{F}(\tilde{\mathbf{v}}_h^-)\right)\cdot\mathbf{n}dS
- \dfrac{1}{2}\dfrac{1}{\Delta t} \int\limits_{t^n}^{t^{n+1}}\int\limits_{\partial \omega_i}
(\bar\T^{n+1/2})^{-1}|\mathsf{D}_{\mathbf{n}_{i\ell}}| ( \tilde\v_h^+-\tilde\v_h^-)dS+\dfrac{1}{\Delta t} \int\limits_{t^n}^{t^{n+1}} 
\int\limits_{ \omega_i} \nabla\cdot \mathbb{F}(\tilde{\mathbf{v}}_h^-)d\x \\
=&  \dfrac{1}{\Delta t} \int\limits_{t^n}^{t^{n+1}} \int\limits_{\partial \omega_i} (\bar\T^{n+1/2})^{-1} \mathbf{D}( \tilde\v_h^+-\tilde\v_h^-)\cdot\mathbf{n}dS +
\dfrac{1}{\Delta t} \int\limits_{t^n}^{t^{n+1}} \int\limits_{ \omega_i} \nabla\cdot \mathbb{F}(\tilde{\mathbf{v}}_h^-)d\x  \\
= & 
  \dfrac{1}{\Delta t} \int\limits_{t^n}^{t^{n+1}} \int\limits_{\partial \omega_i} (\bar\T^{n+1/2})^{-1} \mathbf{D}( \tilde\v_h^+-\tilde\v_h^-)\cdot\mathbf{n}dS +
\dfrac{1}{\Delta t} \int\limits_{t^n}^{t^{n+1}} \int\limits_{ \omega_i}  \T^{-1}(\tilde{\mathbf{v}}_h^-)\H(\tilde{\mathbf{v}}_h^-) d\x 
\\
= &
  \dfrac{1}{\Delta t} \int\limits_{t^n}^{t^{n+1}} \int\limits_{\partial \omega_i} (\bar\T^{n+1/2})^{-1} \mathbf{D}( \tilde\v_h^+-\tilde\v_h^-)\cdot\mathbf{n}dS +
\dfrac{1}{\Delta t} \int\limits_{t^n}^{t^{n+1}} \int\limits_{ \omega_i} (\bar\T^{n+1/2})^{-1} \H(\tilde{\mathbf{v}}_h^-) d\x 
+\mathcal{O}(h^{3}), %+\beta})
\end{aligned}
$$
as long as the boundary quadrature is exact at least for linear polynomials, and 
having used again the regularity of the discrete solution and \eqref{eqn.Jac_prim_cont}
to deduce the estimation   $\| [(\bar\T^{n+1/2})^{-1}- \T^{-1}(\tilde{\mathbf{v}}_h^-)]\H(\tilde{\mathbf{v}}_h^-)\| = \mathcal{O}(h)$.

Comparing this last expression to \eqref{eqn.smooth-proof1}, we readily deduce that 
$$
\Delta_i =  \dfrac{\bar{\mathbf{q}}_{i}^{n+1} -\bar{\mathbf{q}}_{i}^{n} }{\Delta t}   + \dfrac{1}{|\omega_i|}\Phi_{i} = \mathcal{O}(h)
$$
which {concludes the proof}.

\section{Proof of Proposition \ref{prop.contact}}
\label{app:prop-contact-proof}
To  prove the 3 statements  we write the  conservation defect and define a numerical flux that matches the variation of $\mathbf{Q}(\v_h)$ in time up to
appropriately small terms. We start by recalling that due to the hypotheses made on the non-conservative update we have that 
$$
\v_h^{n} = \left(\begin{array}{c}
\rho^n_h\\
u_0 + C_u(x)h^{1+\epsilon}\\
v_0 + C_v(x)h^{1+\epsilon}\\
p_0 + C_p(x)h^{1+\epsilon}\\
\chi_h^n
\end{array}
\right) \;\!,\quad
\tilde\v_h  = \left(\begin{array}{c}
\tilde \rho_h\\
u_0 + C_u(x)h^{1+\epsilon}\\
v_0 + C_v(x)h^{1+\epsilon}\\
p_0 + C_p(x)h^{1+\epsilon}\\
\tilde \chi_h 
\end{array}
\right)  \;\!, \quad
\v_h^{n+1} = \left(\begin{array}{c}
\rho^{n+1}_h\\
u_0 + C_u(x)h^{1+\epsilon}\\
v_0 + C_v(x)h^{1+\epsilon}\\
p_0 + C_p(x)h^{1+\epsilon}\\
\chi_h^{n+1}
\end{array}
\right)\!,
$$
where for simplicity we keep the constants independent of the time iteration (the proof is identical if different values are used).
The second observation is that, due to the previous remarks and to the structure of the Jacobian matrices \eqref{eq.multimaterialPrim} and of their eigenvectors,
at a given interface we have in \eqref{eqn.pathint2}
$$
\mathsf{D}_{\mathbf{n}}^-(\tilde\v_h^+-\tilde\v_h^-) =   \min(0,\pmb{\mathsf{v}}_0\cdot\mathbf{n})\left(\begin{array}{c}
\tilde \rho_h^+-\tilde \rho_h^-\\
0\\
0\\
0\\
\tilde \chi_h^+-\tilde \chi_h^- 
\end{array}
\right) +\mathcal{O}(h^{1+\epsilon}).
$$
Similarly, we can show that 
$$
\mathbf{H}(\tilde\v_h) =   \pmb{\mathsf{v}}_0\cdot   \left(\begin{array}{c}
\nabla \tilde\rho_h \\
0 \\
0\\
0\\
\nabla \tilde\chi_h
\end{array}
\right) +\mathcal{O}(h^{\epsilon}) =
\nabla\cdot (\pmb{\mathsf{v}}_0\tilde \v_h) +\mathcal{O}(h^{\epsilon}).
$$
We now look at the definition of the  conservation defect and evaluate the first term for the mass/momentum/energy equations.
Using the definitions of the conservative variables and the equation of state presented in Section~\ref{sec_pde_Euler} 
\begin{equation}\label{eqn.dQ.contact}
\begin{aligned}
\dfrac{1}{\Delta t}\dfrac{1}{|\omega_i|}\int_{\omega_i} \rho(\v_h^{n+1})d\x -\dfrac{1}{\Delta t}\dfrac{1}{|\omega_i|}\int_{\omega_i} \rho(\v_h^{n})d\x
\, = \, &\dfrac{\bar \rho_i^{n+1}-\bar \rho_i^{n}}{\Delta t} \\
\dfrac{1}{\Delta t}\dfrac{1}{|\omega_i|}\int_{\omega_i} \rho\pmb{\mathsf{v}}(\v_h^{n+1})d\x -\dfrac{1}{\Delta 
t}\dfrac{1}{|\omega_i|}\int_{\omega_i}\rho\pmb{\mathsf{v}}(\v_h^{n})d\x
\, = \, & \pmb{\mathsf{v}}_0\dfrac{\bar \rho_i^{n+1}-\bar \rho_i^{n}}{\Delta t}  +\mathcal{O}(h^{\epsilon}) \\
\dfrac{1}{\Delta t}\dfrac{1}{|\omega_i|}\int_{\omega_i} \rho E(\v_h^{n+1})d\x -\dfrac{1}{\Delta t}\dfrac{1}{|\omega_i|}\int_{\omega_i}\rho E(\v_h^{n})d\x
\, = \, & \kappa_0\dfrac{\bar \rho_i^{n+1}-\bar \rho_i^{n}}{\Delta t}  +  p_0 \dfrac{\bar \lambda_i^{n+1}-\bar \lambda_i^{n}}{\Delta t}  +  
\mathcal{O}(h^{\epsilon}), 
\end{aligned}
\end{equation}
with $\kappa_0 =  \pmb{\mathsf{v}}_0\cdot  \pmb{\mathsf{v}}_0/2$ and $\lambda$ defined in \eqref{eqn.lambda}. The idea is now to use the time 
variations for the cell average density and the marker to define consistent fluxes for each of the three considered cases.

\paragraph{Case 1} In this case, $\chi$ is assumed to be constant in the ball $B_i$, which means that so is $\lambda$ in \eqref{eqn.lambda}.  We  consider 
again  the averaged  numerical flux \eqref{eqn.cons.flux.proof}, {which, under the hypotheses made, can be shown to reduce to}
$$
\sum_{\ell}\mathbf{F}_{\mathbf{n}_{i\ell}} =  \left(\begin{array}{c} 1\\u_0\\v_0\\\kappa_0\\0 \end{array}\right)
\left\{ \dfrac{1}{2\Delta t}\int_{t^n}^{t^{n+1}}\int_{\partial\omega_i} (\tilde\rho_h^++\tilde\rho_h^-)\pmb{\mathsf{v}}_0\cdot\mathbf{n}dSdt -
\dfrac{1}{2\Delta t}\int_{t^n}^{t^{n+1}}\int_{\partial\omega_i} |\pmb{\mathsf{v}}_0\cdot\mathbf{n}| (\tilde\rho_h^+-\tilde\rho_h^-)\,dSdt\right\}
+ \mathcal{O}(h^{2+\epsilon}).
$$ 
Given the hypotheses on the quadrature formulas used, we can write the above equivalently as 
$$
\sum_{\ell}\mathbf{F}_{\mathbf{n}_{i\ell}} =  \left(\begin{array}{c} 1\\u_0\\v_0\\\kappa_0\\0 \end{array}\right)
\left\{ \dfrac{1}{2\Delta t}\int_{t^n}^{t^{n+1}}\int_{\omega_i} \nabla\cdot( \pmb{\mathsf{v}}_0\tilde\rho_h )d\x dt-
\dfrac{1}{2\Delta t}\int_{t^n}^{t^{n+1}}\int_{\partial\omega_i}  \min(0,\pmb{\mathsf{v}}_0\cdot\mathbf{n}) (\tilde\rho_h^+-\tilde\rho_h^-)\,dSdt\right\}
+ \mathcal{O}(h^{2+\epsilon}),
$$ 
having used the identity $(\cdot)^- = [(\cdot)-|(\cdot)|]/2$. Using the estimates on $\mathsf{D}^-_{\mathbf{n}}(\tilde\v_h^+-\tilde\v_h^-)$ and $\mathbf{H}$ 
and the first mode of the density component of the corrector \eqref{eqn.DGscheme},  we deduce
$$
\sum_{\ell}\mathbf{F}_{\mathbf{n}_{i\ell}} = - \left(\begin{array}{c} 1\\u_0\\v_0\\\kappa_0\\0 \end{array}\right) 
|\omega_i| \dfrac{\bar \rho_i^{n+1}-\bar \rho_i^{n}}{\Delta t}+ \mathcal{O}(h^{2+\epsilon}).
$$
Combined with \eqref{eqn.dQ.contact}, this allows to write  for the conservation defect of mass, momentum and energy: 
\begin{equation}\label{eqn.contact.proof1}
\begin{aligned}
(\Delta_i^n)_{\rho}=  &\dfrac{1}{\Delta t}\dfrac{1}{|\omega_i|}\int_{\omega_i} \rho(\v_h^{n+1})d\x -\dfrac{1}{\Delta 
t}\dfrac{1}{|\omega_i|}\int_{\omega_i} \rho(\v_h^{n})d\x
+ \dfrac{1}{|\omega_i|}(\sum_{\ell}\mathbf{F}_{\mathbf{n}_{i\ell}})_{\rho} = \mathcal{O}(h^{\epsilon})\\
(\Delta_i^n)_{\rho\mathsf{v}}= &\dfrac{1}{\Delta t}\dfrac{1}{|\omega_i|}\int_{\omega_i} \rho\pmb{\mathsf{v}}(\v_h^{n+1})d\x 
-\dfrac{1}{\Delta t}\dfrac{1}{|\omega_i|}\int_{\omega_i} \rho\pmb{\mathsf{v}}(\v_h^{n})d\x
+  \dfrac{1}{|\omega_i|}(\sum_{\ell}\mathbf{F}_{\mathbf{n}_{i\ell}})_{\rho\mathsf{v}} = \mathcal{O}(h^{\epsilon})\\
(\Delta_i^n)_{\rho E}=&\dfrac{1}{\Delta t}\dfrac{1}{|\omega_i|}\int_{\omega_i} \rho E(\v_h^{n+1})d\x -\dfrac{1}{\Delta 
t}\dfrac{1}{|\omega_i|}\int_{\omega_i} \rho E(\v_h^{n})d\x
+  \dfrac{1}{|\omega_i|}(\sum_{\ell}\mathbf{F}_{\mathbf{n}_{i\ell}})_{\rho E} = \mathcal{O}(h^{\epsilon})
\end{aligned}
\end{equation}
which is the result sought.
\paragraph{Case 2} We proceed as in the previous case,  and use explicitly  the approximation   \eqref{eqn.eos1a}  to modify the energy increment so that 
\eqref{eqn.dQ.contact} becomes now
\begin{equation}\label{eqn.dQ.contact-a}
\begin{aligned}
\dfrac{1}{\Delta t}\dfrac{1}{|\omega_i|}\int_{\omega_i} \rho(\v_h^{n+1})d\x -\dfrac{1}{\Delta t}\dfrac{1}{|\omega_i|}\int_{\omega_i} \rho(\v_h^{n})d\x
\, = \, &\dfrac{\bar \rho_i^{n+1}-\bar \rho_i^{n}}{\Delta t} \\
\dfrac{1}{\Delta t}\dfrac{1}{|\omega_i|}\int_{\omega_i} \rho\pmb{\mathsf{v}}(\v_h^{n+1})d\x 
-\dfrac{1}{\Delta t}\dfrac{1}{|\omega_i|}\int_{\omega_i}\rho\pmb{\mathsf{v}}(\v_h^{n})d\x
\, = \, & \pmb{\mathsf{v}}_0\dfrac{\bar \rho_i^{n+1}-\bar \rho_i^{n}}{\Delta t}  +\mathcal{O}(h^{\epsilon}) \\
\dfrac{1}{\Delta t}\dfrac{1}{|\omega_i|}\int_{\omega_i} \rho E(\v_h^{n+1})d\x -\dfrac{1}{\Delta t}\dfrac{1}{|\omega_i|}\int_{\omega_i}\rho E(\v_h^{n})d\x
\, = \, & \kappa_0\dfrac{\bar \rho_i^{n+1}-\bar \rho_i^{n}}{\Delta t}  +  \dfrac{p_0}{\gamma 
-1} \dfrac{\bar \chi_i^{n+1}-\bar \chi_i^{n}}{\Delta t}  +  \mathcal{O}(h^{\epsilon}) 
\end{aligned}
\end{equation}
The remainder of the proof is almost identical to the previous case. We consider  again the conservation defect 
associated to the flux \eqref{eqn.cons.flux.proof},
which leads to the same conclusions for both mass and momentum, and notably to the first two relations in \eqref{eqn.contact.proof1}.   
The difference is the extra term in the energy equation for which, {under the current hypotheses},  we can write 
$$
\begin{aligned}
\bigg[\sum_{\ell}\mathbf{F}_{\mathbf{n}_{i\ell}}\bigg]_{\rho E} = \, & \kappa_0 
\left\{ \dfrac{1}{2\Delta t}\int_{t^n}^{t^{n+1}}\int_{\omega_i} \nabla\cdot( \pmb{\mathsf{v}}_0\tilde\rho_h )d\x dt-
\dfrac{1}{2\Delta t}\int_{t^n}^{t^{n+1}}\int_{\partial\omega_i}  \min(0,\pmb{\mathsf{v}}_0\cdot\mathbf{n}) (\tilde\rho_h^+-\tilde\rho_h^-)\,dSdt\right\}\\
& + \dfrac{p_0}{\gamma-1}
\left\{ \dfrac{1}{2\Delta t}\int_{t^n}^{t^{n+1}}\int_{\omega_i} \nabla\cdot( \pmb{\mathsf{v}}_0\tilde\chi_h )d\x dt-
\dfrac{1}{2\Delta t}\int_{t^n}^{t^{n+1}}\int_{\partial\omega_i}  \min(0,\pmb{\mathsf{v}}_0\cdot\mathbf{n}) (\tilde\chi_h^+-\tilde\chi_h^-)\,dSdt\right\}
+ \mathcal{O}(h^{2+\epsilon}).
\end{aligned}
$$ 
Using again the estimates on $\mathsf{D}^-_{\mathbf{n}}(\tilde\v_h^+-\tilde\v_h^-)$ and on $\mathbf{H}$, as well as the first 
mode of both the density and the marker function {corrector}, one easily checks that
$$
\bigg[\sum_{\ell}\mathbf{F}_{\mathbf{n}_{i\ell}}\bigg]_{\rho E}  = - |\omega_i|\left\{\kappa_0\dfrac{\bar \rho_i^{n+1}-\bar \rho_i^{n}}{\Delta t}  + 
 \dfrac{p_0}{\gamma -1} \dfrac{\bar \chi_i^{n+1}-\bar \chi_i^{n}}{\Delta t}
\right\}  +  \mathcal{O}(h^{2+\epsilon}),
$$
which readily allows to prove also the last relation in \eqref{eqn.contact.proof1} for the energy conservation defect.

\paragraph{Case 3} As for the second case, the main issue is now the analysis of the total energy increment, while the validity 
of the first two relations in  \eqref{eqn.contact.proof1} is proved as in Case 1.   
In this case we make use  of hypothesis \eqref{eqn.contact.gammaHP} to pass from the estimates involving $\tilde \lambda$ to the 
estimates involving $\chi/(\gamma -1)$.
More precisely, by  defining
$$
\Gamma_0:= \dfrac{1}{\chi_2-\chi_1}(  \dfrac{1}{\gamma_2-1}  -  \dfrac{1}{\gamma_2-1} ),
$$
 we can write in $\omega_i$ \wb{for some constant $\delta>0$}
$$
\dfrac{\bar\lambda_i^{n+1}-\bar\lambda_i^{n}}{\Delta t} = \Gamma_0 \dfrac{\bar\chi_i^{n+1}-\bar\chi_i^{n}}{\Delta t} + \mathcal{O}(h^{\delta})
$$
and similarly on $\partial \omega_i$
$$
\int_{\partial\omega_i} |\pmb{\mathsf{v}}_0\cdot\mathbf{n}| ( \tilde\lambda(\chi_h^+)  - \tilde\lambda(\chi_h^-)) =\Gamma_0  
 \int_{\partial\omega_i} |\pmb{\mathsf{v}}_0\cdot\mathbf{n}| ( \tilde\chi_h^+  - \tilde\chi_h^-)
+ \mathcal{O}(h^{2+\delta}).
$$
We proceed as in the previous cases and consider the energy component of the flux balanced obtained using  the time integrated 
averaged  numerical flux \eqref{eqn.cons.flux.proof}
which in this case we can show to be 
$$
\begin{aligned}
\bigg[\sum_{\ell}\mathbf{F}_{\mathbf{n}_{i\ell}}\bigg]_{\rho E} =   \kappa_0 &
\left\{ \dfrac{1}{2\Delta t}\int_{t^n}^{t^{n+1}}\int_{\omega_i} \nabla\cdot( \pmb{\mathsf{v}}_0\tilde\rho_h )d\x dt-
\dfrac{1}{2\Delta t}\int_{t^n}^{t^{n+1}}\int_{\partial\omega_i}  \min(0,\pmb{\mathsf{v}}_0\cdot\mathbf{n}) (\tilde\rho_h^+-\tilde\rho_h^-)dSdt\right\}\\
+  p_0&\left\{ \dfrac{1}{2\Delta t}\int_{t^n}^{t^{n+1}}\int_{\partial \omega_i}    \tilde\lambda(\chi_h^-) 
  \pmb{\mathsf{v}}_0\cdot \mathbf{n}dS dt-
\dfrac{1}{2\Delta t}\int_{t^n}^{t^{n+1}}\int_{\partial\omega_i}  \min(0,\pmb{\mathsf{v}}_0\cdot\mathbf{n}) ( \tilde\lambda(\chi_h^+) 
 - \tilde\lambda(\chi_h^-))dSdt\right\}
+ \mathcal{O}(h^{2+\epsilon})\\
=   \kappa_0 &
\left\{ \dfrac{1}{2\Delta t}\int_{t^n}^{t^{n+1}}\int_{\omega_i} \nabla\cdot( \pmb{\mathsf{v}}_0\tilde\rho_h )d\x dt-
\dfrac{1}{2\Delta t}\int_{t^n}^{t^{n+1}}\int_{\partial\omega_i}  \min(0,\pmb{\mathsf{v}}_0\cdot\mathbf{n}) (\tilde\rho_h^+-\tilde\rho_h^-)dSdt\right\}\\
+  \Gamma_0 p_0 &\left\{ \dfrac{1}{2\Delta t}\int_{t^n}^{t^{n+1}}\int_{\partial \omega_i}  \tilde\chi_h^-  \pmb{\mathsf{v}}_0\cdot \mathbf{n}     dS dt-
\dfrac{1}{2\Delta t}\int_{t^n}^{t^{n+1}}\int_{\partial\omega_i}  \min(0,\pmb{\mathsf{v}}_0\cdot\mathbf{n}) (  \tilde\chi_h^+  - \tilde\chi_h^-)dSdt\right\}
+ \mathcal{O}(h^{2+\epsilon}) +\mathcal{O}(h^{2+\delta})  \\
=   \kappa_0 &
\left\{ \dfrac{1}{2\Delta t}\int_{t^n}^{t^{n+1}}\int_{\omega_i} \nabla\cdot( \pmb{\mathsf{v}}_0\tilde\rho_h )d\x dt-
\dfrac{1}{2\Delta t}\int_{t^n}^{t^{n+1}}\int_{\partial\omega_i}  \min(0,\pmb{\mathsf{v}}_0\cdot\mathbf{n}) (\tilde\rho_h^+-\tilde\rho_h^-)dSdt\right\}\\
+  \Gamma_0 p_0 &\left\{ \dfrac{1}{2\Delta t}\int_{t^n}^{t^{n+1}}\int_{ \omega_i}  \nabla\cdot(    \pmb{\mathsf{v}}_0 \tilde\chi_h)     d\x dt-
\dfrac{1}{2\Delta t}\int_{t^n}^{t^{n+1}}\int_{\partial\omega_i}  \min(0,\pmb{\mathsf{v}}_0\cdot\mathbf{n}) (  \tilde\chi_h^+  - \tilde\chi_h^-)dSdt\right\}
+ \mathcal{O}(h^{2+\epsilon}) +\mathcal{O}(h^{2+\delta}),
\end{aligned}
$$ 
having used the fact that contour integrals of a constant times the normal is zero and the previous estimates on 
the jumps of $\tilde \lambda$ on $\partial\omega_i $.
We can now proceed exactly as in the previous case to show that 
$$
(\Delta_i^n)_{\rho E} = \mathcal{O}(h^{\epsilon}) +\mathcal{O}(h^{\delta}),
$$
which proves the final result.

	\bibliographystyle{plain}
	\bibliography{biblio_aderprimitive}

\end{document}